\documentclass[11pt]{article}
\usepackage{amssymb,amsmath,booktabs,hyperref,mathtools}
\usepackage{graphics,authblk}
\usepackage{xcolor}
\usepackage{longtable,nicefrac,diagbox}
\usepackage{graphicx}
\usepackage{algorithm}
\usepackage{algpseudocodex}
\usetikzlibrary{shapes.geometric, arrows.meta, positioning, calc}
\usepackage{pgfplots}
\usepackage{subcaption}
\usepackage{svg}
\usepackage{multirow}
\usepackage{tabularx,enumerate}
\usepackage[inline]{enumitem}
\usepackage{rotating}	
\usepackage{dcolumn}
\usepackage{pdflscape}
\usepackage{array}
\usepackage{lscape}
\usepackage{anysize}
\usepackage{color}
\usepackage{amsthm}
\usepackage{listings}
\newcommand{\mathleft}{\@fleqntrue\@mathmargin0pt}
\newcommand{\mathcenter}{\@fleqnfalse}

\providecommand{\keywords}[1]{\textbf{\textit{Keywords: }} #1}
\usepackage{float}
\usepackage{soulutf8}
\usepackage{adjustbox}
\usepackage{colortbl}
\usepackage{blindtext,nameref}

\newtheorem{theorem}{Theorem}
\newtheorem{proposition}{Proposition}

\theoremstyle{definition}
\newtheorem{definition}{Definition}
\theoremstyle{remark}
\newtheorem*{remark}{Remark}
\allowdisplaybreaks
\begin{document}
	\title{Taxicab Distance Based Best-Worst Method for Multi-Criteria Decision-Making: An Analytical Approach}
	\author{Harshit M. Ratandhara, Mohit Kumar}
	\date{}
	\affil{Department of Basic Sciences,\\ Institute of Infrastructure, Technology, Research And Management, Ahmedabad, Gujarat-380026, India\\ Email: harshitratandhara1999@gmail.com, mohitkumar@iitram.ac.in}
	\maketitle	
	\begin{abstract}
		The best-worst method is a well-known distance based multi-criteria decision-making method used for computing the weights of decision criteria. This article provides a comprehensive analytical examination of the taxicab distance based model of the method, with the objectives of investigating the uniqueness of these solutions, and performing a rigorous consistency analysis. To achieve this, an optimal modification based optimization problem, equivalent to the original one, is first formulated. This reformulated problem is then solved analytically, and the optimal weight sets are derived from its solutions. Contrary to the prevailing understanding derived from numerical experiments with the taxicab model, our analytical framework proves that the model can, in fact, lead to multiple optimal weight sets, and we formally establish the conditions for this occurrence. A mixed-integer linear programming model is then employed to compute the consistency index. A decision-maker-aided selection strategy is also proposed for addressing non-uniqueness of optimal weight sets. In addition, threshold values of the consistency ratio to assess the admissibility of given preferences are also established. This framework provides a solid mathematical foundation that enhances the understanding of the model and eliminates the requirement for optimization software. By significantly improving the model's computational efficiency, it enables implementation in large-scale, dynamic real-world applications such as electricity market bidding strategies and portfolio rebalancing under market volatility. The effectiveness of the proposed framework is demonstrated through numerical examples, and its practical applicability is illustrated via a smartphone selection problem.
	\end{abstract}
	
		
	\section{Introduction}	
	Decision-making is an essential part of daily life. Decision situations involving numerous decision criteria pose significant challenges, particularly when many of these criteria are in conflict. Multi-Criteria Decision-Making (MCDM) is a specialized branch of operations research that assists decision-makers in addressing such complex issues. A fundamental step in resolving an MCDM problem is deriving the weights of decision criteria and determining the priority of alternatives when the values of alternatives concerning a criterion are unknown.\cite{lei2022preference} The methods employed in MCDM for this purpose are known as weighting methods or weight calculation methods. Some of these weighting methods include Analytic Hierarchy Process (AHP),\cite{saaty1994make} Analytic Network Process (ANP),\cite{saaty2004decision} Best-Worst Method (BWM),\cite{rezaei2015best} Simple Multi-Attribute Rating Technique (SMART),\cite{edwards1994smarts} FUCOM (FUll COnsistency Method),\cite{pamuvcar2018new} LBWA (Level Based Weight Assessment) method,\cite{vzivzovic2019new} DIBR (Defining Interrelationships Between Ranked criteria) method,\cite{pamucar2021circular} and the trade-off procedure.\cite{keeney1976decision}\\\\
	The weighting methods require different forms of input from the decision-maker. For instance, SMART requires the decision-maker to directly assign ratings to criteria. In contrast, AHP and BWM collect information in the form of matrix called pairwise comparison matrix $A=(a_{ij})_{n\times n}$, where $n$ denotes the number of criteria, and $a_{ij}$ represents the relative preference of the $i^{th}$ criterion over the $j^{th}$ criterion.\\\\
	The AHP has been one of the most extensively utilized MCDM methods for an extended period, with applications in numerous real-world scenarios.\cite{cebi2023consideration,singh2023assessing} It necessitates pairwise comparisons among each pair of criteria, resulting in a total of $\frac{n(n-1)}{2}$ comparisons. This number increases significantly as $n$ increases. Consequently, for a problem with a large number of criteria, AHP becomes less time-efficient and exhibits greater inconsistency. To overcome this issue, Rezaei\cite{rezaei2015best} developed the BWM, which employs structured comparisons in the form of two vectors, the best-to-other vector $A_{B}=(a_{B1},a_{B2},\ldots,a_{Bn})$ and the other-to-worst vector $A_{W}=(a_{1W},a_{2W},\ldots,a_{nW})^T$, where $B^{th}$ criterion is the best (most preferable) criterion and $W^{th}$ criterion is the worst (least preferable) criterion. Using these comparison values, an optimization problem is formulated, which is nonlinear in nature; hence, the model is referred to as nonlinear BWM. Optimal solutions of this problem yield optimal weights.\\\\
	The BWM requires only $2n-3$ pairwise comparisons, which is a significant reduction compared to the AHP. Although methods such as LBWA, FUCOM, and DIBR require even fewer comparisons ($n-1$), the reduction offered by BWM remains highly meaningful due to its distinct advantages. Unlike these methods, BWM offers a consistency check that identifies contradictory judgments and enables their revision, thereby improving the reliability of the results. bias, while the two-stage evaluation process (best-to-other and other-to-worst) helps counter anchoring effects and promotes more discriminative assessments.\cite{rezaei2021balancing,rezaei2022equalizing,rezaei2024analyzing} These advantages have promoted the application of the BWM in numerous real-world scenarios, such as supplier selection,\cite{ahmadi2017assessing,rostami2023goal} location selection,\cite{kheybari2020sustainable,liang2024inland} supply chain management,\cite{ali2024navigating} energy efficiency,\cite{gupta2017developing,wang2019energy} and healthcare service quality assessment,\cite{khanmohammadi2023healthcare} among others in recent times.\\\\
	Wu et al.\cite{wu2023analytical} and Ratandhara and Kumar\cite{ratandhara2024analytical} introduced an analytical framework for the nonlinear BWM and the multiplicative BWM respectively, providing a mathematical foundation that produces an analytical expression for optimal interval weights. These approaches eliminate the dependency on optimization software, thereby enhancing the efficiency of the model. While these works provide a foundation for other BWM variants, the nonlinear goal programming model of BWM, pioneered by Amiri and Emamat,\cite{amiri2020goal} lacks a comparable analytical foundation. This model determines optimal weights by minimizing the taxicab distance (total deviation) of weight ratios from comparison values, thus also known as the taxicab BWM. A dedicated analysis of this model is crucial because its objective function minimizes the overall total inconsistency, a property not addressed by the nonlinear BWM, which is designed exclusively to minimize the single largest violation.\\\\
	In this work, our aim is to derive the optimal weights by deriving analytical solutions to the underlying optimization problem and perform a consistency analysis. Our approach involves formulating an optimal modification based optimization problem, which yields a collection of consistent PCSs, termed optimally modified PCSs. After establishing a one-to-one correspondence between the collection of optimal weight sets and the collection of optimally modified PCSs, we express each optimally modified PCS in terms of given comparison values and the optimal value of the best-to-worst comparison. We then obtain all possible optimal values for best-to-worst comparison, which leads to all optimally modified PCSs, and subsequently, to all optimal weight sets. The key contributions of this work are as follows:
	\begin{itemize}
		\item We analytically characterize the conditions for non-uniqueness, formally demonstrating that the taxicab BWM may yield multiple optimal weight sets, and propose a decision-maker-aided selection strategy to resolve such situations.
		\item We develop a mixed-integer linear programming model to compute the consistency index, which is required to calculate the consistency ratio—a key indicator used to evaluate irrationality in input data.
		\item We determine threshold values of the consistency ratio, defining the maximum acceptable level of inconsistency in pairwise comparisons.
		\item The proposed analytical framework establishes a rigorous mathematical foundation for the taxicab BWM, eliminating its dependency on optimization software and significantly improving its computational accuracy and time efficiency.
		\item This framework, together with the framework for the nonlinear BWM,\cite{wu2023analytical} enables comparison between the taxicab BWM and nonlinear BWM approaches, clarifying their distinct behaviors and providing guidance on model selection based on specific decision contexts.
	\end{itemize}  
	The remainder of this manuscript is structured as follows: Section 2 presents a review of the relevant literature. Section 3 discusses some preliminaries and provides a brief overview of the taxicab BWM. Section 4 details the analytical framework for taxicab BWM, presenting its weight computation, consistency analysis, threshold determination, numerical validation, solution selection mechanism, comparative evaluation with both traditional optimization software and the nonlinear BWM, and a smartphone selection application. Finally, Section 5 presents concluding remarks and outlines potential directions for future research.
	\section{Literature Review}
	This section reviews the literature on the BWM, organized into three key areas: theoretical advancements, extensions for uncertain environments, and integrations with other MCDM techniques.\\\\
	The literature on the BWM has matured along several distinct thematic pathways. A significant body of work has been dedicated to advancing the method's theoretical core, focusing on its foundational principles and mathematical refinements, as cataloged in Table \ref{lr_theory}. Concurrently, substantial research efforts have been directed towards extending its applicability to scenarios involving imprecise or uncertain decision-maker judgments, with key contributions summarized in Table \ref{lr_fuzzy}. Finally, its role as a powerful component within more complex, integrated decision-making frameworks has been widely explored, a trend detailed in Table \ref{lr_integration}. This structured progression of research, from core theory to practical application under uncertainty and integration, highlights the method's dynamism and provides a comprehensive landscape for situating the current study.
	\begin{longtable}{@{}p{3.5cm}p{8.75cm}@{}}
		\caption{A review of theoretical advancements in BWM}\label{lr_theory}\\
		\toprule
		Model / theoretical concept & Key contribution \\
		\midrule
		\endfirsthead
		\multicolumn{2}{c}{{\it Table \ref{lr_theory}.} ({\it Continued})} \\
		\toprule
		Model / theoretical concept & Key contribution \\
		\midrule
		\endhead
		\bottomrule
		\endfoot
		\bottomrule
		\endlastfoot
		Nonlinear BWM\cite{rezaei2015best}&Established the foundational BWM framework, employed maximum deviation to determine optimal criteria weights and introduced the consistency ratio to evalute the reliability of pairwise comparisons\\
		Interval analysis\cite{rezaei2016best}&Proposed using interval analysis for the nonlinear BWM to determine and rank criteria weights when the model yields multiple optimal solutions\\
		Linear BWM\cite{rezaei2016best}&Developed a BWM model that maintains the nonlinear BWM's philosophy while guaranteeing a unique optimal solution\\
		Euclidean BWM\cite{kocak2018euclidean}&Proposed a BWM variant based on the Euclidean distance, demonstrating its efficiency and developing a computational package for practical application\\
		Group BWM\cite{safarzadeh2018group}&Extended BWM to group decision-making by formulating two mathematical models to aggregate individual preferences and determine optimal criteria weights, demonstrating practical applicability through a real case study\\
		Multiplicative BWM\cite{brunelli2019multiplicative}&Developed a BWM model using a mathematically sound multiplicative metric, incorporating both interval analysis and consistency analysis for this framework\\
		Goal programming-based BWM\cite{amiri2020goal}& Proposed two goal programming models for the BWM—a nonlinear variant (taxicab BWM) using total deviation, and a linear variant—to determine criteria weights\\
		Concentration ratio\cite{rezaei2020concentration}&Introduced the concentration ratio to measure the dispersion of optimal weight intervals from the nonlinear BWM, complementing the consistency ratio to provide enhanced insight into result reliability and flexibility\\
		Consistency measurements and thresholds\cite{liang2020consistency}&Proposed a comprehensive consistency framework featuring an input-based consistency indicator for immediate feedback, an ordinal consistency check, and defined thresholds to determine the reliability of BWM results\\
		Bayesian BWM\cite{mohammadi2020bayesian}& Introduced a Bayesian hierarchical model for group decision-making, deriving aggregated criteria weights and a novel credal ranking with confidence levels to replace simple averaging\\
		Preference rationality analysis\cite{lei2022preference}&Proposed optimization models as interactive advice tools to help decision-makers modify preferences while maintaining both ordinal and cardinal consistency in BWM\\
		Prioritization methods\cite{tu2023priority}&Introduced two new prioritization methods (approximate eigenvalue and logarithmic least squares) for BWM and established their inconsistency thresholds\\
		Analytical framework for nonlinear BWM\cite{wu2023analytical}&Proved solution uniqueness for three criteria and derived closed-form solutions for multiple optimal solutions in larger systems, while introducing a secondary objective for unique solution selection\\
		Nonadditive BWM\cite{liang2023nonadditive}&Proposed a BWM variant using the Choquet integral to model interactions between criteria, introducing both nonlinear and linear optimization models for deriving weights from interdependent pairwise comparisons\\
		Priority derivation methods\cite{xu2024some}&Developed new approaches for BWM priority derivation, extending efficiency concepts and introducing novel inconsistency measures\\
		Parsimonious BWM\cite{corrente2024better}&Extended BWM to handle large alternative sets using reference alternatives and a barycenter approach, significantly reducing required pairwise comparisons while improving preference representation\\
		Analytical framework for multiplicative BWM\cite{ratandhara2024analytical}&Developed closed-form solutions for interval weights and consistency measures, eliminating software dependency and introducing a secondary objective for unique solution selection\\
	\end{longtable}
	\begin{longtable}{@{}p{3.5cm}p{8.75cm}@{}}
		\caption{A review of BWM extensions under uncertainty}\label{lr_fuzzy}\\
		\toprule
		Information type & Key contribution \\
		\midrule
		\endfirsthead
		\multicolumn{2}{c}{{\it Table \ref{lr_fuzzy}.} ({\it Continued})} \\
		\toprule
		Information type & Key contribution \\
		\midrule
		\endhead
		\bottomrule
		\endfoot
		\bottomrule
		\endlastfoot
		Intuitionistic fuzzy number\cite{mou2016intuitionistic} & Developed an intuitionistic fuzzy multiplicative BWM with new aggregation operators and consistency measures for group decision-making, applied to pulmonary emphysema severity assessment\\
		Fuzzy number\cite{guo2017fuzzy} & Proposed a fuzzy BWM using triangular fuzzy numbers and graded mean integration representation for criteria weight calculation, introducing a new consistency ratio and validating with three case studies\\
		Hesitant fuzzy number\cite{ali2019hesitant} & Developed a BWM framework utilizing score functions and nonlinear optimization for deriving criteria weights, established a new consistency ratio, and validated through three case studies\\
		Fuzzy number\cite{dong2021fuzzy} & Developed a BWM framework using fuzzy equations and multiple linear programming models for different risk preferences, introducing new consistency indices with application examples\\
		Fuzzy number\cite{mohtashami2021novel} & Proposed a novel fuzzy BWM that directly derives crisp weights from fuzzy pairwise comparisons, eliminating the need for fuzzy aggregation and ranking procedures\\
		Intuitionistic fuzzy number\cite{wan2021novel} & Developed a novel fuzzy BWM model using intuitionistic fuzzy equations and multiple linear programming models for different risk preferences, incorporating a consistency improvement process\\
		Generalized interval-valued trapezoidal fuzzy number\cite{wan2021fuzzy} & Introduced a novel fuzzy BWM framework with new normalized weight vectors and consistency measures, demonstrating effectiveness through three real-world applications\\
		Z$^\text{E}$-number\cite{haseli2023sustainable} & Developed a novel Z$^\text{E}$-number-based BWM combined with Combined Compromise Solution (CoCoSo) method for sustainable and resilient recycling partner selection, implementing it in a real-world urban case study\\
		Interval-valued intuitionistic fuzzy number\cite{dong2024interval} & Developed an interval-valued intuitionistic fuzzy BWM with additive consistency using graph theory and linear goal programming, introducing new consistency indices and validation through three applications\\
		Fuzzy number\cite{ratandhara2024alpha} & Proposed an $\alpha$-cut based Fuzzy BWM using exact fuzzy arithmetic to optimize comparison value shapes, introducing approximation error analysis and new consistency measures for supply chain risk assessment\\
		Spherical fuzzy number\cite{haseli2024extension} & Proposed a spherical fuzzy BWM using a nonlinear optimization model to determine weight coefficients, demonstrating a threefold improvement in consistency ratio compared to existing methods\\
		Z$^\text{E}$-number\cite{haseli2025logistic} & Developed a Z$^\text{E}$-number-based BWM integrated with Multi-Attribute Border Approximation Area Comparison (MABAC) for logistics hub location selection under uncertainty in a group decision-making context\\
	\end{longtable}
	 
	\begin{table}[t!]
		\caption{A review of hybrid models integrating BWM with other MCDM methods}\label{lr_integration}
		\begin{tabular}{@{}p{5.5cm}p{6.75cm}@{}}
			\toprule
			Integrated method(s) & Key contribution\\
			\midrule
			Tradeoff procedure\cite{liang2022best} & Unified the axiomatic value-based foundation of the tradeoff procedure with the structured, consistency-checked comparisons of BWM, explicitly incorporating attribute ranges into the weighting process through a novel optimization model \\
			Technique for Order of Preference by Similarity to Ideal Solution (TOPSIS)\cite{varchandi2024integrated} &Balanced theoretical rigor with practical simplicity in sustainable and resilient supplier selection, reducing comparison burden while maintaining ranking consistency and accuracy for e-commerce applications \\
			VlseKriterijumska Optimizacija I Kompromisno Resenje (VIKOR)\cite{gao2024novel} & Developed a novel Fermatean fuzzy entropy measure and hybrid weighting scheme to evaluate medical waste treatment technologies, providing a robust framework for handling high uncertainty in critical environmental health decisions \\
			Multiplicative Multi-Objective Optimization by Ratio Analysis (MULTIMOORA)\cite{yucesan2024evaluating} & Enhanced the UN's Sustainable Urban Transport Index with a robust weighting and ranking framework using interval type-2 fuzzy sets, enabling sensitivity analysis for urban mobility planning \\
			ELimination Et Choix Traduisant la REalité (ELECTRE)\cite{chen2024selecting} & Introduced a Bayesian BWM to robustly determine criteria weights under uncertainty, integrated with ELECTRE III to identify optimal honeycomb structural materials for clean rooms through outranking comparisons\\
			Method based on the Removal Effects of Criteria (MEREC) and TOPSIS\cite{kousar2025multi} & Addressed smog mitigation policy by objectively selecting key criteria using MEREC and determining their reliable weights through an entropy-based BWM, enhancing decision robustness \\
			Decision-Making Trial and Evaluation Laboratory (DEMATEL)\cite{kolour2026enhancing} & Integrated Fuzzy DEMATEL for causal criteria analysis with BWM for precise weight allocation, creating a robust decision framework for public sector supplier selection \\
			\bottomrule
		\end{tabular}
	\end{table}
	 
	\section{Basic Concepts and Brief Introduction to Taxicab Best-Worst Method}
	In this section, we first discuss some foundational definitions and results relevant to our study. We then briefly introduce the taxicab distance based BWM, an equivalent formulation to the nonlinear goal programming model for BWM proposed by Amiri and Emamat.\cite{amiri2020goal}
	\subsection{Preliminaries}
	The following definitions and results are essential for the development of an analytical framework for the taxicab BWM.
	\begin{definition}\cite{vincze2019average}
		Let $x=(x_1,x_2,\ldots,x_n)$, $y=(y_1,y_2,\ldots,y_n)$ be elements of $\mathbb{R}^n$. Then the function $d:\mathbb{R}^n\times \mathbb{R}^n\rightarrow \mathbb{R}_{\geq 0}$ defined by $$d(x,y)=\sum_{i=1}^{n} |x_i-y_i|$$ is called the taxicab distance function on $\mathbb{R}^n$.
	\end{definition}
	\textbf{Notations}: Throughout the article, $C=\{c_1,c_2,\ldots,c_n\}$ denotes the set of criteria and $D=\{c_1,c_2,\ldots,c_n\}\setminus\{c_{B},c_{W}\}$ denotes the set of criteria other than the best and worst ones. Whenever there is no ambiguity, these sets are simply referred to as the sets of indices, i.e., $C=\{1,2,\ldots,n\}$ and $D=\{1,2,\ldots,n\}\setminus\{B,W\}$.
	\begin{definition}\cite{rezaei2015best}
		A Pairwise Comparison System (PCS) $(A_{B},A_{W})$, where $A_{B}$ and $A_{W}$ are the best-to-other and the other-to-worst vector respectively, is said to be consistent if $a_{Bi}\times a_{iW}=a_{BW}$ for all $i\in D$.
	\end{definition}
	\begin{theorem}\label{accurate}\cite{wu2023analytical}
		The system of equations
		\begin{equation}\label{system}
			\frac{w_{B}}{w_i}=a_{Bi},\quad \frac{w_i}{w_{W}}=a_{iW},\quad \frac{w_{B}}{w_{W}}=a_{BW},\ i\in D
		\end{equation}
		has a solution if and only if $(A_{B},A_{W})$ is consistent. Moreover, if solution exists, then it is unique and is given by
		\begin{equation}\label{accurate_weights}
			w_j=\frac{a_{jW}}{\displaystyle\sum_{i\in C}a_{iW}}=\frac{1}{a_{Bj}\displaystyle\sum_{i\in C}\frac{1}{a_{Bi}}},\ j\in C.
		\end{equation}
	\end{theorem}
	Theorem \ref{accurate} assigns a unique weight set to each consistent PCS.
	\begin{definition}\cite{liang2020consistency}
		A PCS $(A_{B},A_{W})$ is said to be ordinal-consistent if $$(a_{Bi}-a_{Bj})\times(a_{jW}-a_{iW})>0 \text{ or } (a_{Bi}=a_{Bj}\ \&\ a_{jW}=a_{iW})$$ for all $i,j\in C$.
	\end{definition}
	\subsection{Taxicab BWM}
	In the taxicab BWM, optimal weights are those that minimize the taxicab distance, i.e., the Total Deviation (TD), of weight ratios from the comparison values. For a given PCS $(A_{B}, A_{W})$, optimal weights are computed by solving the following minimization problem.
	\begin{align}
		&\min \text{ TD=}\sum_{i\in D} \bigg(\bigg|\frac{w_{B}}{w_i}-a_{Bi}\bigg|+ \bigg|\frac{w_i}{w_{W}}-a_{iW}\bigg|\bigg)+ \bigg|\frac{w_{B}}{w_{W}}-a_{BW}\bigg|\nonumber\\
		&\text{sub to: \quad}w_1+w_2+\ldots+w_n=1,\nonumber\\
		\label{optimization_1}
		&\quad \quad \quad\ \quad   w_j\geq0 \text{ for all }j\in C.
	\end{align}
	Problem \eqref{optimization_1} is a nonlinear problem with $n$ variables $w_1, w_2, \ldots, w_n$. So, it has optimal solution(s) of the form $(w_1^*, w_2^*, \ldots, w_n^*)$. Each optimal solution gives an optimal weight set $W^*=\{w_1^*,w_2^*,\ldots,w_n^*\}$, and the optimal objective value is the minimum possible TD of weight ratios from the comparison values. Now, consider the following minimization problem.
	\begin{align}
		&\min \epsilon=\sum_{i\in D} (\epsilon_{Bi}+\epsilon_{iW})+ \epsilon_{BW}\nonumber\\
		&\text{sub to: \quad} \bigg|\frac{w_{B}}{w_i}-a_{Bi}\bigg|=\epsilon_{Bi},\quad \bigg|\frac{w_i}{w_{W}}-a_{iW}\bigg|=\epsilon_{iW},\quad \bigg|\frac{w_{B}}{w_{W}}-a_{BW}\bigg| =\epsilon_{BW},\nonumber\\
		&\quad \quad \quad\ \quad   w_1+w_2+\ldots+w_n=1,\nonumber\\
		\label{optimization_2}
		&\quad \quad \quad\ \quad   w_j\geq0 \text{ for all }j\in C.   
	\end{align}
	Problem \eqref{optimization_2} is an equivalent formulation of problem \eqref{optimization_1}. It has optimal solution(s) of the form $(w_j^*,\epsilon_{Bi}^*,\epsilon_{iW}^*,\epsilon_{BW}^*)$, where $i\in D$ and $j\in C$, with the optimal objective value $\epsilon^*$. For each optimal solution, $w_j^*$ forms an optimal weight set. The value $\epsilon^*$ represents the optimal TD of weight ratios from the comparison values. Now, consider the nonlinear goal programming model for BWM developed by Amiri and Emamat.\cite{amiri2020goal}
	\begin{align}
		&\min \sum_{i\in D} (\epsilon_{Bi}^++\epsilon_{Bi}^-+\epsilon_{iW}^++\epsilon_{iW}^-)+ \epsilon_{BW}^++\epsilon_{BW}^-\nonumber\\
		&\text{sub to: \quad} \frac{w_{B}}{w_i}-a_{Bi}=\epsilon_{Bi}^+-\epsilon_{Bi}^-,\nonumber\\
		&\quad \quad \quad\ \quad\frac{w_i}{w_{W}}-a_{iW}=\epsilon_{iW}^+-\epsilon_{iW}^-,\nonumber\\
		&\quad \quad \quad\ \quad\frac{w_{B}}{w_{W}}-a_{BW} =\epsilon_{BW}^+-\epsilon_{BW}^-,\nonumber\\
		&\quad \quad \quad\ \quad   w_1+w_2+\ldots+w_n=1,\nonumber\\
		\label{optimization_3}
		&\quad \quad \quad\ \quad   \epsilon_{Bi}^+,\epsilon_{Bi}^-,\epsilon_{iW}^+,\epsilon_{iW}^-,\epsilon_{BW}^+,\epsilon_{BW}^-,w_j\geq0 \text{ for all }i\in D\text{ and }j\in C.   
	\end{align}
	Problem \eqref{optimization_3} has optimal solution(s) of the form $\left(w_j^*,{\epsilon_{Bi}^+}^*,{\epsilon_{Bi}^-}^*,{\epsilon_{iW}^+}^*,{\epsilon_{iW}^-}^*,{\epsilon_{BW}^+}^*,\right.$ $\left.{\epsilon_{BW}^-}^*\right)$, where $i\in D$ and $j\in C$. Note that the function $f$ from the collection of optimal solutions of problem \eqref{optimization_3} to the collection of optimal solutions of problem \eqref{optimization_2}, defined by 
	\begin{align*}
		&f(w_j^*,{\epsilon_{Bi}^+}^*,{\epsilon_{Bi}^-}^*,{\epsilon_{iW}^+}^*,{\epsilon_{iW}^-}^*,{\epsilon_{BW}^+}^*,{\epsilon_{BW}^-}^*)\\
		&\quad\quad\quad\quad\quad\quad\quad\quad\quad\quad\quad\quad\quad=(w_j^*,{\epsilon_{Bi}^+}^*+{\epsilon_{Bi}^-}^*,{\epsilon_{iW}^+}^*+{\epsilon_{iW}^-}^*,{\epsilon_{BW}^+}^*+{\epsilon_{BW}^-}^*)
	\end{align*}
	is a well-defined, one-to-one correspondence. This indicates that the taxicab BWM and the nonlinear goal programming model for BWM are equivalent.
	\section{Analytical Framework}
	This section establishes the analytical framework for the taxicab BWM, detailing the process for calculating optimal weights, performing consistency analysis, and determining threshold values. The framework is illustrated and validated through numerical examples. A decision-maker-aided strategy for selecting among multiple optimal weight sets is also proposed. Furthermore, a comparative evaluation is conducted against both traditional optimization software and the nonlinear BWM. The practical applicability of the framework is demonstrated through a smartphone selection problem.
	\subsection{Calculation of optimal weights}
	To compute optimal weights for the taxicab BWM analytically, we first consider the following minimization problem, formulated based on the optimal modification of the given PCS.
	\begin{align}
		&\min \sum_{i\in D} (|\tilde{a}_{Bi}-a_{Bi}|+|\tilde{a}_{iW}-a_{iW}|)+ |\tilde{a}_{BW}-a_{BW}|\nonumber\\
		\label{optimization_4}
		&\text{sub to: \quad} \tilde{a}_{Bi}\times \tilde{a}_{iW}=\tilde{a}_{BW},\quad   \tilde{a}_{Bi},\tilde{a}_{iW},\tilde{a}_{BW}\geq 0\text{ for all }i\in D.
	\end{align}
	Note that problem \eqref{optimization_4} is a nonlinear problem having $2n-3$ variables $\tilde{a}_{Bi}$, $\tilde{a}_{iW}$ and $\tilde{a}_{BW}$, where $i\in D$. So, it has optimal solution(s) of the form $(\tilde{a}_{Bi}^*,\tilde{a}_{iW}^*,\tilde{a}_{BW}^*)$, where $i\in D$. For each optimal solution, the optimal comparison values, along with $\tilde{a}_{BB}^*=\tilde{a}_{WW}^*=1$, form a consistent PCS, referred to as an optimally modified PCS. The optimal objective value indicates the total deviation between the optimal and the given comparison values. Now, observe that this problem is equivalent to the following minimization problem.
	\begin{align}
		&\min \eta=\sum_{i\in D} (\eta_{Bi}+\eta_{iW})+ \eta_{BW}\nonumber\\
		&\text{sub to: \quad} |\tilde{a}_{Bi}-a_{Bi}|=\eta_{Bi},\quad |\tilde{a}_{iW}-a_{iW}|=\eta_{iW},\quad |\tilde{a}_{BW}-a_{BW}|=\eta_{BW},\nonumber\\
		\label{optimization_5}
		&\quad \quad \quad\ \quad   \tilde{a}_{Bi}\times \tilde{a}_{iW}=\tilde{a}_{BW},\quad \tilde{a}_{Bi},\tilde{a}_{iW},\tilde{a}_{BW}\geq 0\text{ for all }i\in D.
	\end{align}
	This problem has optimal solution(s) of the form $(\tilde{a}_{Bi}^*,\tilde{a}_{iW}^*,\tilde{a}_{BW}^*,\eta_{Bi}^*,\eta_{iW}^*,\eta_{BW}^*)$, where $i\in D$, with the optimal objective value $\eta^*$. Similar to problem \eqref{optimization_4}, for each optimal solution, $\tilde{a}_{Bi}^*,\tilde{a}_{iW}^*$ and $\tilde{a}_{BW}^*$, along with $\tilde{a}_{BB}^*=\tilde{a}_{WW}^*=1$, form an optimally modified PCS and $\eta^*$ is the total deviation between the optimal and the given comparison values.\\\\
	Now, we establish a one-to-one correspondence between the collections of optimal solutions of problem \eqref{optimization_2} and problem \eqref{optimization_5}.\\\\
	Let $(w_j^*,\epsilon_{Bi}^*,\epsilon_{iW}^*,\epsilon_{BW}^*)$, where $i\in D$ and $j\in C$, be an optimal solution of problem \eqref{optimization_2}. So, we have $\bigg|\frac{w_{B}^*}{w_i^*}-a_{Bi}\bigg|=\epsilon_{Bi}^*$, $\bigg|\frac{w_i^*}{w_{W}^*}-a_{iW}\bigg|=\epsilon_{iW}^*$ and $\bigg|\frac{w_{B}^*}{w_{W}^*}-a_{BW}\bigg| =\epsilon_{BW}^*$ for all $i\in D$. Take 
	\begin{equation}\label{weight_to_pcs}
		\tilde{a}_{Bi}=\frac{w_{B}^*}{w_i^*}, \quad \tilde{a}_{iW}=\frac{w_i^*}{w_{W}^*}\quad \text{and} \quad \tilde{a}_{BW}=\frac{w_{B}^*}{w_{W}^*}
	\end{equation} 
	for all $i\in D$. Thus, we get $|\tilde{a}_{Bi}-a_{Bi}|=\epsilon_{Bi}^*$, $|\tilde{a}_{iW}-a_{iW}|=\epsilon_{iW}^*$ and $|\tilde{a}_{BW}-a_{BW}| =\epsilon_{BW}^*$ for all $i\in D$. Also, $\tilde{a}_{Bi}\times\tilde{a}_{iW}=\tilde{a}_{BW}$ for all $i\in D$. This gives $\eta^*\leq \displaystyle\sum_{i\in D}(\epsilon_{Bi}^*+\epsilon_{iW}^*)+\epsilon_{BW}^*=\epsilon^*$.\\\\
	Let $(\tilde{a}_{Bi}^*,\tilde{a}_{iW}^*,\tilde{a}_{BW}^*,\eta_{Bi}^*,\eta_{iW}^*,\eta_{BW}^*)$, where $i\in D$, be an optimal solution of problem \eqref{optimization_5}. So, we have $|\tilde{a}_{Bi}^*-a_{Bi}|=\eta_{Bi}^*$, $|\tilde{a}_{iW}^*-a_{iW}|=\eta_{iW}^*$ and $|\tilde{a}_{BW}^*-a_{BW}| =\eta_{BW}^*$ for all $i\in D$. This gives $\tilde{a}_{Bi}^*\in\{a_{Bi}+\eta_{Bi}^*,a_{Bi}-\eta_{Bi}^*\}$, $\tilde{a}_{iW}^*\in\{a_{iW}+\eta_{iW}^*,a_{iW}-\eta_{iW}^*\}$ and $\tilde{a}_{BW}^*\in\{a_{BW}+\eta_{BW}^*,a_{BW}-\eta_{BW}^*\}$. Since $\tilde{a}_{Bi}^*,\tilde{a}_{iW}^*$ and $\tilde{a}_{BW}^*$, along with $\tilde{a}_{BB}^*=\tilde{a}_{WW}^*=1$, form a consistent PCS, by Theorem \ref{accurate}, 
	\begin{equation}\label{pcs_to_weights}
		w_j=\frac{\tilde{a}_{jW}^*}{\displaystyle\sum_{i\in C}\tilde{a}_{iW}^*},\ j\in C
	\end{equation} 
	is the solution of the system of equations $$\frac{w_{B}}{w_i}=\tilde{a}_{Bi}^*,\quad  \frac{w_i}{w_{W}}=\tilde{a}_{iW}^*,\quad  \frac{w_{B}}{w_{W}}=\tilde{a}_{BW}^*,\quad  i\in D.$$Using the values of $\tilde{a}_{Bi}^*$, $\tilde{a}_{iW}^*$ and $\tilde{a}_{BW}^*$, we get $$\frac{w_{B}}{w_i}-a_{Bi}\in\{\eta_{Bi}^*,-\eta_{Bi}^*\},\  \frac{w_i}{w_{W}}-a_{iW}\in\{\eta_{iW}^*,-\eta_{iW}^*\}, \  \frac{w_{B}}{w_{W}}-a_{BW}\in\{\eta_{BW}^*,-\eta_{BW}^*\}$$ for all $i\in D$. So,  $$\bigg|\frac{w_{B}}{w_i}-a_{Bi}\bigg|=\eta_{Bi}^* \quad  \bigg|\frac{w_i}{w_{W}}-a_{iW}\bigg|=\eta_{iW}^*, \quad \bigg|\frac{w_{B}}{w_{W}}-a_{BW}\bigg|=\eta_{BW}^*$$ for all $i\in D$. This gives $\epsilon^*\leq \displaystyle\sum_{i\in D}(\eta_{Bi}^*+\eta_{iW}^*)+\eta_{BW}^*=\eta^*$.\\\\
	From the above discussion, it follows that $\epsilon^*=\eta^*$. Therefore, $\tilde{a}_{Bi}$, $\tilde{a}_{iW}$ and $\tilde{a}_{BW}$ given by Eq. \eqref{weight_to_pcs}, along with $\epsilon_{Bi}^*$, $\epsilon_{iW}^*$ and $\epsilon_{BW}^*$, form an optimal solution of problem \eqref{optimization_5}. Similarly, $w_j$ defined by Eq. \eqref{pcs_to_weights}, along with $\eta_{Bi}^*$, $\eta_{iW}^*$ and $\eta_{BW}^*$, form an optimal solution of problem \eqref{optimization_2}. So, for every $(\tilde{a}_{Bi}^*,\tilde{a}_{iW}^*,\tilde{a}_{BW}^*,\eta_{Bi}^*,\eta_{iW}^*,\eta_{BW}^*)$, there exists unique $(w_j^*,\epsilon_{Bi}^*,\epsilon_{iW}^*,\epsilon_{BW}^*)$ such that 
	\begin{equation*}
		\tilde{a}_{Bi}^*=\frac{w_{B}^*}{w_i^*}, \quad \tilde{a}_{iW}^*=\frac{w_i^*}{w_{W}^*},\quad \tilde{a}_{BW}^*=\frac{w_{B}^*}{w_{W}^*}, \quad \epsilon_{Bi}^*=\eta_{Bi}^*,\quad \epsilon_{iW}^*=\eta_{iW}^*,\quad \epsilon_{BW}^*=\eta_{BW}^*
	\end{equation*}
	for all $i\in D$. Thus, to obtain an analytical expression for the optimal solution(s) of problem \eqref{optimization_2}, it is sufficient to derive an analytical expression for the optimal solution(s) of problem \eqref{optimization_5}.
	\begin{remark}
		Based on Theorem \ref{accurate}, which establishes a unique weight set for every consistent PCS, a fundamental insight emerges: the process of obtaining weight sets by minimizing the TD between weight ratios and a given PCS (problem \eqref{optimization_2}) should be equivalent to the process of converting the given PCS into a consistent one by minimizing the TD between them (problem \eqref{optimization_5}). This equivalence between the two optimization approaches has been formally demonstrated above.
	\end{remark}
	\begin{proposition}\label{geq1}
		Let $(A_{B},A_{W})$ be a given PCS, and let $(\tilde{A}_{B},\tilde{A}_{W})$ be a consistent PCS having $\tilde{a}_{BW}<1$. Then there exist a consistent $(\tilde{A}_{B}',\tilde{A}_{W}')$ having $\tilde{a}_{BW}'=1$ such that $|\tilde{a}_{Bi}'-a_{Bi}|\leq|\tilde{a}_{Bi}-a_{Bi}|$, $|\tilde{a}_{iW}'-a_{iW}|\leq|\tilde{a}_{iW}-a_{iW}|$ and $|\tilde{a}_{BW}'-a_{BW}|<|\tilde{a}_{BW}-a_{BW}|$ for all $i\in D$.
	\end{proposition}
	Proofs of all Propositions and Theorems are provided in the Appendix.\\\\
	Let $(A_{B}^*,A_{W}^*)$ be an optimally modified PCS. Then, by Proposition \ref{geq1}, we get $\tilde{a}_{BW}^*\geq 1$. 
	\begin{definition}\cite{wu2023analytical}
		Let $i\in D$. Then $i$ is said to be consistent criterion if $a_{Bi}\times a_{iW}=a_{BW}$. Similarly, $i$ is called downside criterion if $a_{Bi}\times a_{iW}<a_{BW}$ and upside criterion if $a_{Bi}\times a_{iW}>a_{BW}$.
	\end{definition}
	\begin{definition}
		An optimal modification strategy for $(a_{Bi},a_{iW},a_{BW})$, $i\in D$, is $(x^*,y^*,z^*)\in \mathbb{R}^3$ such that $(a_{Bi}+x^*)\times (a_{iW}+y^*)=a_{BW}+z^*$ and $|x^*|+|y^*|+|z^*|=\inf\{|x|+|y|+|z|: (a_{Bi}+x)\times (a_{iW}+y)=a_{BW}+z\}$.
	\end{definition}
	It is clear that if $i$ is consistent criterion, then the only optimal modification strategy for $(a_{Bi},a_{iW},a_{BW})$ is $(x^*,y^*,z^*)=(0,0,0)$, and $\inf\{|x|+|y|+|z|: (a_{Bi}+x)\times (a_{iW}+y)=a_{BW}+z\}=0$. Also, the optimally modified $(a_{Bi},a_{iW},a_{BW})$ is $(a_{Bi}+0,a_{iW}+0,a_{BW}+0)=(a_{Bi},a_{iW},a_{BW})$.\\\\
	Now, we deal with downside criteria.
	\begin{proposition}\label{min_1}
		Let $a,b \in \{1,2,\ldots,9\}$ and $c\geq 1$ be such that $a\times b< c$, and let $(x,y,z)\in \mathbb{R}^3$ be such that $(a+x)\times (b+y)=c+z$. Then at least one of the following statements holds.
		\begin{enumerate}
			\item $x,y\geq 0$, $z\leq 0$.
			\item There exist $(x',y',z')\in \mathbb{R}^3$ such that $x',y'\geq 0$, $z'\leq 0$, $(a+x')\times(b+y')=c+z'$ and $|x'|+|y'|+|z'|<|x|+|y|+|z|$.
		\end{enumerate}
	\end{proposition}
	\begin{theorem}\label{a=b}
		Let $a\in\{1,2,\ldots,9\}$ and $c\geq 1$ be such that $a\times a<c$, let $x'>0$ be such that $(a+x')\times (a+x')=c$, i.e., $x'=\sqrt{c}-a$, and let $(x,y,z)\neq(x',x',0)$ be such that $x,y,z\geq 0$ and $(a+x)\times(a+y)=c-z$. Then $2x'<x+y+z$.
	\end{theorem}
	From Proposition \ref{min_1} and Theorem \ref{a=b}, it follows that for a downside criterion $i$, if $a_{Bi}=a_{iW}$, then the only optimal modification strategy for $(a_{Bi},a_{iW},a_{BW})$ is $(x^*,y^*,z^*)=(\sqrt{a_{BW}}-a_{Bi},\sqrt{a_{BW}}-a_{iW},0)$, and thus, $\inf\{|x|+|y|+|z|: (a_{Bi}+x)\times (a_{iW}+y)=a_{BW}+z\}=2\sqrt{a_{BW}}-a_{Bi}-a_{iW}$. Also, the optimally modified $(a_{Bi},a_{iW},a_{BW})$ is $(a_{Bi}+\sqrt{a_{BW}}-a_{Bi},a_{iW}+\sqrt{a_{BW}}-a_{iW},a_{BW}+0)=(\sqrt{a_{BW}},\sqrt{a_{BW}},a_{BW})$.
	\begin{theorem}\label{a<b}
		Let $a,b\in\{1,2,\ldots,9\}$ and $c\geq 1$ be such that $a<b$ and $a\times b<c$, and let $(x,y,z)$ be such that $x,y,z\geq 0$ and $(a+x)\times (b+y)=c-z$. Then the following statements hold.
		\begin{enumerate}
			\item If $b\geq \sqrt{c}$, then $x'<x+y+z$ for $(x,y,z)\neq (x',0,0)$, where $x'>0$ is such that $(a+x')\times b=c$, i.e., $x'=\frac{c}{b}-a$.
			\item If $b< \sqrt{c}$, then $b-a+2y'<x+y+z$ for $(x,y,z)\neq (b-a+y',y',0)$, where $y'>0$ is such that $(b+y')\times (b+y')=c$, i.e., $y'=\sqrt{c}-b$.
		\end{enumerate}
	\end{theorem}
	From Proposition \ref{min_1} and Theorem \ref{a<b}, for a downside criterion $i$, the following conclusions can be drawn.
	\begin{enumerate}
		\item If $a_{Bi}<a_{iW}$ and $\sqrt{a_{BW}}\leq a_{iW}$, then the only optimal modification strategy for $(a_{Bi},a_{iW},a_{BW})$ is $(x^*,y^*,z^*)=(\frac{a_{BW}}{a_{iW}}-a_{Bi},0,0)$, and thus, $\inf\{|x|+|y|+|z|: (a_{Bi}+x)\times (a_{iW}+y)=a_{BW}+z\}=\frac{a_{BW}}{a_{iW}}-a_{Bi}$. Also, the optimally modified $(a_{Bi},a_{iW},a_{BW})$ is $(a_{Bi}+\frac{a_{BW}}{a_{iW}}-a_{Bi},a_{iW},a_{BW})=(\frac{a_{BW}}{a_{iW}},a_{iW},a_{BW})$.
		\item If $a_{Bi}>a_{iW}$ and $\sqrt{a_{BW}}\leq a_{Bi}$, then the only optimal modification strategy for $(a_{Bi},a_{iW},a_{BW})$ is $(x^*,y^*,z^*)=(0,\frac{a_{BW}}{a_{Bi}}-a_{iW},0)$, and thus, $\inf\{|x|+|y|+|z|: (a_{Bi}+x)\times (a_{iW}+y)=a_{BW}+z\}=\frac{a_{BW}}{a_{Bi}}-a_{iW}$. Also, the optimally modified $(a_{Bi},a_{iW},a_{BW})$ is $(a_{Bi},a_{iW}+\frac{a_{BW}}{a_{Bi}}-a_{iW},a_{BW})=(a_{Bi},\frac{a_{BW}}{a_{Bi}},a_{BW})$.
		\item If $a_{Bi}<a_{iW}<\sqrt{a_{BW}}$ or $a_{iW}<a_{Bi}<\sqrt{a_{BW}}$, then the only optimal modification strategy for $(a_{Bi},a_{iW},a_{BW})$ is $(x^*,y^*,z^*)=(\sqrt{a_{BW}}-a_{Bi},\sqrt{a_{BW}}-a_{iW},0)$, and thus, $\inf\{|x|+|y|+|z|: (a_{Bi}+x)\times (a_{iW}+y)=a_{BW}+z\}=2\sqrt{a_{BW}}-a_{Bi}-a_{iW}$. Also, the optimally modified $(a_{Bi},a_{iW},a_{BW})$ is $(a_{Bi}+\sqrt{a_{BW}}-a_{Bi},a_{iW}+\sqrt{a_{BW}}-a_{iW},a_{BW})=(\sqrt{a_{BW}},\sqrt{a_{BW}},a_{BW})$.
	\end{enumerate}
	\begin{proposition}\label{min_2}
		Let $a,b \in \{1,2,\ldots,9\}$ and $c\geq 1$ be such that $a\times b> c$, and let $(x,y,z)\in \mathbb{R}^3$ be such that $(a+x)\times (b+y)=c+z$. Then at least one of the following statements holds.
		\begin{enumerate}
			\item $x,y\leq 0$, $z\geq 0$, $a+x,b+y>0$.
			\item There exist $(x',y',z')\in \mathbb{R}^3$ such that $x',y'\leq 0$, $z'\geq 0$, $a+x',b+y'>0$, $(a+x')\times(b+y')=c+z'$ and $|x'|+|y'|+|z'|<|x|+|y|+|z|$.
		\end{enumerate}
	\end{proposition}
	\begin{theorem}\label{a<=b}
		Let $a,b\in\{1,2,\ldots,9\}$ and $c\geq 1$ be such that $a\times b>c$ and $a\leq b\leq c$, let $x'>0$ be such that $(a-x')\times b=c$, i.e., $x'=a-\frac{c}{b}$, and let $(x,y,z)$ be such that $x,y,z\geq 0$, $a-x,b-y>0$ and $(a-x)\times(b-y)=c+z$. Then the following statements hold.
		\begin{enumerate}
			\item If $a<b$, then $x'<x+y+z$ for $(x,y,z)\neq (x',0,0)$.
			\item If $a=b$, then $x'<x+y+z$ for $(x,y,z)\neq (x',0,0)\neq (0,x',0)$.
		\end{enumerate}
	\end{theorem}
	From Proposition \ref{min_2} and Theorem \ref{a<=b}, for an upside criterion $i$, the following conclusions can be drawn.
	\begin{enumerate}
		\item If $a_{Bi}<a_{iW}$, then the only optimal modification strategy for $(a_{Bi},a_{iW},a_{BW})$ is $(x^*,y^*,z^*)=(\frac{a_{BW}}{a_{iW}}-a_{Bi},0,0)$, and thus, $\inf\{|x|+|y|+|z|: (a_{Bi}+x)\times (a_{iW}+y)=a_{BW}+z\}=a_{Bi}-\frac{a_{BW}}{a_{iW}}$. Also, the optimally modified $(a_{Bi},a_{iW},a_{BW})$ is $(a_{Bi}+\frac{a_{BW}}{a_{iW}}-a_{Bi},a_{iW},a_{BW})=(\frac{a_{BW}}{a_{iW}},a_{iW},a_{BW})$.
		\item If $a_{Bi}>a_{iW}$, then only optimal modification strategy for $(a_{Bi},a_{iW},a_{BW})$ is $(x^*,y^*,z^*)=(0,\frac{a_{BW}}{a_{Bi}}-a_{iW},0)$, and thus, $\inf\{|x|+|y|+|z|: (a_{Bi}+x)\times (a_{iW}+y)=a_{BW}+z\}=a_{iW}-\frac{a_{BW}}{a_{Bi}}$. Also, the optimally modified $(a_{Bi},a_{iW},a_{BW})$ is $(a_{Bi},a_{iW}+\frac{a_{BW}}{a_{Bi}}-a_{iW},a_{BW})=(a_{Bi},\frac{a_{BW}}{a_{Bi}},a_{BW})$.
		\item If $a_{Bi}=a_{iW}$, then $(x^*,y^*,z^*)=(\frac{a_{BW}}{a_{iW}}-a_{Bi},0,0)$ and $(x^*,y^*,z^*)=(0,\frac{a_{BW}}{a_{Bi}}-a_{iW},0)$ are the only optimal modification strategies for $(a_{Bi},a_{iW},a_{BW})$. Note that, for both strategies, we have $\inf\{|x|+|y|+|z|: (a_{Bi}+x)\times (a_{iW}+y)=a_{BW}+z\}=a_{Bi}-\frac{a_{BW}}{a_{iW}}$. Also, optimally modified $(a_{Bi},a_{iW},a_{BW})$ are $(a_{Bi}+\frac{a_{BW}}{a_{iW}}-a_{Bi},a_{iW},a_{BW})=(\frac{a_{BW}}{a_{iW}},a_{iW},a_{BW})$ and $(a_{Bi},a_{iW}+\frac{a_{BW}}{a_{Bi}}-a_{iW},a_{BW})=(a_{Bi},\frac{a_{BW}}{a_{Bi}},a_{BW})$.
	\end{enumerate}
	Note that for all the aforementioned optimally modified $(a_{Bi},a_{iW},a_{BW})$, $a_{BW}$ remains unchanged. Therefore, an optimally modified PCS can be expressed in terms of its $\tilde{a}_{BW}^*$ as follows:
	\begin{align}
		\begin{cases}
			\begin{cases}
				\tilde{a}_{Bi}^*=a_{Bi}\\
				\tilde{a}_{iW}^*=a_{iW}
			\end{cases}\quad\quad \quad \quad \quad \quad \quad \quad & \text{if } a_{Bi}\times a_{iW}=\tilde{a}_{BW}^*,\nonumber\\
			\begin{cases}
				\tilde{a}_{Bi}^*=\sqrt{\tilde{a}_{BW}^*}\\
				\tilde{a}_{iW}^*=\sqrt{\tilde{a}_{BW}^*}
			\end{cases}\quad \quad \quad \quad \quad \quad \quad& \parbox[t]{0.45\textwidth}{\text{if } $a_{Bi}\times a_{iW}< \tilde{a}_{BW}^*$ \text{ and } $a_{Bi}, a_{iW}<\sqrt{\tilde{a}_{BW}^*}$,}\nonumber\\
			\begin{cases}
				\tilde{a}_{Bi}^*=\frac{\tilde{a}_{BW}^*}{a_{iW}}\\
				\tilde{a}_{iW}^*=a_{iW}
			\end{cases}\quad \quad \quad \quad \quad \quad \quad \quad& \parbox[t]{0.45\textwidth}{\text{if } ($a_{Bi}\times a_{iW}<\tilde{a}_{BW}^*$ and $a_{Bi}<\sqrt{\tilde{a}_{BW}^*}\leq a_{iW}$) or ($a_{Bi}\times a_{iW}>\tilde{a}_{BW}^*$ and $a_{Bi}<a_{iW}$),}\nonumber\\
			\begin{cases}
				\tilde{a}_{Bi}^*=a_{Bi}\\
				\tilde{a}_{iW}^*=\frac{\tilde{a}_{BW}^*}{a_{Bi}}
			\end{cases}\quad \quad \quad \quad \quad \quad \quad \quad& \parbox[t]{0.45\textwidth}{\text{if } ($a_{Bi}\times a_{iW}<\tilde{a}_{BW}^*$ and $a_{iW}<\sqrt{\tilde{a}_{BW}^*}\leq a_{Bi}$) or ($a_{Bi}\times a_{iW}>\tilde{a}_{BW}^*$ and $a_{iW}<a_{Bi}$),}\nonumber\\
			\begin{cases}
				\tilde{a}_{Bi}^*=\frac{\tilde{a}_{BW}^*}{a_{iW}}\\
				\tilde{a}_{iW}^*=a_{iW}
			\end{cases} \text{ or }
			\begin{cases}
				\tilde{a}_{Bi}^*=a_{Bi}\\
				\tilde{a}_{iW}^*=\frac{\tilde{a}_{BW}^*}{a_{Bi}}
			\end{cases} \quad &\text{if } a_{Bi}\times a_{iW}>\tilde{a}_{BW}^* \text{ and } a_{Bi}=a_{iW},\nonumber\\
			\begin{cases}
				\tilde{a}_{BB}^*=\tilde{a}_{WW}^*=1
			\end{cases}
		\end{cases}\\
		\label{optimal_pcs}
	\end{align}
	where $i\in D$; therefore, we get
	\begin{align}
		\begin{cases}
			\begin{cases}
				\eta_{bi}^*=0\\
				\eta_{iw}^*=0
			\end{cases}\quad\quad \quad \quad \quad \quad \quad \quad \quad \quad \quad \quad\quad\ & \text{if } a_{Bi}\times a_{iW}=\tilde{a}_{BW}^*,\nonumber\\
			\begin{cases}
				\eta_{bi}^*=\sqrt{\tilde{a}_{BW}^*}-a_{Bi}\\
				\eta_{iw}^*=\sqrt{\tilde{a}_{BW}^*}-a_{iW}
			\end{cases}\quad \quad \quad \quad\quad \quad \quad \quad\quad& \parbox[t]{0.35\textwidth}{\text{if } $a_{Bi}\times a_{iW}< \tilde{a}_{BW}^*$ \text{ and } $a_{Bi},a_{iW} <\sqrt{\tilde{a}_{BW}^*}$,}\nonumber\\
			\begin{cases}
				\eta_{bi}^*=\bigg|a_{Bi}-\frac{\tilde{a}_{BW}^*}{a_{iW}}\bigg|\\
				\eta_{iw}^*=0
			\end{cases}\quad \quad \quad \quad \quad \quad \quad \quad\quad\ \ & \parbox[t]{0.35\textwidth}{\text{if } ($a_{Bi}\times a_{iW}<\tilde{a}_{BW}^*$ and $a_{Bi}<\sqrt{\tilde{a}_{BW}^*}\leq a_{iW}$) or ($a_{Bi}\times a_{iW}>\tilde{a}_{BW}^*$ and $a_{Bi}<a_{iW}$),}\nonumber\\
			\begin{cases}
				\eta_{bi}^*=0\\
				\eta_{iw}^*=\bigg|a_{iW}-\frac{\tilde{a}_{BW}^*}{a_{Bi}}\bigg|
			\end{cases}\quad \quad \quad \quad \quad \quad \quad \quad\quad\ & \parbox[t]{0.35\textwidth}{\text{if } ($a_{Bi}\times a_{iW}<\tilde{a}_{BW}^*$ and $a_{iW}<\sqrt{\tilde{a}_{BW}^*}\leq a_{Bi}$) or ($a_{Bi}\times a_{iW}>\tilde{a}_{BW}^*$ and $a_{iW}<a_{Bi}$),}\nonumber\\
			\begin{cases}
				\eta_{bi}^*=a_{Bi}-\frac{\tilde{a}_{BW}^*}{a_{iW}}\\
				\eta_{iw}^*=0
			\end{cases} \text{ or }
			\begin{cases}
				\eta_{bi}^*=0\\
				\eta_{iw}^*=a_{iW}-\frac{\tilde{a}_{BW}^*}{a_{Bi}}
			\end{cases} \quad &\parbox[t]{0.35\textwidth}{\text{if } $a_{Bi}\times a_{iW}>\tilde{a}_{BW}^*$ \text{ and } $a_{Bi}=a_{iW}$,}\nonumber\\
			\begin{cases}
				\eta_{bw}^*=|a_{BW}-\tilde{a}_{BW}^*|	
			\end{cases}
		\end{cases}\nonumber\\
		\label{optimal_modification}
	\end{align}
	for all $i\in D$. Thus, to obtain analytical form of optimally modified PCS, it is sufficient to determine all possible values of $\tilde{a}_{BW}^*$. Also, the analytical expression of optimal objective value of problem \eqref{optimization_5}, and thus of problem \eqref{optimization_2}, is 
	\begin{equation}\label{optimal_obj_value}
		\epsilon^*=\eta^*=\displaystyle\sum_{i\in D} (\eta_{Bi}^*+\eta_{iW}^*)+\eta_{BW}^*.
	\end{equation}
	For $x\in [1,\infty)$ and $i\in D$, define
	\begin{align}
		&f_i(x)=
		\begin{cases}
			\left|a_{iW}-\frac{x}{a_{Bi}}\right|\quad \ \quad \quad \text{if } 1\leq x\leq a_{Bi}^2\text{ and }a_{iW}\leq a_{Bi},\\
			\left|a_{Bi}-\frac{x}{a_{iW}}\right|\quad \ \quad \quad \text{if } 1\leq x\leq a_{iW}^2\text{ and }a_{Bi}\leq a_{iW},\\
			2\sqrt{x}-a_{Bi}-a_{iW}\quad \text{otherwise},
		\end{cases}\nonumber\\
		&f_{B}(x)=|a_{BW}-x|\quad \text{ and }\nonumber\\
		\label{function}
		&f(x)=\displaystyle\sum_{i\in D}f_i(x)+f_{B}(x).
	\end{align}
	Note that $f_i,f_{B}$ and $f$ are continuous functions. Furthermore, it can be observed that the global minimum value of $f$ is the same as the optimal objective value of problem \eqref{optimization_2}, and the points at which $f$ attains this global minimum represent all possible values of $\tilde{a}_{BW}^*$.\\\\
	Let $u=\max\{a_{Bi}\times a_{iW}, a_{BW}: i\in D\}$. Consider
	\begin{equation}\label{set}
		X=\{a_{Bi}\times a_{iW}, a_{BW}: i\in D\}\cup \{\max\{a_{Bi}^2,a_{iW}^2\}:\max\{a_{Bi}^2,a_{iW}^2\}\leq u,i\in D\}.
	\end{equation}
	Since $X$ is finite, it can be expressed as $X=\{x_0,x_1,\ldots,x_m\}$, where $x_0<x_1<\ldots<x_m$. Now, $a_{Bi},a_{iW},a_{BW}\geq 1$ for all $i\in D$ implies that $x_0\geq 1$. Thus,
	\begin{equation}\label{decomposition}
		[1,\infty)=[1,x_0]\cup [x_0,x_1]\cup \ldots \cup [x_{m-1},x_m]\cup [x_m,\infty).
	\end{equation}
	\begin{theorem}\label{minima}
		Let $f$ and $X=\{x_0,x_1,\ldots,x_m\}$ be as in Eq. \eqref{function} and \eqref{set} respectively. Then $f$ attains its global minimum at some $x_j\in X$. Furthermore, if $f$ is nonconstant on each interval $[x_{j-1},x_j]$ for $j=1,2,\ldots,m$, then this global minimum is achieved only at some $x_j\in X$. 
	\end{theorem}
	From Theorem \eqref{minima}, it follows that if $f$ attains its global minimum at $x_{j-1},x_{j}\in X$ for some $j$ and $f$ is constant on $[x_{j-1},x_j]$, then the interval $(x_{j-1},x_j)$, along with all points of $X$ where f achieves its global minimum, constitute the possible values of $\tilde{a}_{BW}^*$. Otherwise, the only possible values of $\tilde{a}_{BW}^*$ are the points of $X$ where f achieves its global minimum. After obtaining all possible values of $\tilde{a}_{BW}^*$, the collection of optimally modified PCS is obtained using Eq. \eqref{optimal_pcs}. Subsequently, the collection of optimal weight sets is determined using Eq. \eqref{pcs_to_weights}, and the optimal TD is calculated using Eq. \eqref{optimal_modification} and Eq. \eqref{optimal_obj_value}.
	\subsection{Consistency index}
	The resultant weights depend on pairwise comparisons, which may exhibit inconsistency due to human involvement. This inconsistency is estimated using a ratio known as the Consistency Ratio (CR) defined as
	\begin{equation}\label{cr}
		\text{CR}=\frac{\epsilon^*}{\text{Consistency Index (CI)}},
	\end{equation}
	\begin{align*}
		&\text{where CI} =\sup\left\{\epsilon^*: \epsilon^*\text{ is the optimal objective value of problem \eqref{optimization_2} for some }\right.\\ 
		&\ \ \left.(A_{B},A_{W})  \text{ having the given number of criteria } n \text{ and the given value of } a_{BW}\right\}.\text{\cite{rezaei2015best}}
	\end{align*} So, CI is a function of $a_{BW}$ and $n$. In this subsection, our goal is to obtain the values of CI$_{a_{BW}}(n)$ in the context of the taxicab BWM.\\\\
	Fix $a_{BW}$ and $n$. Let $\mathcal{A}_{a_{BW},n}$ be the collection of all PCSs having the given $a_{BW}$ with $n$ criteria. Consider
	\begin{equation}\label{set1}
		X'=\{a\times b: a,b=1,2,\ldots,a_{BW}\}.
	\end{equation} Let $(A_{B},A_{W})\in\mathcal{A}_{a_{BW},n}$, let $\epsilon^*$ be the corresponding optimal objective value of problem \eqref{optimization_2}, and let $X$ be the set defined by Eq. \eqref{set}. Note that $X\subseteq X'$. This, along with the fact that $f$ attains its global minimum value at some point in $X$, implies that $\epsilon^*=\displaystyle\min_{x\in X'}f(x)$.\\\\
	For all $i\in D$, $(a_{Bi},a_{iW})$ is of the form $(a,b)$, where $a,b\in\{1,2,\ldots,a_{BW}\}$. Eq. \eqref{optimal_modification} and Eq. \eqref{optimal_obj_value} imply that interchanging $a_{Bi}$ and $a_{iW}$ does not affect the optimal objective value. Thus, without loss of generality, we may assume $a\leq b$. Let $n_{a,b}$ denote the number of pairs $(a,b)$ in $(A_{B},A_{W})$. Then, by Eq. \eqref{function},
	\begin{equation}\label{function1}
		f(x)=\displaystyle\sum_{\substack{a,b=1\\a\leq b}}^{a_{BW}}n_{a,b}\times f_{a,b}(x)+|a_{BW}-x|, \text{ where } f_{a,b}(x)=\begin{cases}
			\left|a-\frac{x}{b}\right| \ \quad \quad \text{if } 1\leq x\leq b^2,\\
			2\sqrt{x}-a-b\  \text{otherwise.}
		\end{cases}
	\end{equation}
	This gives $\epsilon^*=\displaystyle\min_{x\in X'}\biggl\{\displaystyle\sum_{\substack{a,b=1\\a\leq b}}^{a_{BW}}n_{a,b}\times f_{a,b}(x)+|a_{BW}-x|\biggr\}$. By definition of CI, we have $\text{CI}_{a_{BW}}(n)=\displaystyle\max_{(A_{B},A_{W})\in\mathcal{A}_{a_{BW},n}}\biggl\{\displaystyle\min_{x\in X'}\biggl\{\displaystyle\sum_{\substack{a,b=1\\a\leq b}}^{a_{BW}}n_{a,b}\times f_{a,b}(x)+|a_{BW}-x|\biggr\}\biggr\}$. To obtain CI$_{a_{BW}}(n)$, consider the following mixed-integer linear optimization problem.
	\begin{align}
		&\max\biggl\{\displaystyle\min_{x\in X'}\displaystyle\sum_{\substack{a,b=1\\a\leq b}}^{a_{BW}}n_{a,b}\times f_{a,b}(x)+|a_{BW}-x|\biggr\} \nonumber\\
		\label{optimization_6}
		&\text{sub to: \quad}\displaystyle\sum_{\substack{a,b=1\\a\leq b}}^{a_{BW}}n_{a,b}+2=n,\quad   n_{a,b}\in \mathbb{N}\cup\{0\} \text{ for all }a,b.
	\end{align}
	Problem \eqref{optimization_6} has $\frac{n(n+1)}{2}$ variables $n_{a,b}$, where $a,b=1,2,\ldots,a_{BW}$ and $a\leq b$. Observe that the optimal objective value is precisely CI$_{a_{BW}}(n)$, and an optimal solution $n_{a,b}^*$ yields a PCS with $n$ criteria and the given $a_{BW}$ having $\epsilon^*=$ CI$_{a_{BW}}(n)$. Now, consider another mixed-integer linear optimization problem.
	\begin{align}
		&\max z \nonumber\\
		&\text{sub to: \quad}\displaystyle\sum_{\substack{a,b=1\\a\leq b}}^{a_{BW}}n_{a,b}\times f_{a,b}(x)+|a_{BW}-x|\geq z\text{ for all }x\in X',\nonumber\\
		\label{optimization_7}
		&\quad \quad \quad\ \quad  \displaystyle\sum_{\substack{a,b=1\\a\leq b}}^{a_{BW}}n_{a,b}+2=n,\quad   n_{a,b}\in \mathbb{N}\cup\{0\} \text{ for all }a,b.
	\end{align}
	Problem \eqref{optimization_7} is an equivalent formulation of problem \eqref{optimization_6} with $\frac{n(n+1)}{2}+1$ variables $n_{a,b}$ and $z$, where $a,b=1,2,\ldots,a_{BW}$ and $a\leq b$. So, it has an optimal solution of the form $(n_{a,b}^*,z^*)$. Here, $z^*$, which is also the optimal objective value, is precisely CI$_{a_{BW}}(n)$, while $n_{a,b}^*$ defines a PCS with $n$ criteria and the given $a_{BW}$ having $\epsilon^*=$ CI$_{a_{BW}}(n)$.\\\\
	We now compute CI$_{a_{BW}}(n)$ for two parameter sets to demonstrate the methodology. First, consider $n=5$ and $a_{BW}=2$. Then, by Eq. \eqref{set1}, $X'=\{1,2,4\}$. Using Eq. \eqref{function1}, we calculate the values of $f_{a,b}(x)$ for $a,b=1,2$ with $a\leq b$ and $x\in X'$, which are given in Table \ref{3ci_ex_1}.
	  
	\begin{table}[t!]
		\caption{The values of $f_{a,b}(x)$ for $a_{BW}=2$}\label{3ci_ex_1}
		\centering
		\begin{tabular}{@{}ccccc@{}}
			\toprule
			\multirow{2}{*}{$x\downarrow$}&\phantom{}&\multicolumn{3}{c}{$f_{a,b}\rightarrow$}\\
			\cmidrule{3-5}
			&&$f_{1,1}$&$f_{1,2}$&$f_{2,2}$\\
			\midrule
			$1$&&$0$&$0.5$&$1.5$\\
			$2$&&$2\sqrt{2}-2$&$0$&$1$\\
			$4$&&$2$&$1$&$0$\\
			\bottomrule				
		\end{tabular}
		\par\smallskip \small \textit{These values serve as the coefficients in problem \eqref{optimization_7} for computing CI$_{a_{BW}}(n)$ when $a_{BW}=2$.}
	\end{table}
	 
	\hspace{-0.7cm}
	Thus, problem \eqref{optimization_7} takes the form 
	\begin{align}
		&\max z \nonumber\\
		&\text{sub to: \quad}0.5 n_{1,2}+1.5 n_{2,2}+1\geq z,\quad (2\sqrt{2}-2) n_{1,1}+n_{2,2}\geq z,\nonumber\\
		&\quad \quad \quad\ \quad 2 n_{1,1}+n_{1,2}+2\geq z,\quad n_{1,1}+n_{1,2}+n_{2,2}+2=5,\nonumber\\
		\label{optimization_8}
		&\quad \quad \quad\ \quad n_{1,1},n_{1,2},n_{2,2}\in \mathbb{N}\cup\{0\}.
	\end{align}
	The optimal solution of this problem is $(n_{1,1}^*,n_{1,2}^*,n_{2,2}^*,z^*)=(1,0,2,2.8284)$. Thus, CI$_2(5)=2.8284$, and the best-to-other vector $A_{B}=(1,1,2,2,2)$ and the other-to-worst vector $A_{W}=(2,1,2,2,1)^T$ with $c_1$ as the best and $c_5$ as the worst criterion forms a PCS with $\epsilon^*=2.8284$.\\\\
	Now, consider $n=15$ and $a_{BW}=7$. Then, by Eq. \eqref{set1}, $X'=\left\{1,2,3,4,5,\right.$ $\left.6,7,8,9,10,12,14,15,16,18,20,21,24,25,28,30,35,36,42,49\right\}$. Using Eq. \eqref{function1}, we calculate the values of $f_{a,b}(x)$ for $a,b=1,2,\ldots,7$ with $a\leq b$ and $x\in X'$, which are given in Table \ref{3ci_ex_2}. 
	\begin{sidewaystable}[ph]
		\caption{The values of $f_{a,b}(x)$ for $a_{BW}=7$}\label{3ci_ex_2}
		\centering
		\begin{adjustbox}{width=\textwidth}
			\small		
			\begin{tabular}{@{}cccccccccccccccccccccccccccccc@{}}
				\toprule
				\multirow{2}{*}{$x\downarrow$}&\phantom{}&\multicolumn{28}{c}{$f_{a,b}\rightarrow$}\\
				\cmidrule{3-30}
				&&$f_{1,1}$&$f_{1,2}$&$f_{1,3}$&$f_{1,4}$&$f_{1,5}$&$f_{1,6}$&$f_{1,7}$&$f_{2,2}$&$f_{2,3}$&$f_{2,4}$&$f_{2,5}$&$f_{2,6}$&$f_{2,7}$&$f_{3,3}$&$f_{3,4}$&$f_{3,5}$&$f_{3,6}$&$f_{3,7}$&$f_{4,4}$&$f_{4,5}$&$f_{4,6}$&$f_{4,7}$&$f_{5,5}$&$f_{5,6}$&$f_{5,7}$&$f_{6,6}$&$f_{6,7}$&$f_{7,7}$\\
				\midrule
				$1$&&$0$&$0.5$&$0.6667$&$0.75$&$0.8$&$0.8333$&$0.8571$&$1.5$&$1.6667$&$1.75$&$1.8$&$1.8333$&$1.8571$&$2.6667$&$2.75$&$2.8$&$2.8333$&$2.8571$&$3.75$&$3.8$&$3.8333$&$3.8571$&$4.8$&$4.8333$&$4.8571$&$5.8333$&$5.8571$&$6.8571$\\
				$2$&&$2\sqrt{2}-2$&$0$&$0.3333$&$0.5$&$0.6$&$0.6667$&$0.7143$&$1$&$1.3333$&$1.5$&$1.6$&$1.6667$&$1.7143$&$2.3333$&$2.5$&$2.6$&$2.6667$&$2.7143$&$3.5$&$3.6$&$3.6667$&$3.7143$&$4.6$&$4.6667$&$4.7143$&$5.6667$&$5.7143$&$6.7143$\\
				$3$&&$2\sqrt{3}-2$&$0.5$&$0$&$0.25$&$0.4$&$0.5$&$0.5714$&$0.5$&$1$&$1.25$&$1.4$&$1.5$&$1.5714$&$2$&$2.25$&$2.4$&$2.5$&$2.5714$&$3.25$&$3.4$&$3.5$&$3.5714$&$4.4$&$4.5$&$4.5714$&$5.5$&$5.5714$&$6.5714$\\
				$4$&&$2$&$1$&$0.3333$&$0$&$0.2$&$0.3333$&$0.4286$&$0$&$0.6667$&$1$&$1.2$&$1.3333$&$1.4286$&$1.6667$&$2$&$2.2$&$2.3333$&$2.4286$&$3$&$3.2$&$3.3333$&$3.4286$&$4.2$&$4.3333$&$4.4286$&$5.3333$&$5.4286$&$6.4286$\\
				$5$&&$2\sqrt{5}-2$&$2\sqrt{5}-3$&$0.6667$&$0.25$&$0$&$0.1667$&$0.2857$&$2\sqrt{5}-4$&$0.3333$&$0.75$&$1$&$1.1667$&$1.2857$&$1.3333$&$1.75$&$2$&$2.1667$&$2.2857$&$2.75$&$3$&$3.1667$&$3.2857$&$4$&$4.1667$&$4.2857$&$5.1667$&$5.2857$&$6.2857$\\
				$6$&&$2\sqrt{6}-2$&$2\sqrt{6}-3$&$1$&$0.5$&$0.2$&$0$&$0.1429$&$2\sqrt{6}-4$&$0$&$0.5$&$0.8$&$1$&$1.1429$&$1$&$1.5$&$1.8$&$2$&$2.1429$&$2.5$&$2.8$&$3$&$3.1429$&$3.8$&$4$&$4.1429$&$5$&$5.1429$&$6.1429$\\
				$7$&&$2\sqrt{7}-2$&$2\sqrt{7}-3$&$1.3333$&$0.75$&$0.4$&$0.1667$&$0$&$2\sqrt{7}-4$&$0.3333$&$0.25$&$0.6$&$0.8333$&$1$&$0.6667$&$1.25$&$1.6$&$1.8333$&$2$&$2.25$&$2.6$&$2.8333$&$3$&$3.6$&$3.8333$&$4$&$4.8333$&$5$&$6$\\
				$8$&&$2\sqrt{8}-2$&$2\sqrt{8}-3$&$1.6667$&$1$&$0.6$&$0.3333$&$0.1429$&$2\sqrt{8}-4$&$0.6667$&$0$&$0.4$&$0.6667$&$0.8571$&$0.3333$&$1$&$1.4$&$1.6667$&$1.8571$&$2$&$2.4$&$2.6667$&$2.8571$&$3.4$&$3.6667$&$3.8571$&$4.6667$&$4.8571$&$5.8571$\\
				$9$&&$4$&$3$&$2$&$1.25$&$0.8$&$0.5$&$0.2857$&$2$&$1$&$0.25$&$0.2$&$0.5$&$0.7143$&$0$&$0.75$&$1.2$&$1.5$&$1.7143$&$1.75$&$2.2$&$2.5$&$2.7143$&$3.2$&$3.5$&$3.7143$&$4.5$&$4.7143$&$5.7143$\\
				$10$&&$2\sqrt{10}-2$&$2\sqrt{10}-3$&$2\sqrt{10}-4$&$1.5$&$1$&$0.6667$&$0.4286$&$2\sqrt{10}-4$&$2\sqrt{10}-5$&$0.5$&$0$&$0.3333$&$0.5714$&$2\sqrt{10}-6$&$0.5$&$1$&$1.3333$&$1.5714$&$1.5$&$2$&$2.3333$&$2.5714$&$3$&$3.3333$&$3.5714$&$4.3333$&$4.5714$&$5.5714$\\
				$12$&&$2\sqrt{12}-2$&$2\sqrt{12}-3$&$2\sqrt{12}-4$&$2$&$1.4$&$1$&$0.7143$&$2\sqrt{12}-4$&$2\sqrt{12}-5$&$1$&$0.4$&$0$&$0.2857$&$2\sqrt{12}-6$&$0$&$0.6$&$1$&$1.2857$&$1$&$1.6$&$2$&$2.2857$&$2.6$&$3$&$3.2857$&$4$&$4.2857$&$5.2857$\\
				$14$&&$2\sqrt{14}-2$&$2\sqrt{14}-3$&$2\sqrt{14}-4$&$2.5$&$1.8$&$1.3333$&$1$&$2\sqrt{14}-4$&$2\sqrt{14}-5$&$1.5$&$0.8$&$0.3333$&$0$&$2\sqrt{14}-6$&$0.5$&$0.2$&$0.6667$&$1$&$0.5$&$1.2$&$1.6667$&$2$&$2.2$&$2.6667$&$3$&$3.6667$&$4$&$5$\\
				$15$&&$2\sqrt{15}-2$&$2\sqrt{15}-3$&$2\sqrt{15}-4$&$2.75$&$2$&$1.5$&$1.1429$&$2\sqrt{15}-4$&$2\sqrt{15}-5$&$1.75$&$1$&$0.5$&$0.1429$&$2\sqrt{15}-6$&$0.75$&$0$&$0.5$&$0.8571$&$0.25$&$1$&$1.5$&$.8571$&$2$&$2.5$&$2.8571$&$3.5$&$3.8571$&$4.8571$\\
				$16$&&$6$&$5$&$4$&$3$&$2.2$&$1.6667$&$1.2857$&$4$&$3$&$2$&$1.2$&$0.6667$&$0.2857$&$2$&$1$&$0.2$&$0.3333$&$0.7143$&$0$&$0.8$&$1.3333$&$1.7143$&$1.8$&$2.3333$&$2.7143$&$3.3333$&$3.7143$&$4.7143$\\
				$18$&&$2\sqrt{18}-2$&$2\sqrt{18}-3$&$2\sqrt{18}-4$&$2\sqrt{18}-5$&$2.6$&$2$&$1.5714$&$2\sqrt{18}-4$&$2\sqrt{18}-5$&$2\sqrt{18}-6$&$1.6$&$1$&$0.5714$&$2\sqrt{18}-6$&$2\sqrt{18}-7$&$0.6$&$0$&$0.4286$&$2\sqrt{18}-8$&$0.4$&$1$&$1.4286$&$1.4$&$2$&$2.4286$&$3$&$3.4286$&$4.4286$\\
				$20$&&$2\sqrt{20}-2$&$2\sqrt{20}-3$&$2\sqrt{20}-4$&$2\sqrt{20}-5$&$3$&$2.3333$&$1.8571$&$2\sqrt{20}-4$&$2\sqrt{20}-5$&$2\sqrt{20}-6$&$2$&$1.3333$&$0.8571$&$2\sqrt{20}-6$&$2\sqrt{20}-7$&$1$&$0.3333$&$0.1429$&$2\sqrt{20}-8$&$0$&$0.6667$&$1.1429$&$1$&$1.6667$&$2.1429$&$2.6667$&$3.1429$&$4.1429$\\
				$21$&&$2\sqrt{21}-2$&$2\sqrt{21}-3$&$2\sqrt{21}-4$&$2\sqrt{21}-5$&$3.2$&$2.5$&$2$&$2\sqrt{21}-4$&$2\sqrt{21}-5$&$2\sqrt{21}-6$&$2.2$&$1.5$&$1$&$2\sqrt{21}-6$&$2\sqrt{21}-7$&$1.2$&$0.5$&$0$&$2\sqrt{21}-8$&$0.2$&$0.5$&$1$&$0.8$&$1.5$&$2$&$2.5$&$3$&$4$\\
				$24$&&$2\sqrt{24}-2$&$2\sqrt{24}-3$&$2\sqrt{24}-4$&$2\sqrt{24}-5$&$3.8$&$3$&$2.4286$&$2\sqrt{24}-4$&$2\sqrt{24}-5$&$2\sqrt{24}-6$&$2.8$&$2$&$1.4286$&$2\sqrt{24}-6$&$2\sqrt{24}-7$&$1.8$&$1$&$0.4286$&$2\sqrt{24}-8$&$0.8$&$0$&$0.5714$&$0.2$&$1$&$1.5714$&$2$&$2.5714$&$3.5714$\\
				$25$&&$8$&$7$&$6$&$5$&$4$&$3.1667$&$2.5714$&$6$&$5$&$4$&$3$&$2.1667$&$1.5714$&$4$&$3$&$2$&$1.1667$&$0.5714$&$2$&$1$&$0.1667$&$4286$&$0$&$0.8333$&$1.4286$&$1.8333$&$2.4286$&$3.4286$\\
				$28$&&$2\sqrt{28}-2$&$2\sqrt{28}-3$&$2\sqrt{28}-4$&$2\sqrt{28}-5$&$2\sqrt{28}-6$&$3.6667$&$3$&$2\sqrt{28}-4$&$2\sqrt{28}-5$&$2\sqrt{28}-6$&$2\sqrt{28}-7$&$2.6667$&$2$&$2\sqrt{28}-6$&$2\sqrt{28}-7$&$2\sqrt{28}-8$&$1.6667$&$1$&$2\sqrt{28}-8$&$2\sqrt{28}-9$&$0.6667$&$0$&$2\sqrt{28}-10$&$0.3333$&$1$&$1.3333$&$2$&$3$\\
				$30$&&$2\sqrt{30}-2$&$2\sqrt{30}-3$&$2\sqrt{30}-4$&$2\sqrt{30}-5$&$2\sqrt{30}-6$&$4$&$3.2857$&$2\sqrt{30}-4$&$2\sqrt{30}-5$&$2\sqrt{30}-6$&$2\sqrt{30}-7$&$3$&$2.2857$&$2\sqrt{30}-6$&$2\sqrt{30}-7$&$2\sqrt{30}-8$&$2$&$1.2857$&$2\sqrt{30}-8$&$2\sqrt{30}-9$&$1$&$0.2857$&$2\sqrt{30}-10$&$0$&$0.7143$&$1$&$1.7143$&$2.7143$\\
				$35$&&$2\sqrt{35}-2$&$2\sqrt{35}-3$&$2\sqrt{35}-4$&$2\sqrt{35}-5$&$2\sqrt{35}-6$&$4.8333$&$4$&$2\sqrt{35}-4$&$2\sqrt{35}-5$&$2\sqrt{35}-6$&$2\sqrt{35}-7$&$3.8333$&$3$&$2\sqrt{35}-6$&$2\sqrt{35}-7$&$2\sqrt{35}-8$&$2.8333$&$2$&$2\sqrt{35}-8$&$2\sqrt{35}-9$&$1.8333$&$1$&$2\sqrt{35}-10$&$0.8333$&$0$&$0.1667$&$1$&$2$\\
				$36$&&$10$&$9$&$8$&$7$&$6$&$5$&$4.1429$&$8$&$7$&$6$&$5$&$4$&$3.1429$&$6$&$5$&$4$&$3$&$2.1429$&$4$&$3$&$2$&$1.1429$&$2$&$1$&$0.1429$&$0$&$0.8571$&$1.8571$\\
				$42$&&$2\sqrt{42}-2$&$2\sqrt{42}-3$&$2\sqrt{42}-4$&$2\sqrt{42}-5$&$2\sqrt{42}-6$&$2\sqrt{42}-7$&$5$&$2\sqrt{42}-4$&$2\sqrt{42}-5$&$2\sqrt{42}-6$&$2\sqrt{42}-7$&$2\sqrt{42}-8$&$4$&$2\sqrt{42}-6$&$2\sqrt{42}-7$&$2\sqrt{42}-8$&$2\sqrt{42}-9$&$3$&$2\sqrt{42}-8$&$2\sqrt{42}-9$&$2\sqrt{42}-10$&$2$&$2\sqrt{42}-10$&$2\sqrt{42}-11$&$1$&$2\sqrt{42}-12$&$0$&$1$\\
				$49$&&$12$&$11$&$10$&$9$&$8$&$7$&$6$&$10$&$9$&$8$&$7$&$6$&$5$&$8$&$7$&$6$&$5$&$4$&$6$&$5$&$4$&$3$&$4$&$3$&$2$&$2$&$1$&$0$\\
				\bottomrule				
			\end{tabular}
		\end{adjustbox}
		\par\smallskip \small \textit{These values serve as the coefficients in problem \eqref{optimization_7} for computing CI$_{a_{BW}}(n)$ when $a_{BW}=7$.}
	\end{sidewaystable}
	Thus, problem \eqref{optimization_7} takes the form 
	
	{\small\begin{align}
			&\max z \nonumber\\
			&\text{sub to: } \nonumber\\
			&0.5n_{1,2} + 0.6667n_{1,3} + 0.75n_{1,4} + 0.8n_{1,5} + 0.8333n_{1,6} + 0.8571n_{1,7} + 1.5n_{2,2}  \nonumber\\
			&\quad + 1.6667n_{2,3}+ 1.75n_{2,4} + 1.8n_{2,5}+ 1.8333n_{2,6}+ 1.8571n_{2,7} + 2.6667n_{3,3}+ 2.75n_{3,4}   \nonumber\\
			&\quad + 2.8n_{3,5}+ 2.8333n_{3,6}+ 2.8571n_{3,7} + 3.75n_{4,4} + 3.8n_{4,5}+ 3.8333n_{4,6}+ 3.8571n_{4,7}    \nonumber\\
			&\quad + 4.8n_{5,5}+ 4.8333n_{5,6}+ 4.8571n_{5,7}+ 5.8333n_{6,6}+ 5.8571n_{6,7} + 6.8571n_{7,7} +6\geq z, \nonumber\\
			&(2\sqrt{2}-2)n_{1,1} + 0.3333n_{1,3} + 0.5n_{1,4} + 0.6n_{1,5} + 0.6667n_{1,6} + 0.7143n_{1,7} + n_{2,2}  \nonumber\\
			&\quad + 1.3333n_{2,3}+ 1.5n_{2,4} + 1.6n_{2,5} + 1.6667n_{2,6} + 1.7143n_{2,7} + 2.3333n_{3,3} + 2.5n_{3,4}  \nonumber\\
			&\quad + 2.6n_{3,5} + 2.6667n_{3,6}+ 4.7143n_{5,7}+ 2.7143n_{3,7} + 3.5n_{4,4} + 3.6n_{4,5} + 3.6667n_{4,6}  \nonumber\\
			&\quad+ 3.7143n_{4,7} + 4.6n_{5,5} + 4.6667n_{5,6}+ 5.6667n_{6,6} + 5.7143n_{6,7} + 6.7143n_{7,7}+5 \geq z, \nonumber\\
			&(2\sqrt{3}-2)n_{1,1} + 0.5n_{1,2} + 0.25n_{1,4} + 0.4n_{1,5} + 0.5n_{1,6} + 0.5714n_{1,7}+ 0.5n_{2,2} + n_{2,3} \nonumber\\
			&\quad  + 1.25n_{2,4} + 1.4n_{2,5} + 1.5n_{2,6} + 1.5714n_{2,7} + 2n_{3,3}+ 2.25n_{3,4} + 2.4n_{3,5} + 2.5n_{3,6} \nonumber\\
			&\quad  + 2.5714n_{3,7} + 3.25n_{4,4} + 3.4n_{4,5} + 3.5n_{4,6} + 3.5714n_{4,7} + 4.4n_{5,5} + 4.5n_{5,6} \nonumber\\
			&\quad + 4.5714n_{5,7}+ 5.5n_{6,6} + 5.5714n_{6,7} + 6.5714n_{7,7} +4\geq z, \nonumber\\
			&2n_{1,1} + n_{1,2} + 0.3333n_{1,3} + 0.2n_{1,5} + 0.3333n_{1,6} + 0.4286n_{1,7} + + 0.6667n_{2,3} + n_{2,4}  \nonumber\\
			&\quad + 1.2n_{2,5}+ 1.3333n_{2,6} + 1.4286n_{2,7} + 1.6667n_{3,3} + 2n_{3,4} + 2.2n_{3,5} + 2.3333n_{3,6}  \nonumber\\
			&\quad + 2.4286n_{3,7} + 3n_{4,4}+ 3.2n_{4,5} + 3.3333n_{4,6} + 3.4286n_{4,7} + 4.2n_{5,5} + 4.3333n_{5,6}  \nonumber\\
			&\quad + 4.4286n_{5,7} + 5.3333n_{6,6}+ 5.4286n_{6,7} + 6.4286n_{7,7}+3 \geq z, \nonumber\\
			&(2\sqrt{5}-2)n_{1,1} + (2\sqrt{5}-3)n_{1,2} + 0.6667n_{1,3} + 0.25n_{1,4} + 0.1667n_{1,6} + 0.2857n_{1,7} \nonumber\\
			&\quad + (2\sqrt{5}-4)n_{2,2}+ 0.3333n_{2,3} + 0.75n_{2,4} + n_{2,5} + 1.1667n_{2,6} + 1.2857n_{2,7} \nonumber \\
			&\quad +1.3333n_{3,3}+ 1.75n_{3,4} + 2n_{3,5}+ 2.1667n_{3,6} + 2.2857n_{3,7} + 2.75n_{4,4} + 3n_{4,5}  \nonumber\\
			&\quad + 3.1667n_{4,6} + 3.2857n_{4,7}+ 4n_{5,5} + 4.1667n_{5,6} + 4.2857n_{5,7} + 5.1667n_{6,6}\nonumber\\
			&\quad + 5.2857n_{6,7} + 6.2857n_{7,7}+2 \geq z, \nonumber\\
			&(2\sqrt{6}-2)n_{1,1} + (2\sqrt{6}-3)n_{1,2} + n_{1,3} + 0.5n_{1,4} + 0.2n_{1,5}  + 0.1429n_{1,7} + (2\sqrt{6}-4)n_{2,2}\nonumber\\
			&\quad + 0.5n_{2,4} + 0.8n_{2,5} + n_{2,6} + 1.1429n_{2,7} + n_{3,3} + 1.5n_{3,4} + 1.8n_{3,5} + 2n_{3,6} \nonumber\\
			&\quad + 2.1429n_{3,7}+ 2.5n_{4,4} + 2.8n_{4,5} + 3n_{4,6} + 3.1429n_{4,7} + 3.8n_{5,5} + 4n_{5,6} + 4.1429n_{5,7}  \nonumber\\
			&\quad + 5n_{6,6}+ 5.1429n_{6,7} + 6.1429n_{7,7}+1 \geq z, \nonumber\\
			&(2\sqrt{7}-2)n_{1,1} + (2\sqrt{7}-3)n_{1,2} + 1.3333n_{1,3} + 0.75n_{1,4} + 0.4n_{1,5} + 0.1667n_{1,6} \nonumber\\
			&\quad + (2\sqrt{7}-4)n_{2,2}+ 0.3333n_{2,3} + 0.25n_{2,4} + 0.6n_{2,5} + 0.8333n_{2,6} + n_{2,7} + 0.6667n_{3,3}  \nonumber\\
			&\quad + 1.25n_{3,4} + 1.6n_{3,5}+ 1.8333n_{3,6} + 2n_{3,7} + 2.25n_{4,4} + 2.6n_{4,5} + 2.8333n_{4,6} + 3n_{4,7}  \nonumber\\
			&\quad + 3.6n_{5,5} + 3.8333n_{5,6}+ 4n_{5,7}+ 4.8333n_{6,6} + 5n_{6,7} + 6n_{7,7} \geq z, \nonumber\\
			&(2\sqrt{8}-2)n_{1,1} + (2\sqrt{8}-3)n_{1,2} + 1.6667n_{1,3} + n_{1,4} + 0.6n_{1,5} + 0.3333n_{1,6} + 0.1429n_{1,7} \nonumber\\
			&\quad + (2\sqrt{8}-4)n_{2,2} + 0.6667n_{2,3} + 0.4n_{2,5} + 0.6667n_{2,6} + 0.8571n_{2,7} + 0.3333n_{3,3} + n_{3,4} \nonumber\\
			&\quad + 1.4n_{3,5}+ 1.6667n_{3,6} + 1.8571n_{3,7} + 2n_{4,4} + 2.4n_{4,5} + 2.6667n_{4,6} + 2.8571n_{4,7}  \nonumber\\
			&\quad + 3.4n_{5,5}+ 3.6667n_{5,6}+ 3.8571n_{5,7} + 4.6667n_{6,6} + 4.8571n_{6,7} + 5.8571n_{7,7}+1 \geq z, \nonumber\\
			&4n_{1,1} + 3n_{1,2} + 2n_{1,3} + 1.25n_{1,4} + 0.8n_{1,5} + 0.5n_{1,6} + 0.2857n_{1,7} + 2n_{2,2} + n_{2,3} \nonumber\\
			&\quad + 0.25n_{2,4}+ 0.2n_{2,5} + 0.5n_{2,6} + 0.7143n_{2,7} + 0.75n_{3,4} + 1.2n_{3,5} + 1.5n_{3,6} \nonumber\\
			&\quad + 1.7143n_{3,7} + 1.75n_{4,4}+ 2.2n_{4,5} + 2.5n_{4,6} + 2.7143n_{4,7} + 3.2n_{5,5} + 3.5n_{5,6} \nonumber\\
			&\quad + 3.7143n_{5,7} + 4.5n_{6,6} + 4.7143n_{6,7}+ 5.7143n_{7,7}+2 \geq z, \nonumber\\
			&(2\sqrt{10}-2)n_{1,1} + (2\sqrt{10}-3)n_{1,2} + (2\sqrt{10}-4)n_{1,3} + 1.5n_{1,4} + n_{1,5} + 0.6667n_{1,6} \nonumber\\
			&\quad + 0.4286n_{1,7}+ (2\sqrt{10}-4)n_{2,2} + (2\sqrt{10}-5)n_{2,3} + 0.5n_{2,4} + 0.3333n_{2,6} + 0.5714n_{2,7} \nonumber\\
			&\quad + (2\sqrt{10}-6)n_{3,3}+ 0.5n_{3,4} + n_{3,5} + 1.3333n_{3,6} + 1.5714n_{3,7} + 1.5n_{4,4} + 2n_{4,5} \nonumber\\
			&\quad + 2.3333n_{4,6} + 2.5714n_{4,7}+ 3n_{5,5} + 3.3333n_{5,6} + 3.5714n_{5,7} + 4.3333n_{6,6} \nonumber\\
			&\quad + 4.5714n_{6,7}+ 5.5714n_{7,7}+3 \geq z, \nonumber\\
			&(2\sqrt{12}-2)n_{1,1} + (2\sqrt{12}-3)n_{1,2} + (2\sqrt{12}-4)n_{1,3} + 2n_{1,4} + 1.4n_{1,5} + n_{1,6} + 0.7143n_{1,7}\nonumber\\
			&\quad + (2\sqrt{12}-4)n_{2,2} + (2\sqrt{12}-5)n_{2,3} + n_{2,4} + 0.4n_{2,5} + 0.2857n_{2,7} + (2\sqrt{12}-6)n_{3,3} \nonumber\\
			&\quad + 0.6n_{3,5} + n_{3,6} + 1.2857n_{3,7} + n_{4,4} + 1.6n_{4,5} + 2n_{4,6} + 2.2857n_{4,7} + 2.6n_{5,5} + 3n_{5,6}\nonumber\\
			&\quad + 3.2857n_{5,7} + 4n_{6,6} + 4.2857n_{6,7} + 5.2857n_{7,7} +5\geq z, \nonumber\\
			&(2\sqrt{14}-2)n_{1,1} + (2\sqrt{14}-3)n_{1,2} + (2\sqrt{14}-4)n_{1,3} + 2.5n_{1,4} + 1.8n_{1,5} + 1.3333n_{1,6} \nonumber\\
			&\quad + n_{1,7}+ (2\sqrt{14}-4)n_{2,2} + (2\sqrt{14}-5)n_{2,3} + 1.5n_{2,4} + 0.8n_{2,5} + 0.3333n_{2,6} \nonumber\\
			&\quad + (2\sqrt{14}-6)n_{3,3}+ 0.5n_{3,4} + 0.2n_{3,5} + 0.6667n_{3,6} + n_{3,7} + 0.5n_{4,4} + 1.2n_{4,5} \nonumber\\
			&\quad + 1.6667n_{4,6} + 2n_{4,7} + 2.2n_{5,5}+ 2.6667n_{5,6} + 3n_{5,7} + 3.6667n_{6,6} + 4n_{6,7} + 5n_{7,7}\nonumber\\
			&\quad +7 \geq z, \nonumber\\
			&(2\sqrt{15}-2)n_{1,1} + (2\sqrt{15}-3)n_{1,2} + (2\sqrt{15}-4)n_{1,3} + 2.75n_{1,4} + 2n_{1,5} + 1.5n_{1,6} \nonumber\\
			&\quad + 1.1429n_{1,7}+ (2\sqrt{15}-4)n_{2,2} + (2\sqrt{15}-5)n_{2,3} + 1.75n_{2,4} + n_{2,5} + 0.5n_{2,6} \nonumber\\
			&\quad + 0.1429n_{2,7}+ (2\sqrt{15}-6)n_{3,3} + 0.75n_{3,4} + 0.5n_{3,6} + 0.8571n_{3,7} + 0.25n_{4,4} + n_{4,5}  \nonumber\\
			&\quad + 1.5n_{4,6}+ 1.8571n_{4,7}+ 2n_{5,5} + 2.5n_{5,6} + 2.8571n_{5,7} + 3.5n_{6,6} + 3.8571n_{6,7} \nonumber\\
			&\quad + 4.8571n_{7,7}+8 \geq z, \nonumber\\
			&6n_{1,1} + 5n_{1,2} + 4n_{1,3} + 3n_{1,4} + 2.2n_{1,5} + 1.6667n_{1,6} + 1.2857n_{1,7} + 4n_{2,2} + 3n_{2,3} \nonumber\\
			&\quad + 2n_{2,4}+ 1.2n_{2,5} + 0.6667n_{2,6} + 0.2857n_{2,7} + 2n_{3,3} + n_{3,4} + 0.2n_{3,5} + 0.3333n_{3,6} \nonumber\\
			&\quad + 0.7143n_{3,7}+ 0.8n_{4,5} + 1.3333n_{4,6} + 1.7143n_{4,7} + 1.8n_{5,5} + 2.3333n_{5,6} + 2.7143n_{5,7} \nonumber\\
			&\quad + 3.3333n_{6,6}+ 3.7143n_{6,7} + 4.7143n_{7,7}+9 \geq z, \nonumber\\
			&(2\sqrt{18}-2)n_{1,1} + (2\sqrt{18}-3)n_{1,2} + (2\sqrt{18}-4)n_{1,3} + (2\sqrt{18}-5)n_{1,4} + 2.6n_{1,5} + 2n_{1,6} \nonumber\\
			&\quad + 1.5714n_{1,7} + (2\sqrt{18}-4)n_{2,2} + (2\sqrt{18}-5)n_{2,3} + (2\sqrt{18}-6)n_{2,4} + 1.6n_{2,5} + n_{2,6}  \nonumber\\
			&\quad + 0.5714n_{2,7} + (2\sqrt{18}-6)n_{3,3} + (2\sqrt{18}-7)n_{3,4} + 0.6n_{3,5} + 0.4286n_{3,7}  \nonumber\\
			&\quad + (2\sqrt{18}-8)n_{4,4}+ 0.4n_{4,5} + n_{4,6} + 1.4286n_{4,7} + 1.4n_{5,5} + 2n_{5,6} + 2.4286n_{5,7} \nonumber\\
			&\quad + 3n_{6,6} + 3.4286n_{6,7}+ 4.4286n_{7,7}+11 \geq z, \nonumber\\
			&(2\sqrt{20}-2)n_{1,1} + (2\sqrt{20}-3)n_{1,2} + (2\sqrt{20}-4)n_{1,3} + (2\sqrt{20}-5)n_{1,4} + 3n_{1,5} \nonumber\\
			&\quad + 2.3333n_{1,6}+ 1.8571n_{1,7} + (2\sqrt{20}-4)n_{2,2} + (2\sqrt{20}-5)n_{2,3} + (2\sqrt{20}-6)n_{2,4} \nonumber\\
			&\quad + 2n_{2,5} + 1.3333n_{2,6}+ 0.8571n_{2,7} + (2\sqrt{20}-6)n_{3,3} + (2\sqrt{20}-7)n_{3,4} + n_{3,5} \nonumber\\
			&\quad + 0.3333n_{3,6} + 0.1429n_{3,7}+ (2\sqrt{20}-8)n_{4,4} + 0.6667n_{4,6} + 1.1429n_{4,7} + n_{5,5}   \nonumber\\
			&\quad + 1.6667n_{5,6}+ 2.1429n_{5,7} + 2.6667n_{6,6}+ 3.1429n_{6,7} + 4.1429n_{7,7}+13 \geq z, \nonumber\\
			&(2\sqrt{21}-2)n_{1,1} + (2\sqrt{21}-3)n_{1,2} + (2\sqrt{21}-4)n_{1,3} + (2\sqrt{21}-5)n_{1,4} + 3.2n_{1,5} + 2.5n_{1,6} \nonumber\\
			&\quad + 2n_{1,7}+ (2\sqrt{21}-4)n_{2,2} + (2\sqrt{21}-5)n_{2,3} + (2\sqrt{21}-6)n_{2,4} + 2.2n_{2,5} + 1.5n_{2,6}  \nonumber\\
			&\quad + n_{2,7}+ (2\sqrt{21}-6)n_{3,3} + (2\sqrt{21}-7)n_{3,4} + 1.2n_{3,5} + 0.5n_{3,6} + (2\sqrt{21}-8)n_{4,4}  \nonumber\\
			&\quad + 0.2n_{4,5}+ 0.5n_{4,6}+ n_{4,7} + 0.8n_{5,5} + 1.5n_{5,6} + 2n_{5,7} + 2.5n_{6,6} + 3n_{6,7} + 4n_{7,7}\nonumber\\
			&\quad +14 \geq z, \nonumber\\
			&(2\sqrt{24}-2)n_{1,1} + (2\sqrt{24}-3)n_{1,2} + (2\sqrt{24}-4)n_{1,3} + (2\sqrt{24}-5)n_{1,4} + 3.8n_{1,5} + 3n_{1,6}\nonumber\\
			&\quad + 2.4286n_{1,7} + (2\sqrt{24}-4)n_{2,2} + (2\sqrt{24}-5)n_{2,3} + (2\sqrt{24}-6)n_{2,4} + 2.8n_{2,5} + 2n_{2,6}\nonumber\\
			&\quad + 1.4286n_{2,7} + (2\sqrt{24}-6)n_{3,3} + (2\sqrt{24}-7)n_{3,4} + 1.8n_{3,5} + n_{3,6} + 0.4286n_{3,7} \nonumber\\
			&\quad + (2\sqrt{24}-8)n_{4,4} + 0.8n_{4,5} + 0.5714n_{4,7} + 0.2n_{5,5} + n_{5,6} + 1.5714n_{5,7} + 2n_{6,6} \nonumber\\
			&\quad + 2.5714n_{6,7} + 3.5714n_{7,7}+17 \geq z, \nonumber\\
			&8n_{1,1} + 7n_{1,2} + 6n_{1,3} + 5n_{1,4} + 4n_{1,5} + 3.1667n_{1,6} + 2.5714n_{1,7} + 6n_{2,2} + 5n_{2,3} + 4n_{2,4}\nonumber\\
			&\quad + 3n_{2,5} + 2.1667n_{2,6} + 1.5714n_{2,7} + 4n_{3,3} + 3n_{3,4} + 2n_{3,5} + 1.1667n_{3,6} + 0.5714n_{3,7}\nonumber\\
			&\quad + 2n_{4,4}+ n_{4,5} + 0.1667n_{4,6} + 0.4286n_{4,7} + 0.8333n_{5,6} + 1.4286n_{5,7} + 1.8333n_{6,6} \nonumber\\
			&\quad + 2.4286n_{6,7} + 3.4286n_{7,7}+18 \geq z, \nonumber\\
			&(2\sqrt{28}-2)n_{1,1} + (2\sqrt{28}-3)n_{1,2} + (2\sqrt{28}-4)n_{1,3} + (2\sqrt{28}-5)n_{1,4} + (2\sqrt{28}-6)n_{1,5}\nonumber\\
			&\quad + 3.6667n_{1,6} + 3n_{1,7} + (2\sqrt{28}-4)n_{2,2} + (2\sqrt{28}-5)n_{2,3} + (2\sqrt{28}-6)n_{2,4}\nonumber\\
			&\quad + (2\sqrt{28}-7)n_{2,5} + 2.6667n_{2,6} + 2n_{2,7} + (2\sqrt{28}-6)n_{3,3} + (2\sqrt{28}-7)n_{3,4} \nonumber\\
			&\quad + (2\sqrt{28}-8)n_{3,5} + 1.6667n_{3,6} + n_{3,7} + (2\sqrt{28}-8)n_{4,4} + (2\sqrt{28}-9)n_{4,5} \nonumber\\
			&\quad + 0.6667n_{4,6}+ (2\sqrt{28}-10)n_{5,5} + 0.3333n_{5,6} + n_{5,7} + 1.3333n_{6,6} + 2n_{6,7} + 3n_{7,7}\nonumber\\
			&\quad +21 \geq z, \nonumber\\
			&(2\sqrt{30}-2)n_{1,1} + (2\sqrt{30}-3)n_{1,2} + (2\sqrt{30}-4)n_{1,3} + (2\sqrt{30}-5)n_{1,4} + (2\sqrt{30}-6)n_{1,5}\nonumber\\
			&\quad + 4n_{1,6} + 3.2857n_{1,7} + (2\sqrt{30}-4)n_{2,2} + (2\sqrt{30}-5)n_{2,3} + (2\sqrt{30}-6)n_{2,4} \nonumber\\
			&\quad + (2\sqrt{30}-7)n_{2,5}+ 3n_{2,6} + 2.2857n_{2,7} + (2\sqrt{30}-6)n_{3,3} + (2\sqrt{30}-7)n_{3,4} \nonumber\\
			&\quad + (2\sqrt{30}-8)n_{3,5} + 2n_{3,6}+ 1.2857n_{3,7} + (2\sqrt{30}-8)n_{4,4} + (2\sqrt{30}-9)n_{4,5} +n_{4,6} \nonumber\\
			&\quad + 0.2857n_{4,7} + (2\sqrt{30}-10)n_{5,5}+ 0.7143n_{5,7} + n_{6,6} + 1.7143n_{6,7} + 2.7143n_{7,7}\nonumber\\
			&\quad +23 \geq z, \nonumber\\
			&(2\sqrt{35}-2)n_{1,1} + (2\sqrt{35}-3)n_{1,2} + (2\sqrt{35}-4)n_{1,3} + (2\sqrt{35}-5)n_{1,4} + (2\sqrt{35}-6)n_{1,5}\nonumber\\
			&\quad + 4.8333n_{1,6} + 4n_{1,7} + (2\sqrt{35}-4)n_{2,2} + (2\sqrt{35}-5)n_{2,3} + (2\sqrt{35}-6)n_{2,4}\nonumber\\
			&\quad + (2\sqrt{35}-7)n_{2,5} + 3.8333n_{2,6} + 3n_{2,7} + (2\sqrt{35}-6)n_{3,3} + (2\sqrt{35}-7)n_{3,4} \nonumber\\
			&\quad + (2\sqrt{35}-8)n_{3,5} + 2.8333n_{3,6} + 2n_{3,7} + (2\sqrt{35}-8)n_{4,4} + (2\sqrt{35}-9)n_{4,5}\nonumber\\
			&\quad + 1.8333n_{4,6}+ n_{4,7} + (2\sqrt{35}-10)n_{5,5} + 0.8333n_{5,6} + 0.1667n_{6,6} + 1n_{6,7} + 2n_{7,7}\nonumber\\
			&quad +28 \geq z, \nonumber\\
			&10n_{1,1} + 9n_{1,2} + 8n_{1,3} + 7n_{1,4} + 6n_{1,5} + 5n_{1,6} + 4.1429n_{1,7} + 8n_{2,2} + 7n_{2,3} + 6n_{2,4}\nonumber\\
			&\quad + 5n_{2,5} + 4n_{2,6} + 3.1429n_{2,7} + 6n_{3,3} + 5n_{3,4} + 4n_{3,5} + 3n_{3,6} + 2.1429n_{3,7} + 4n_{4,4}\nonumber\\
			&\quad + 3n_{4,5} + 2n_{4,6} + 1.1429n_{4,7} + 2n_{5,5} + n_{5,6} + 0.1429n_{5,7} + 0.8571n_{6,7} + 1.8571n_{7,7}\nonumber\\
			&\quad +29 \geq z, \nonumber\\
			&(2\sqrt{42}-2)n_{1,1} + (2\sqrt{42}-3)n_{1,2} + (2\sqrt{42}-4)n_{1,3} + (2\sqrt{42}-5)n_{1,4} + (2\sqrt{42}-6)n_{1,5}\nonumber\\
			&\quad + (2\sqrt{42}-7)n_{1,6} + 5n_{1,7} + (2\sqrt{42}-4)n_{2,2} + (2\sqrt{42}-5)n_{2,3} + (2\sqrt{42}-6)n_{2,4}\nonumber\\
			&\quad + (2\sqrt{42}-7)n_{2,5} + (2\sqrt{42}-8)n_{2,6} + 4n_{2,7} + (2\sqrt{42}-6)n_{3,3} + (2\sqrt{42}-7)n_{3,4}\nonumber\\
			&\quad + (2\sqrt{42}-8)n_{3,5} + (2\sqrt{42}-9)n_{3,6} + 3n_{3,7} + (2\sqrt{42}-8)n_{4,4} + (2\sqrt{42}-9)n_{4,5}\nonumber\\
			&\quad + (2\sqrt{42}-10)n_{4,6} + 2n_{4,7} + (2\sqrt{42}-10)n_{5,5} + (2\sqrt{42}-11)n_{5,6} + n_{5,7} \nonumber\\
			&\quad + (2\sqrt{42}-12)n_{6,6} + n_{7,7}+35 \geq z, \nonumber\\
			&12n_{1,1} + 11n_{1,2} + 10n_{1,3} + 9n_{1,4} + 8n_{1,5} + 7n_{1,6} + 6n_{1,7} + 10n_{2,2} + 9n_{2,3} + 8n_{2,4}\nonumber\\
			&\quad + 7n_{2,5} + 6n_{2,6} + 5n_{2,7} + 8n_{3,3} + 7n_{3,4} + 6n_{3,5} + 5n_{3,6} + 4n_{3,7} + 6n_{4,4} + 5n_{4,5}\nonumber\\
			&\quad + 4n_{4,6} + 3n_{4,7} + 4n_{5,5} + 3n_{5,6} + 2n_{5,7} + 2n_{6,6} + n_{6,7} +42 \geq z, \nonumber\\
			\label{optimization_9}
			&\displaystyle\sum_{\substack{a,b=1\\a\leq b}}^{7}n_{a,b}+2=15,\quad   n_{a,b}\in \mathbb{N}\cup\{0\} \text{ for all }a,b. 
	\end{align}}
	An optimal solution of this problem is $n_{1,1}^*=3$, $n_{7,7}^*=10$, $n_{a,b}^*=0$ for all other $a,b$, and $z^*=69.8745$. Thus, CI$_7(15)=69.8745$, and the best-to-other vector $A_{B}=\left(1,1,1,1,7,7,7,7,7,7,7,7,7,7,7\right)$ and the other-to-worst vector $A_{W}=(7,1,1,1,7,7,7,7,7,7,7,7,7,7,1)^T$ with $c_1$ as the best and $c_{15}$ as the worst criterion forms a PCS with $\epsilon^*=69.8745$.\\\\
	The values of CI$_{a_{BW}}(n)$ for $a_{BW}=2,3,\ldots,9$ and $n=3,4,\ldots,15$ are given in Table \ref{3ci}.\\
	 
	\begin{table}[t!]
		\caption{The values of CI$_{a_{BW}}(n)$}\label{3ci}
		\centering		
		\begin{tabular}{@{}cccccccccc@{}}
			\toprule
			\multirow{2}{*}{$n\downarrow$}&\phantom{}&\multicolumn{8}{c}{$a_{BW}\rightarrow$}\\
			\cmidrule{3-10}
			&&$2$&$3$&$4$&$5$&$6$&$7$&$8$&$9$\\
			\midrule
			$3$&&$1$&$2$&$3$&$4$&$5$&$6$&$7$&$8$\\
			$4$&&$2$&$4$&$6$&$8$&$10$&$12$&$14$&$16$\\
			$5$&&$2.8284$&$6$&$9$&$12$&$15$&$18$&$21$&$24$\\
			$6$&&$3.8284$&$7.4641$&$12$&$16$&$20$&$24$&$28$&$32$\\
			$7$&&$4.6569$&$9.4641$&$14$&$20$&$25$&$30$&$35$&$40$\\
			$8$&&$5.6569$&$10.9282$&$17$&$22.4721$&$30$&$36$&$42$&$48$\\
			$9$&&$6.4853$&$12.9282$&$19$&$26.4721$&$32.899$&$42$&$49$&$56$\\
			$10$&&$7.4853$&$14.3923$&$22$&$28.9443$&$37.899$&$45.2915$&$56$&$64$\\
			$11$&&$8.3137$&$16.3923$&$24$&$32.9443$&$40.798$&$51.2915$&$59.6569$&$71.9996$\\
			$12$&&$9.3137$&$18$&$27$&$36$&$45.798$&$54.5830$&$66.6569$&$76$\\
			$13$&&$10$&$19.8564$&$30$&$39.4164$&$50$&$60.5830$&$70.3137$&$84$\\
			$14$&&$11.1421$&$21.8564$&$32$&$43.4164$&$53.6969$&$66$&$77.3137$&$88$\\
			$15$&&$12$&$23.3205$&$35$&$45.8885$&$58.6969$&$69.8745$&$84$&$96$\\			
			\bottomrule		
		\end{tabular}
		\par\smallskip \small \textit{The CI increases monotonically with both parameters $n$ and $a_{BW}$, reflecting the higher maximum possible TD in larger or more dispersed PCSs.}
	\end{table}
	 
	\hspace{-0.7cm}
	Some of the key points regarding the CI and the CR are as follows.
	\begin{enumerate}
		\item For the nonlinear BWM, the CI can be calculated by simply finding the roots of linear and quadratic equations, while the calculation is more complex in the case of the taxicab BWM, which requires the formulation and solution of a mixed-integer linear optimization problem.
		\item The computational complexity of optimization model \eqref{optimization_7} is highly sensitive to the parameter $a_{BW}$, as a larger value increases the number of variables and constraints. In contrast, the complexity remains unaffected by the number of criteria $n$, as an increase in $n$ does not increase the number of variables or constraints, making the model scalable for large-scale MCDM problems. This distinction is clearly demonstrated in the examples above.
		\item The nonlinear BWM's reliance on the maximum deviation means its CR captures only the most extreme violation, disregarding information from other pairwise comparisons. Consequently, the CI for the nonlinear BWM is independent of $n$ and depends solely on $a_{BW}$. In contrast, the taxicab BWM's use of total deviation causes its CR to aggregate all pairwise comparison deviations, resulting in a more comprehensive assessment. This fundamental difference in objective functions also explains the computational complexity gap between the two models' CIs.
		\item A lower value of the CR indicates a higher level of consistency in the PCS, with a value of CR $= 0$ representing a consistent PCS.
		\item The overall CR provides a global measure of inconsistency for the entire PCS. To estimate the local inconsistency at the $i^{th}$ criterion, where $i\in D$, we can consider a three-criteria PCS involving the best criterion ($B$), the worst criterion ($W$), and the criterion in question ($i$). The CR calculated for this subsystem, denoted as CR$_i$, serves as a measure of the local inconsistency at the $i^{th}$ criterion.\\\\
		Note that for this three-criteria PCS, the global minima of $f$ is  $x^* + y^* + z^*$, where $(x^*, y^*, z^*)$ is an optimal modification strategy for $(a_{Bi}, a_{iW}, a_{BW})$. Since for any optimally modified $(a_{Bi}, a_{iW}, a_{BW})$, $a_{BW}$ remains unchanged, we get $z^* = 0$, and thus the only possible value of $\tilde{a}_{BW}^*$ is $a_{BW}$. So, by Eq. \eqref{optimal_modification}, we get  
		{\small\begin{align}
				\epsilon^*= \begin{cases}
					0\quad\quad\quad\quad\quad\quad\quad\quad\ \ \ \ \   \text{if } a_{Bi}\times a_{iW}=a_{BW},\\
					2\sqrt{a_{BW}}-a_{Bi}-a_{iW}\quad \text{if } a_{Bi}\times a_{iW}< a_{BW} \text{ and } a_{Bi},a_{iW} <\sqrt{a_{BW}},\\
					\bigg|a_{Bi}-\frac{a_{BW}}{a_{iW}}\bigg|\quad\quad\quad\quad\ \   \text{if } (a_{Bi}\times a_{iW}<a_{BW} \text{ and } a_{Bi}<\sqrt{a_{BW}}\leq a_{iW})\\ \quad\quad\quad\quad\quad\quad\quad\quad\ \ \ \ \ \text{ or } (a_{Bi}\times a_{iW}>a_{BW} \text{ and } a_{Bi}\leq a_{iW}),\\
					\bigg|a_{iW}-\frac{a_{BW}}{a_{Bi}}\bigg|\quad\quad\quad\quad\ \   \text{if } (a_{Bi}\times a_{iW}<a_{BW} \text{ and } a_{iW}<\sqrt{a_{BW}}\leq a_{Bi})\\ \quad\quad\quad\quad\quad\quad\quad\quad\ \ \ \ \   \text{ or } (a_{Bi}\times a_{iW}>a_{BW} \text{ and } a_{iW}<a_{Bi}).
				\end{cases}\nonumber\\
		\end{align}}
		Also observe that, by Table \ref{3ci}, $\text{CI}_{a_{BW}}(3) = a_{BW}-1$. So,
		\begin{equation}\label{cr_i}
			\text{CR}_i=\frac{\epsilon^*}{a_{BW}-1}.
		\end{equation}
		This analytical form allows $CR_i$ to be used as an input-based consistency indicator.
	\end{enumerate}
	\subsection{Determination of thresholds}
\begin{figure}[th]
		\centering
		\includegraphics[height=6cm,width=13cm]{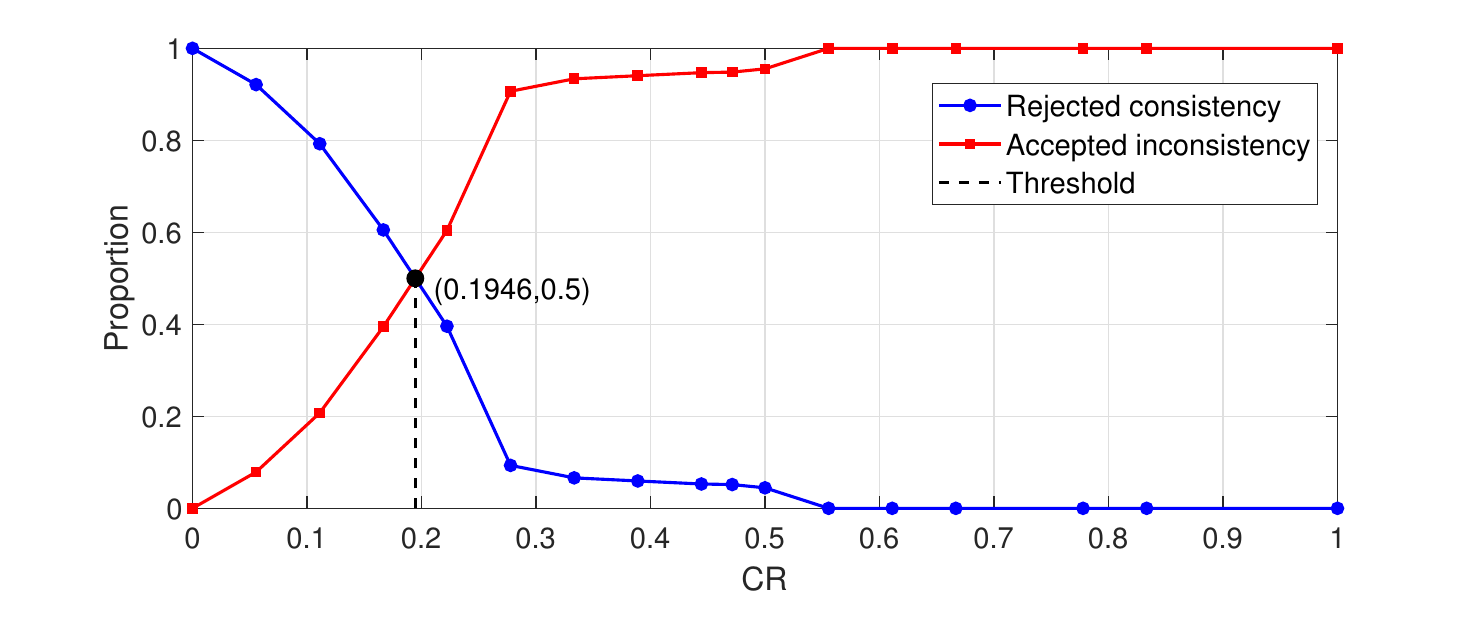}
		\caption{Relative acceptance and rejection proportions ($n=4$, $a_{BW}=4$)}
		\label{n4_4}		
	\end{figure}

\begin{figure}[th]
		\centering
		\includegraphics[height=14cm,width=13cm]{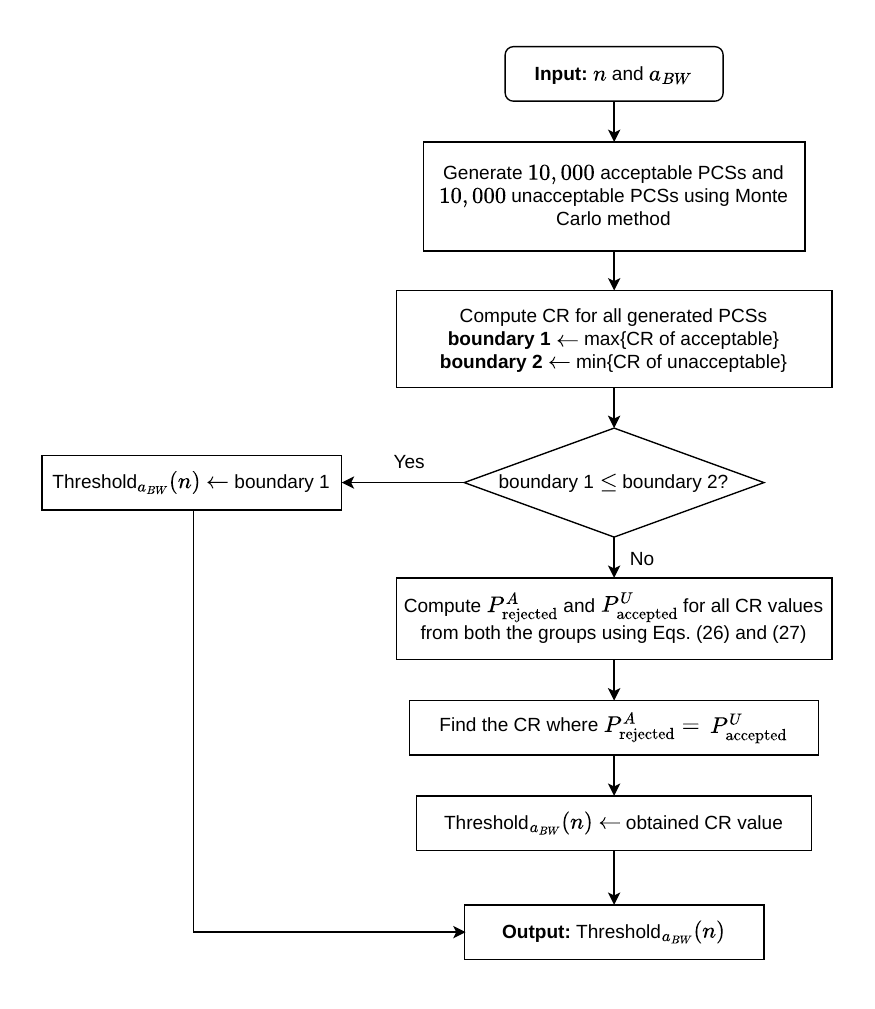}
		\caption{Algorithm for computing $\text{Threshold}_{a_{BW}}(n)$}
		\label{fig:flowchart_threshold}		
	\end{figure}

\begin{figure}[th]
		\centering
		\includegraphics[height=14cm,width=13cm]{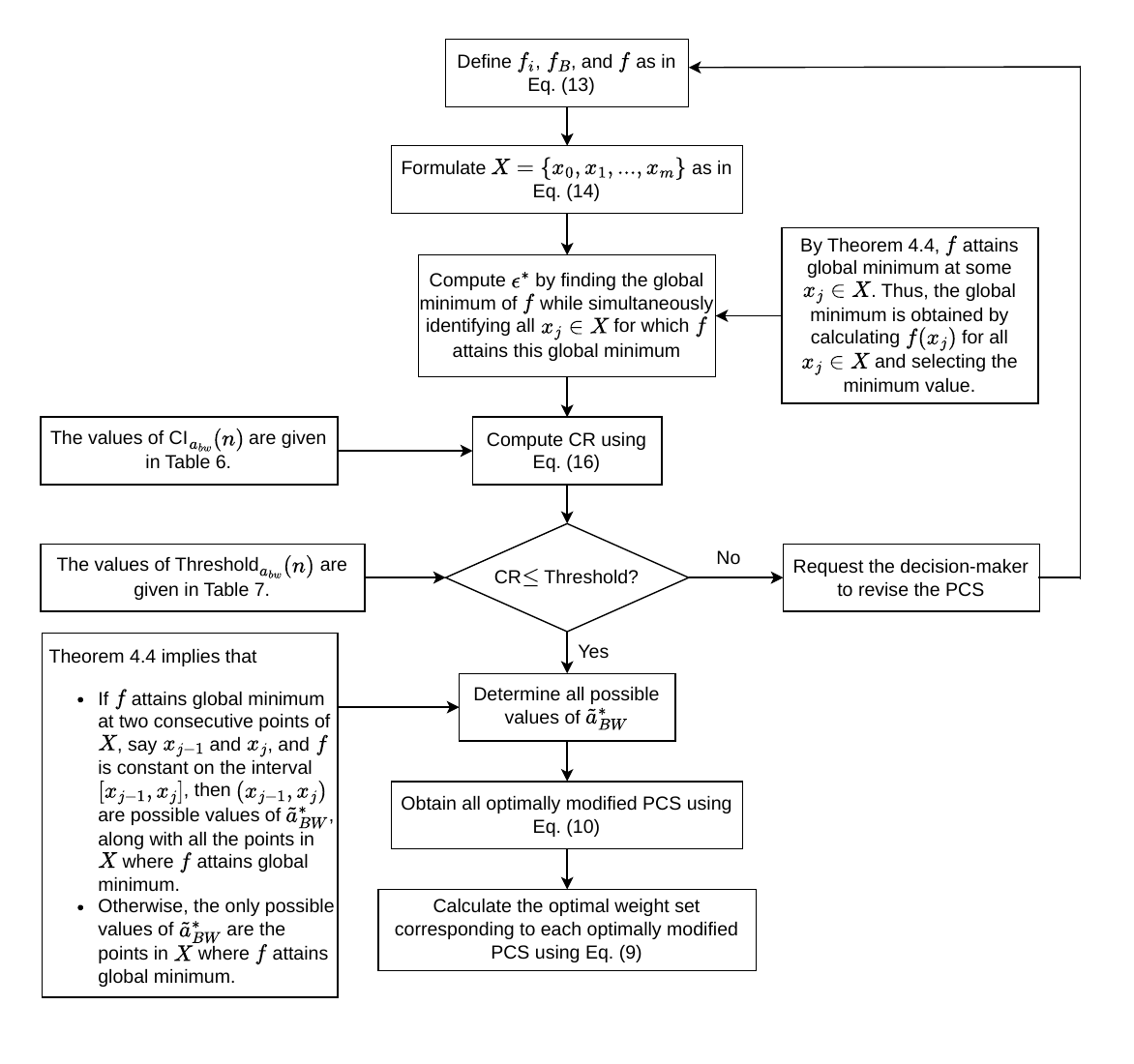}
		\caption{Workflow for weight computation and consistency analysis}
		\label{flowchart}		
	\end{figure}

	In this subsection, we establish the ordinal consistency-based thresholds for the CR to assess the admissibility of a PCS, adopting the methodology developed by Liang et al. \cite{liang2020consistency} in their similar study for the nonlinear BWM.\\\\
	The fundamental philosophy of this approach is to select an ordinal-consistent PCS as admissible, because for such a PCS, the ranking of weights derived solely from $A_B$ and that derived solely from $A_W$ remain the same. The theoretical procedure for determining the threshold for a PCS with $n$ criteria and a given $a_{BW}$, denoted as Threshold$_{a_{BW}}(n)$, is as follows:\\\\
		\textbf{Step 1:} Classification of PCS\\
		Classify all PCSs having $n$ criteria and a given $a_{BW}$ into two groups: ordinal-consistent PCSs (the acceptable group) and ordinal-inconsistent PCSs (the unacceptable group).\\\\
		\textbf{Step 2:} CR for the acceptable group\\
		Compute the CR for every PCS in the acceptable group and denote the maximum CR among them as boundary 1. Thus, a PCS with CR above boundary 1 cannot be ordinal-consistent.\\\\
		\textbf{Step 3:} CR for the unacceptable group\\
		Compute the CR for every PCS in the unacceptable group and denote the minimum CR among them as boundary 2. Thus, a PCS with CR below boundary 2 is always ordinal-consistent.\\\\
		\textbf{Step 4:} Determination of the threshold\\
		If boundary 1 $\leq$ boundary 2, then set Threshold$_{a_{BW}}(n) =$ boundary 1, as it indicates a clear separation between the two groups.\\\\
		If boundary 1 $>$ boundary 2, then there exists an overlapping region between boundary 2 and boundary 1 which contains both ordinal-consistent and ordinal-inconsistent PCSs. In such cases, we select Threshold$_{a_{BW}}(n)$ such that the trade-off between the proportion of false rejections in the acceptable group and the proportion of false acceptances in the unacceptable group is perfectly balanced. This can be achieved using the Cumulative Distribution Function (CDF) of CR, defined as
		\begin{equation}\label{cdf}
			F(\alpha)=\frac{1}{N}\sum_{i=1}^{N}I\{\text{CR}_i\leq \alpha\},
		\end{equation}
		where $N$ is the number of PCSs, $\alpha\in [0,1]$ is the possible threshold, $\text{CR}_i$ is the CR of $i$th PCS, and $I\{\text{CR}_i\leq \alpha\}=\begin{cases}
			1\ \ \text{if } \text{CR}_i\leq \alpha,\\
			0\ \ \text{otherwise}.
		\end{cases}$
		The relative rejected proportion of CRs in the acceptable group ($P^A_{\text{rejected}}$) and the accepted proportion of the CRs in the unacceptable group ($P^U_{\text{accepted}}$) are given by
		\begin{equation}\label{rejected}
			P^A_{\text{rejected}}(\alpha)=\frac{1-F^A(\alpha)}{1-F^A(\alpha)+F^U(\alpha)}
		\end{equation}
		and
		\begin{equation}\label{accepted}
			P^U_{\text{accepted}}(\alpha)=\frac{F^U(\alpha)}{1-F^A(\alpha)+F^U(\alpha)}
		\end{equation}
		respectively, where $F^A$ and $F^U$ are the CDFs of CR for the acceptable group and the unacceptable group respectively.\\\\
		The value of $\alpha$ for which $P^A_{\text{rejected}}(\alpha) = P^U_{\text{accepted}}(\alpha)$ is selected as Threshold$_{a_{BW}}(n)$.\\\\
		Since it is practically infeasible to perform the above procedure for all possible PCSs, we compute the approximate Threshold$_{a_{BW}}(n)$ by randomly generating 10,000 PCSs each for the acceptable and unacceptable groups using the Monte Carlo simulation. After obtaining boundary 1 and boundary 2, if boundary 1 $\leq$ boundary 2, we set Threshold$_{a_{BW}}(n) =$ boundary 1. Otherwise, we compute $P^A_{\text{rejected}}$ and $P^U_{\text{accepted}}$ for all CR values from both groups using the empirical CDFs $\hat{F}^A$ and $\hat{F}^U$. If a CR value exists such that $P^A_{\text{rejected}}= P^U_{\text{accepted}}$, that value is selected as Threshold$_{a_{BW}}(n)$. If no such value exists---which is possible due to the discrete nature of CR---we find the intersection point of the linear interpolation of the two discrete proportion curves, which provides the required Threshold$_{a_{BW}}(n)$.\\\\
		Fig. \ref{n4_4} illustrates the relationship between $P^A_{\text{rejected}}$ and $P^U_{\text{accepted}}$ for the case $n=4$ and $a_{BW}=4$. Note that the equality $P^A_{\text{rejected}} = P^U_{\text{accepted}}$ does not hold for any CR value from either group, and linear interpolation of the two discrete proportion curves gives Threshold$_4(4) = 0.1946$.\\\\
		The computed values of Threshold$_{a_{BW}}(n)$ for $n=3,4,\ldots,10$ and $a_{BW}=2,3,\ldots,9$ are reported in Table \ref{threshold} and portrayed in Fig. \ref{fig:all_thresholds}. The algorithmic procedure for computing the threshold is detailed in Fig. \ref{fig:flowchart_threshold}.
	\begin{table}[t!]
		\caption{Approximate Threshold$_{a_{BW}}(n)$}\label{threshold}
		\centering		
		\begin{tabular}{@{}cccccccccc@{}}
				\toprule
				\multirow{2}{*}{$n\downarrow$}&\phantom{}&\multicolumn{8}{c}{$a_{BW}\rightarrow$}\\
				\cmidrule{3-10}
				&&$2$&$3$&$4$&$5$&$6$&$7$&$8$&$9$\\
				\midrule
				$3$&&$0$&$0.25$&$0.2015$&$0.2278$&$0.1999$&$0.2368$&$0.2380$&$0.2458$\\
				$4$&&$0$&$0.1275$&$0.1946$&$0.2202$&$0.2426$&$0.2549$&$0.2615$&$0.2631$\\
				$5$&&$0$&$0.1198$&$0.18$&$0.2199$&$0.2303$&$0.2398$&$0.2471$&$0.2526$\\
				$6$&&$0$&$0.1440$&$0.1610$&$0.2148$&$0.2250$&$0.2360$&$0.2419$&$0.2472$\\
				$7$&&$0$&$0.1564$&$0.1856$&$0.2081$&$0.2225$&$0.2332$&$0.2393$&$0.2443$\\
				$8$&&$0$&$0.1430$&$0.1928$&$0.2224$&$0.2199$&$0.2301$&$0.2379$&$0.2420$\\
				$9$&&$0$&$0.1532$&$0.1928$&$0.2172$&$0.2330$&$0.2285$&$0.2350$&$0.2409$\\
				$10$&&$0$&$0.1377$&$0.1950$&$0.2330$&$0.2310$&$0.2414$&$0.2340$&$0.2394$\\		
				\bottomrule		
		\end{tabular}
		\par\smallskip \small \textit{A PCS with CR less than or equal to the corresponding threshold value is considered acceptable.}
	\end{table}
	\begin{figure}[th]
		\centering
		\begin{tikzpicture}
			\begin{axis}[
				width=10.5cm,
				height=7cm,
				xmin=3, xmax=10,
				ymin=0.10, ymax=0.28,
				xlabel={$n$},
				ylabel={Threshold},
				xtick={3,4,5,6,7,8,9,10},
				ytick={0.10, 0.14, 0.18, 0.22, 0.26},
				yticklabels={0, 0.14, 0.18, 0.22, 0.26}, 
				thin, 
				clip=false, 
				legend pos=outer north east, 
				legend cell align={left},
				legend style={draw=black, fill=white, font=\scriptsize}
				]
				
				\fill[white] ([xshift=-2pt] axis cs:3,0.115) rectangle ([xshift=2pt] axis cs:3,0.125);
				\draw[thin] ([xshift=-4pt,yshift=-2pt] axis cs:3,0.117) -- ([xshift=4pt,yshift=2pt] axis cs:3,0.117);
				\draw[thin] ([xshift=-4pt,yshift=-2pt] axis cs:3,0.123) -- ([xshift=4pt,yshift=2pt] axis cs:3,0.123);
				
				\addplot[color=black, solid, mark=square*, thin, mark size=1pt] coordinates {
					(3,0.10) (4,0.10) (5,0.10) (6,0.10) (7,0.10) (8,0.10) (9,0.10) (10,0.10)
				};
				\addlegendentry{$a_{BW}=2$}
				
				\addplot[color=blue, solid, mark=*, thin, mark size=1pt] coordinates {
					(3,0.25) (4,0.1275) (5,0.1198) (6,0.1440) (7,0.1564) (8,0.1430) (9,0.1532) (10,0.1377)
				};
				\addlegendentry{$a_{BW}=3$}
				
				\addplot[color=red, solid, mark=triangle*, thin, mark size=1pt] coordinates {
					(3,0.2015) (4,0.1946) (5,0.18) (6,0.1610) (7,0.1856) (8,0.1928) (9,0.1928) (10,0.1950)
				};
				\addlegendentry{$a_{BW}=4$}
				
				\addplot[color=green!60!black, solid, mark=diamond*, thin, mark size=1pt] coordinates {
					(3,0.2278) (4,0.2202) (5,0.2199) (6,0.2148) (7,0.2081) (8,0.2224) (9,0.2172) (10,0.2330)
				};
				\addlegendentry{$a_{BW}=5$}
				
				\addplot[color=orange, solid, mark=x, thin, mark size=1.5pt] coordinates {
					(3,0.1999) (4,0.2426) (5,0.2303) (6,0.2250) (7,0.2225) (8,0.2199) (9,0.2330) (10,0.2310)
				};
				\addlegendentry{$a_{BW}=6$}
				
				\addplot[color=purple, solid, mark=+, thin, mark size=1.5pt] coordinates {
					(3,0.2368) (4,0.2549) (5,0.2398) (6,0.2360) (7,0.2332) (8,0.2301) (9,0.2285) (10,0.2414)
				};
				\addlegendentry{$a_{BW}=7$}
				
				\addplot[color=cyan!80!black, solid, mark=o, thin, mark size=1pt] coordinates {
					(3,0.2380) (4,0.2615) (5,0.2471) (6,0.2419) (7,0.2393) (8,0.2379) (9,0.2350) (10,0.2340)
				};
				\addlegendentry{$a_{BW}=8$}
				
				\addplot[color=magenta, solid, mark=pentagon*, thin, mark size=1pt] coordinates {
					(3,0.2458) (4,0.2631) (5,0.2526) (6,0.2472) (7,0.2443) (8,0.2420) (9,0.2409) (10,0.2394)
				};
				\addlegendentry{$a_{BW}=9$}
				
			\end{axis}
		\end{tikzpicture}
		\caption{Threshold$_{a_{BW}}(n)$}
		\label{fig:all_thresholds}
	\end{figure}
	The flowchart outlining the entire framework is presented in Fig. \ref{flowchart}, and a detailed step-by-step pseudocode summarizing the analytical procedure is provided in Algorithm \ref{alg:taxicab_bwm}.
	\begin{remark}
		It is important to mention that the proposed algorithm, while demonstrated using Saaty's scale, is general and applicable to any ratio scale.
	\end{remark}
	\begin{algorithm}[t!]
		\caption{Analytical approach for weight calculation and consistency analysis}
		\label{alg:taxicab_bwm}
		\begin{algorithmic}[1]
				\renewcommand{\algorithmicrequire}{\textbf{Input:}}
				\renewcommand{\algorithmicensure}{\textbf{Output:}}
				
				\Require Set of decision criteria $C$, best criterion $c_B$, worst criterion $c_W$, best-to-other vector $A_B$, other-to-worst vector $A_W$
				\Ensure Optimal weight set(s) $W^*$, optimal TD $\epsilon^*$, CR
				\Statex \textbf{Step 0:} Define required variables and notations
				\State Define $x \in [1, \infty)$
				\State $D \gets C \setminus \{c_B, c_W\}$
				\State $u \gets a_{BW}$
				\For{each $i \in D$}
				\If{$a_{Bi} \times a_{iW} > u$}
				\State $u \gets a_{Bi} \times a_{iW}$
				\EndIf
				\EndFor
				\Statex \textbf{Step 1:} Define $f(x)$
				\State $f(x) \gets |a_{BW} - x|$
				\For{each $i \in D$}
				\If{$1\leq x \leq a_{Bi}^2$ \textbf{and} $a_{iW}\leq a_{Bi}$ }
				\State $f_i(x) \gets |a_{iW} - x/a_{Bi}|$
				\ElsIf{$1\leq x \leq a_{iW}^2$ \textbf{and} $a_{Bi}\leq a_{iW}$ }
				\State $f_i(x) \gets |a_{Bi} - x/a_{iW}|$
				\Else
				\State $f_i(x) \gets 2\sqrt{x}-a_{Bi}-a_{iW}$
				\EndIf
				\State $f(x) \gets f(x) + f_i(x)$
				\EndFor
				\Statex \textbf{Step 2:} Construct $X$
				\State $X \gets \{a_{BW}\}$
				\For{each $i \in D$}
				\State $X \gets X \cup \{a_{Bi} \times a_{iW}\}$
				\If{$\max\{a_{Bi}^2,\; a_{iW}^2\} \leq u$}
				\State $X \gets X \cup \{\max\{a_{Bi}^2,\; a_{iW}^2\}\}$
				\EndIf
				\EndFor
				\State Sort $X$ in ascending order: $X = \{x_0, x_1, \ldots, x_m\}$
				\Statex \textbf{Step 3:} Compute $\epsilon^*$ and identify all $x\in X$ where $f(x)=\epsilon^*$
				\State $\epsilon^* \gets \infty$
				\State $X^* \gets \emptyset$
				\For{each $x \in X$}
				\State $y \gets f(x)$
				\If{$y < \epsilon^*$}
				\State $\epsilon^* \gets y$
				\State $X^* \gets \{x\}$
				\ElsIf{$y = \epsilon^*$}
				\State $X^* \gets X^* \cup \{x\}$
				\EndIf
				\EndFor
				\Statex \textbf{Step 4:} Calculate CR
				\State Obtain CI from Table \ref{3ci}
				\State $\text{CR} \gets \frac{\epsilon^*}{\text{CI}}$	
				\State Obtain Threshold$_{a_{BW}}(n)$ from Table \ref{threshold}
				\If{$\text{CR} \leq \text{Threshold}_{a_{BW}}(n)$}
				\State Continue to Step 5
				\Else
				\State \textbf{return:} PCS is inadmissible. Revise $A_B$ and $A_W$ and restart the process
				\EndIf
				\algstore{mybreak}
		\end{algorithmic}
	\end{algorithm}
	
	\addtocounter{algorithm}{-1}
	\begin{algorithm*}[t!]
		\caption{(Continued)}
		\begin{algorithmic}[1]
				\algrestore{mybreak}
				\Statex \textbf{Step 5:} Determine all possible values of $\tilde{a}_{BW}^*$
				\For{$j \gets 1$ \textbf{to} $m$}
				\If{$x_{j-1} \in X^*$ \textbf{and} $x_j \in X^*$}
				\If{$f(\frac{x_{j-1} + x_j}{2}) = \epsilon^*$}
				\State $X^* \gets X^* \cup [x_{j-1}, x_j]$
				\EndIf
				\EndIf
				\EndFor
				\Statex \textbf{Step 6:} Obtain all optimally modified PCSs and corresponding optimal weight sets
				\If{$X^*$ is finite}
				\For{each $x^* \in X^*$}
				\State Substitute $\tilde{a}_{BW}^* = x^*$ into Eq. \eqref{optimal_pcs} to obtain an optimally modified PCS
				\State Compute the corresponding optimal weight set $W^*$ using Eq. \eqref{pcs_to_weights}
				\EndFor
				\ElsIf{$X^*$ contains an interval $[x_{j-1}, x_j]$}
				\State Introduce a parameter $a \in [x_{j-1}, x_j]$
				\State Substitute $a$ into Eq. \eqref{optimal_pcs} to obtain a parameterized family of optimally modified PCSs
				\State Compute the corresponding family of optimal weight sets $W^*$ using Eq. \eqref{pcs_to_weights}
				\EndIf
				\Statex \textbf{Results}
				\State \Return All optimal weight sets $W^*$, $\epsilon^*$, and CR
		\end{algorithmic}
	\end{algorithm*}
	\begin{remark}
		By Theorem \ref{minima}, if $f$ attains its global minimum at $x_{j-1}, x_j \in X$, then it attains the minimum at any point in $(x_{j-1}, x_j)$ only if $f$ is constant on $[x_{j-1}, x_j]$. Therefore, to check whether $f$ is constant on the interval, it suffices to verify whether $f$ achieves the global minimum at the midpoint $\frac{x_{j-1} + x_j}{2}$ (line 32-33 in Algorithm \ref{alg:taxicab_bwm}).
	\end{remark}
	\subsection{Numerical examples}
	In this subsection, we demonstrate the proposed framework using numerical examples covering all possible solution scenarios. While the admissibility of each PCS is verified using the established threshold values, we compute the weights in all cases---regardless of the admissibility status---to fully illustrate the model's behavior.\\\\
	\textbf{Example 1:} Let $C=\{c_1,c_2,\ldots,c_5\}$ be the set of decision criteria with $c_1$ as the best and $c_5$ as the worst criterion, and let $A_{B}=(1,2,3,5,8)$ and $A_{W}=(8,3,4,3,1)^T$ be the best-to-other and the other-to-worst vectors respectively.\\\\
	Step 1: By \eqref{function}, we have 
	\begin{align*}
		f_1(x)&=|8-x|,\\
		f_2(x)&=
		\begin{cases}
			\left|2-\frac{x}{3}\right|\quad\ \  \text{if } 1\leq x\leq 9,\\
			2\sqrt{x}-5\quad \text{otherwise},
		\end{cases}\\
		f_3(x)&=
		\begin{cases}
			\left|3-\frac{x}{4}\right|\quad\ \  \text{if } 1\leq x\leq 16,\\
			2\sqrt{x}-7\quad \text{otherwise},
		\end{cases}\\
		f_4(x)&=
		\begin{cases}
			\left|3-\frac{x}{5}\right|\quad\ \  \text{if } 1\leq x\leq 25,\\
			2\sqrt{x}-8\quad \text{otherwise},
		\end{cases}\\
		f(x)&=f_1(x)+f_2(x)+f_3(x)+f_4(x)\text{ for }x\in [1,\infty).
	\end{align*}\\
	Step 2: From \eqref{set}, we get $X=\{x_0,x_1,x_2,x_3,x_4\}=\{6,8,9,12,15\}$.\\\\
	Step 3: Theorem \ref{minima} implies that 
	\begin{align*}
		\displaystyle\min_{x\in [1,\infty)}f(x) &=\min\{f(6),f(8),f(9),f(12),f(15)\}\\
		&=\min\{5.3,3.0667,3.95,6.5282,10.4960\}\\
		&=3.0667\\
		&=f(8).
	\end{align*}\\
	So, the global minimum value of $f$ is $3.0667$, attained at $x_1=8$. Thus, $\epsilon^*=3.0667$. Fig. \ref{example_1_figure} shows the graph of $f$ in the interval $[1,26]$, which supports this conclusion and validates Theorem \ref{minima}.\\\\
	Step 4: Using \eqref{cr}, we get CR $=\frac{3.0667}{21}=0.1460$.\\
	Since CR $\leq$ Threshold$_{8}(5)=0.2471$ (from Table \ref{threshold}), the given PCS is admissible for real-world scenarios.\\\\
	Step 5: There are no consecutive $x_j$ at which $f$ attains its global minimum value. Therefore, the only possible value of $\tilde{a}_{BW}^*$ is $8$.\\\\
	Step 6: From \eqref{optimal_pcs}, the optimally modified PCS is given by $\tilde{A}_{B}^*=(1,2.6667,2,5,8)$, $\tilde{A}_{W}^*=(8,3,4,1.6,1)^T$.\\\\
	Step 7: By \eqref{pcs_to_weights}, the optimal weight set is $$W^*=\{0.4545,0.1705,0.2273,0.0909,0.0568\}.$$
	In this example, we get a unique optimal weight set. It is important to note that $w_1^* > w_5^*$, confirming that the intended best-worst ordering is preserved. Furthermore, the constraint deviations for the optimally modified PCS are as follows:
	\begin{align*}
		&|\tilde{a}_{12}^* - a_{12}| + |\tilde{a}_{25}^* - a_{25}| = 0.6667, \quad
		|\tilde{a}_{13}^* - a_{13}| + |\tilde{a}_{35}^* - a_{35}| = 1, \\
		&|\tilde{a}_{14}^* - a_{14}| + |\tilde{a}_{45}^* - a_{45}| = 1.4, \quad\quad\ 
		|\tilde{a}_{15}^* - a_{15}| = 0.
	\end{align*}
	Thus, the largest constraint deviation occurs for the pair $(a_{14}, a_{45})$.\\\\
	\begin{figure}[th]
		\centering
		\begin{tikzpicture}[
			declare function={
				f1(\x) = abs(8-\x);
				f2(\x) = (\x>=1 && \x<=9) ? abs(2-\x/3) : 2*sqrt(\x)-5;
				f3(\x) = (\x>=1 && \x<=16) ? abs(3-\x/4) : 2*sqrt(\x)-7;
				f4(\x) = (\x>=1 && \x<=25) ? abs(3-\x/5) : 2*sqrt(\x)-8;
				f(\x) = f1(\x) + f2(\x) + f3(\x) + f4(\x);
			}
			]
			\begin{axis}[
				width=12cm,
				height=7cm,
				xmin=1, xmax=26,
				ymin=0, ymax=35, 
				xlabel={$x$}, 
				ylabel={$f(x)$}, 
				legend style={at={(0.05,0.95)}, anchor=north west, cells={anchor=west}},
				ymajorgrids=true,
				xmajorgrids=true,
				grid style=dashed,
				thick,
				xtick={1, 6, 11, 16, 21, 26},
				extra x ticks={8},
				extra x tick labels={{\color{violet}$x_1=8$}},
				extra x tick style={tick label style={font=\bfseries}},
				extra y ticks={3.0667},
				extra y tick labels={{\color{violet}$3.0667$}},
				extra y tick style={tick label style={font=\bfseries}}
				]
				
				\addplot[
				color=blue,
				mark=none,
				domain=1:26,
				samples=300, 
				very thick 
				] {f(x)}; 
				\addlegendentry{$f(x)$}
				
				\addplot[
				color=red,
				dashed, 
				mark=none,
				domain=1:26,
				thick
				] {3.0667};
				\addlegendentry{Global minimum value}
				
				\node[circle, fill=black, inner sep=1.5pt] at (axis cs:8, 3.0667) {};
				\node[below, yshift=-2pt, font=\scriptsize, fill=white, inner sep=1pt] at (axis cs:8, 3.0667) {{\color{violet}$(8, 3.0667)$}};

			\end{axis}
		\end{tikzpicture}
		\caption{Graph of $f$ in $[1, 26]$ for Example 1}
		\label{example_1_figure}
	\end{figure}
	\textbf{Example 2:} Let $C=\{c_1,c_2,\ldots,c_5\}$ be the set of decision criteria with $c_1$ as the best and $c_5$ as the worst criterion, and let $A_{B}=(1,2,4,5,8)$ and $A_{W}=(8,3,4,2,1)^T$ be the best-to-other and the other-to-worst vectors respectively.\\\\
	Step 1: By \eqref{function}, we have 
	\begin{align*}
		f_1(x)&=|8-x|,\\
		f_2(x)&=
		\begin{cases}
			\left|2-\frac{x}{3}\right|\quad\ \  \text{if } 1\leq x\leq 9,\\
			2\sqrt{x}-5\quad \text{otherwise},
		\end{cases}\\
		f_3(x)&=
		\begin{cases}
			\left|4-\frac{x}{4}\right|\quad\ \  \text{if } 1\leq x\leq 16,\\
			2\sqrt{x}-8\quad \text{otherwise},
		\end{cases}\\
		f_4(x)&=
		\begin{cases}
			\left|2-\frac{x}{5}\right|\quad\ \  \text{if } 1\leq x\leq 25,\\
			2\sqrt{x}-7\quad \text{otherwise},
		\end{cases}\\
		f(x)&=f_1(x)+f_2(x)+f_3(x)+f_4(x)\text{ for }x\in [1,\infty).
	\end{align*}\\
	Step 2: From \eqref{set}, we get $X=\{x_0,x_1,x_2,x_3,x_4\}=\{6,8,9,10,16\}$.\\\\
	Step 3: Theorem \ref{minima} implies that 
	\begin{align*}
		\displaystyle\min_{x\in [1,\infty)}f(x) &=\min\{f(6),f(8),f(9),f(10),f(16)\}\\
		&=\min\{5.3,3.0667,3.95,4.8246,12.2\}\\
		&=3.0667\\
		&=f(8).
	\end{align*}\\
	So, the global minimum value of $f$ is $3.0667$, attained at $x_1=8$. Thus, $\epsilon^*=3.0667$. Fig. \ref{example_2_figure} shows the graph of $f$ in the interval $[1,26]$, which supports this conclusion and validates Theorem \ref{minima}.\\\\
	Step 4: Using \eqref{cr}, we get CR $=\frac{3.0667}{21}=0.1460$.\\
	Since CR $\leq$ Threshold$_{8}(5)=0.2471$ (from Table \ref{threshold}), the given PCS is admissible for real-world scenarios.\\\\
	Step 5: There are no consecutive $x_j$ at which $f$ attains its global minimum value. Therefore, the only possible value of $\tilde{a}_{BW}^*$ is $8$.\\\\
	Step 6: From \eqref{optimal_pcs}, we get two optimally modified PCSs as follows:
	\begin{enumerate}
		\item $(\tilde{A}_{B}^*)_1=(1,2.6667,4,5,8)$, $(\tilde{A}_{W}^*)_1=(8,3,2,1.6,1)^T$
		\item $(\tilde{A}_{B}^*)_2=(1,2.6667,2,5,8)$, $(\tilde{A}_{W}^*)_2=(8,3,4,1.6,1)^T$.
	\end{enumerate}
	\vspace{0.25cm}
	Step 7: Using \eqref{pcs_to_weights}, we get the corresponding optimal weight sets as follows:
	\begin{enumerate}
		\item $W_1^*=\{0.5128,0.1923,0.1282,0.1026,0.0641\}$
		\item $W_2^*=\{0.4545,0.1705,0.2273,0.0909,0.0568\}$.
	\end{enumerate}
	A visual comparison of these optimal weight sets is provided in Fig. \ref{fig:bar_example2}.\\\\
	In this example, we get two optimal weight sets. It is important to note that for both weight sets, the condition $w_1^* > w_5^*$ holds, which confirms that the intended best-worst ordering is preserved. Furthermore, the constraint deviations for both optimally modified PCSs are as follows:
	\begin{align*}
		&|\tilde{a}_{12}^* - a_{12}| + |\tilde{a}_{25}^* - a_{25}| = 0.6667, \quad
		|\tilde{a}_{13}^* - a_{13}| + |\tilde{a}_{35}^* - a_{35}| = 2, \\
		&|\tilde{a}_{14}^* - a_{14}| + |\tilde{a}_{45}^* - a_{45}| = 0.4, \quad\quad\ 
		|\tilde{a}_{15}^* - a_{15}| = 0.
	\end{align*}
	Thus, the largest constraint deviation occurs for the pair $(a_{13}, a_{35})$.\\\\
	\begin{figure}[th]
		\centering
		\begin{tikzpicture}[
			declare function={
				f1(\x) = abs(8-\x);
				f2(\x) = (\x>=1 && \x<=9) ? abs(2-\x/3) : 2*sqrt(\x)-5;
				f3(\x) = (\x>=1 && \x<=16) ? abs(4-\x/4) : 2*sqrt(\x)-8;
				f4(\x) = (\x>=1 && \x<=25) ? abs(2-\x/5) : 2*sqrt(\x)-7;
				f(\x) = f1(\x) + f2(\x) + f3(\x) + f4(\x);
			}
			]
			\begin{axis}[
				width=12cm,
				height=7cm,
				xmin=1, xmax=26,
				ymin=0, ymax=35, 
				xlabel={$x$}, 
				ylabel={$f(x)$}, 
				legend style={at={(0.05,0.95)}, anchor=north west, cells={anchor=west}},
				ymajorgrids=true,
				xmajorgrids=true,
				grid style=dashed,
				thick,
				xtick={1, 6, 11, 16, 21, 26},
				extra x ticks={8},
				extra x tick labels={{\color{violet}$x_1=8$}},
				extra x tick style={tick label style={font=\bfseries}},
				extra y ticks={3.0667},
				extra y tick labels={{\color{violet}$3.0667$}},
				extra y tick style={tick label style={font=\bfseries}}
				]
				
				\addplot[
				color=blue,
				mark=none,
				domain=1:26,
				samples=300, 
				very thick 
				] {f(x)}; 
				\addlegendentry{$f(x)$}
				
				\addplot[
				color=red,
				dashed, 
				mark=none,
				domain=1:26,
				thick
				] {3.0667};
				\addlegendentry{Global minimum value}
				
				\node[circle, fill=black, inner sep=1.5pt] at (axis cs:8, 3.0667) {};
				\node[below, yshift=-2pt, font=\scriptsize, fill=white, inner sep=1pt] at (axis cs:8, 3.0667) {{\color{violet}$(8, 3.0667)$}};

			\end{axis}
		\end{tikzpicture}
		\caption{Graph of $f$ in $[1, 26]$ for Example 2}
		\label{example_2_figure}
	\end{figure}
	\begin{figure}[th]
		\centering
		\begin{tikzpicture}
			\begin{axis}[
				ybar=3pt,
				bar width=12pt,
				area legend, 
				width=12cm,
				height=7cm,
				ymin=0, 
				ymax=0.65, 
				ylabel={Optimal Weights},
				xlabel={Criteria},
				symbolic x coords={$c_1$, $c_2$, $c_3$, $c_4$, $c_5$},
				xtick=data,
				nodes near coords,
				every node near coord/.append style={
					font=\scriptsize, 
					rotate=90, 
					anchor=west, 
					/pgf/number format/fixed, 
					/pgf/number format/precision=4
				},
				legend style={at={(0.95,0.95)}, anchor=north east}, 
				enlarge x limits=0.15,
				ymajorgrids=true,
				grid style=dashed
				]
				
				\addplot[fill=blue!60, draw=black] coordinates {
					($c_1$, 0.5128) ($c_2$, 0.1923) ($c_3$, 0.1282) ($c_4$, 0.1026) ($c_5$, 0.0641)
				};
				
				\addplot[fill=red!60, draw=black] coordinates {
					($c_1$, 0.4545) ($c_2$, 0.1705) ($c_3$, 0.2273) ($c_4$, 0.0909) ($c_5$, 0.0568)
				};
				
				\legend{$W_1^*$, $W_2^*$}
			\end{axis}
		\end{tikzpicture}
		\caption{Comparison of optimal weight sets for Example 2}
		\label{fig:bar_example2}
	\end{figure}
	\textbf{Example 3:} Let $C=\{c_1,c_2,\ldots,c_5\}$ be the set of decision criteria with $c_1$ as the best and $c_5$ as the worst criterion, and let $A_{B}=(1,1,1,2,4)$ and $A_{W}=(4,1,1,3,1)^T$ be the best-to-other and the other-to-worst vectors respectively.\\\\
	Step 1: By \eqref{function}, we have 
	\begin{align*}
		f_1(x)&=|4-x|,\\
		f_2(x)&=2\sqrt{x}-2,\\
		f_3(x)&=2\sqrt{x}-2,\\
		f_4(x)&=
		\begin{cases}
			\left|2-\frac{x}{3}\right|\quad\ \  \text{if } 1\leq x\leq 9,\\
			2\sqrt{x}-5\quad \text{otherwise},
		\end{cases}\\
		f(x)&=f_1(x)+f_2(x)+f_3(x)+f_4(x)\text{ for }x\in [1,\infty).
	\end{align*}\\
	Step 2: From \eqref{set}, we get $X=\{x_0,x_1,x_2\}=\{1,4,6\}$.\\\\
	Step 3: Theorem \ref{minima} implies that 
	\begin{align*}
		\displaystyle\min_{x\in [1,\infty)}f(x) &=\min\{f(1),f(4),f(6)\}\\
		&=\min\{4.6667,4.6667,7.7980\}\\
		&=4.6667\\
		&=f(1)\\
		&=f(4).
	\end{align*}\\
	So, the global minimum value of $f$ is $4.6667$, attained at $x_0=1$ and $x_1=4$. Thus, $\epsilon^*=4.6667$. Fig. \ref{example_3_figure} shows the graph of $f$ in the interval $[1,16]$, which supports this conclusion and validates Theorem \ref{minima}.\\\\
	Step 4: Using \eqref{cr}, we get CR $=\frac{4.6667}{9}=0.5185$.\\
	Since CR $>$ Threshold$_{4}(5)=0.18$ (from Table \ref{threshold}), the given PCS is not admissible for real-world scenarios.\\\\
	Step 5: $f$ attains its global minimum value at $x_0=1$ and $x_1=4$. From $f(x)=4\sqrt{x}-\frac{4}{3}x+2$ for $1\leq x\leq 4$, it follows that $f$ is nonconstant on $[1,4]$. Thus, $(\tilde{a}_{BW}^*)_1=1$ and $(\tilde{a}_{BW}^*)_2=4$ are two possible values of $\tilde{a}_{BW}^*$.\\\\
	Step 6: From \eqref{optimal_pcs}, we get two optimally modified PCSs, one for each value of $\tilde{a}_{BW}^*$, as follows:
	\begin{enumerate}
		\item $(\tilde{A}_{B}^*)_1=(1,1,1,0.3333,1)$, $(\tilde{A}_{W}^*)_1=(1,1,1,3,1)^T$
		\item $(\tilde{A}_{B}^*)_2=(1,2,2,1.3333,4)$, $(\tilde{A}_{W}^*)_2=(4,2,2,3,1)^T$.
	\end{enumerate}
	\vspace{0.25cm}
	Step 7: Using \eqref{pcs_to_weights}, we get the corresponding optimal weight sets as follows:
	\begin{enumerate}
		\item $W_1^*=\{0.1429,0.1429,0.1429,0.4286,0.1429\}$
		\item $W_2^*=\{0.3333,0.1667,0.1667,0.25,0.0833\}$.
	\end{enumerate}
	A visual comparison of these optimal weight sets is provided in Fig. \ref{fig:bar_example3}.\\\\
	In this example, we get two optimal weight sets. It is important to note that for $((\tilde{A}_{B}^*)_1,(\tilde{A}_{W}^*)_1)$, we have $\tilde{a}_{45}^*>\tilde{a}_{15}^*=\tilde{a}_{BW}$. This results in a lower weight for the best criterion $c_1$ compared to $c_4$ in $W_1^*$, which violates the intended ordering and makes $W_1^*$ less preferable than $W_2^*$. Moreover, for $W_2^*$, the condition $w_1^* > w_5^*$ holds, which confirms that the intended best-worst ordering is preserved. Furthermore, the constraint deviations for $((\tilde{A}_B^*)_2,(\tilde{A}_W^*)_2)$ are as follows:
	\begin{align*}
		&|(\tilde{a}_{12}^*)_2 - a_{12}| + |(\tilde{a}_{25}^*)_2 - a_{25}| = 2, \quad\quad\quad\ 
		|(\tilde{a}_{13}^*)_2 - a_{13}| + |(\tilde{a}_{35}^*)_2 - a_{35}| = 2, \\
		&|(\tilde{a}_{14}^*)_2 - a_{14}| + |(\tilde{a}_{45}^*)_2 - a_{45}| = 0.6667, \quad 
		|(\tilde{a}_{15}^*)_2 - a_{15}| = 0.
	\end{align*}
	Thus, the largest constraint deviation occurs for the pairs $(a_{12}, a_{25})$ and $(a_{13}, a_{35})$.\\\\
	\begin{figure}[th]
		\centering
		\begin{tikzpicture}[
			declare function={
				f1(\x) = abs(4-\x);
				f2(\x) = 2*sqrt(\x)-2;
				f3(\x) = 2*sqrt(\x)-2;
				f4(\x) = (\x>=1 && \x<=9) ? abs(2-\x/3) : 2*sqrt(\x)-5;
				f(\x) = f1(\x) + f2(\x) + f3(\x) + f4(\x);
			}
			]
			\begin{axis}[
				width=12cm,
				height=7cm,
				xmin=1, xmax=16,
				ymin=0, ymax=30, 
				xlabel={$x$}, 
				ylabel={$f(x)$}, 
				legend style={at={(0.05,0.95)}, anchor=north west, cells={anchor=west}},
				ymajorgrids=true,
				xmajorgrids=true,
				grid style=dashed,
				thick,
				xtick={7, 10, 13, 16},
				extra x ticks={1, 4},
				extra x tick labels={{\color{violet}$x_0=1$}, {\color{violet}$x_1=4$}},
				extra x tick style={tick label style={font=\bfseries}},
				extra y ticks={4.6667},
				extra y tick labels={{\color{violet}$4.6667$}},
				extra y tick style={tick label style={font=\bfseries}}
				]
				
				\addplot[
				color=blue,
				mark=none,
				domain=1:16,
				samples=300, 
				very thick 
				] {f(x)}; 
				\addlegendentry{$f(x)$}
				
				\addplot[
				color=red,
				dashed, 
				mark=none,
				domain=1:16,
				thick
				] {4.6667};
				\addlegendentry{Global minimum value}
				
				\node[circle, fill=black, inner sep=1.5pt] at (axis cs:1, 4.6667) {};
				\node[below right, font=\scriptsize, fill=white, inner sep=1pt, yshift=-2pt] at (axis cs:1, 4.6667) {{\color{violet}$(1, 4.6667)$}};
				
				\node[circle, fill=black, inner sep=1.5pt] at (axis cs:4, 4.6667) {};
				\node[below, font=\scriptsize, fill=white, inner sep=1pt, yshift=-2pt] at (axis cs:4, 4.6667) {{\color{violet}$(4, 4.6667)$}};
				
			\end{axis}
		\end{tikzpicture}
		\caption{Graph of $f$ in $[1, 16]$ for Example 3}
		\label{example_3_figure}
	\end{figure}
	\begin{figure}[th]
		\centering
		\begin{tikzpicture}
			\begin{axis}[
				ybar=3pt,
				bar width=12pt,
				area legend, 
				width=12cm,
				height=7cm,
				ymin=0, 
				ymax=0.55, 
				ylabel={Optimal Weights},
				xlabel={Criteria},
				symbolic x coords={$c_1$, $c_2$, $c_3$, $c_4$, $c_5$},
				xtick=data,
				nodes near coords,
				every node near coord/.append style={
					font=\scriptsize, 
					rotate=90, 
					anchor=west, 
					/pgf/number format/fixed, 
					/pgf/number format/precision=4
				},
				legend style={at={(0.95,0.95)}, anchor=north east}, 
				enlarge x limits=0.15,
				ymajorgrids=true,
				grid style=dashed
				]
				
				\addplot[fill=blue!60, draw=black] coordinates {
					($c_1$, 0.1429) ($c_2$, 0.1429) ($c_3$, 0.1429) ($c_4$, 0.4286) ($c_5$, 0.1429)
				};
				
				\addplot[fill=red!60, draw=black] coordinates {
					($c_1$, 0.3333) ($c_2$, 0.1667) ($c_3$, 0.1667) ($c_4$, 0.2500) ($c_5$, 0.0833)
				};
				
				\legend{$W_1^*$, $W_2^*$}
			\end{axis}
		\end{tikzpicture}
		\caption{Comparison of optimal weight sets for Example 3}
		\label{fig:bar_example3}
	\end{figure}
	\textbf{Example 4:} Let $C=\{c_1,c_2,\ldots,c_5\}$ be the set of decision criteria with $c_1$ as the best and $c_5$ as the worst criterion, and let $A_{B}=(1,1,1,1,9)$ and $A_{W}=(9,1,1,5,1)^T$ be the best-to-other and the other-to-worst vectors respectively.\\\\
	Step 1: By \eqref{function}, we have 
	\begin{align*}
		f_1(x)&=|9-x|,\\
		f_2(x)&=2\sqrt{x}-2,\\
		f_3(x)&=2\sqrt{x}-2,\\
		f_4(x)&=
		\begin{cases}
			\left|1-\frac{x}{5}\right|\quad\ \  \text{if } 1\leq x\leq 25,\\
			2\sqrt{x}-6\quad \text{otherwise},
		\end{cases}\\
		f(x)&=f_1(x)+f_2(x)+f_3(x)+f_4(x)\text{ for }x\in [1,\infty).
	\end{align*}\\
	Step 2: From \eqref{set}, we get $X=\{x_0,x_1,x_2\}=\{1,5,9\}$.\\\\
	Step 3: Theorem \ref{minima} implies that 
	\begin{align*}
		\displaystyle\min_{x\in [1,\infty)}f(x) &=\min\{f(1),f(5),f(9)\}\\
		&=\min\{8.8,8.9443,8.8\}\\
		&=8.8\\
		&=f(1)\\
		&=f(9).
	\end{align*}\\
	So, the global minimum value of $f$ is $8.8$, attained at $x_0=1$ and $x_2=9$. Thus, $\epsilon^*=8.8$. Fig. \ref{example_4_figure} shows the graph of $f$ in the interval $[1,25]$, which supports this conclusion and validates Theorem \ref{minima}.\\\\
	Step 4: Using \eqref{cr}, we get CR $=\frac{8.8}{24}=0.3667$.\\
	Since CR $>$ Threshold$_{9}(5)=0.2526$ (from Table \ref{threshold}), the given PCS is not admissible for real-world scenarios.\\\\
	Step 5: There are no consecutive $x_j$ at which $f$ attains its global minimum value. Therefore, $(\tilde{a}_{BW}^*)_1=1$ and $(\tilde{a}_{BW}^*)_2=9$ are two possible values of $\tilde{a}_{BW}^*$.\\\\
	Step 6: From \eqref{optimal_pcs}, we get two optimally modified PCSs, one for each value of $\tilde{a}_{BW}^*$, as follows:
	\begin{enumerate}
		\item $(\tilde{A}_{B}^*)_1=(1,1,1,0.2,1)$, $(\tilde{A}_{W}^*)_1=(1,1,1,5,1)^T$
		\item $(\tilde{A}_{B}^*)_2=(1,3,3,1.8,9)$, $(\tilde{A}_{W}^*)_2=(9,3,3,5,1)^T$.
	\end{enumerate}
	\vspace{0.25cm}
	Step 7: Using \eqref{pcs_to_weights}, we get the corresponding optimal weight sets as follows:
	\begin{enumerate}
		\item $W_1^*=\{0.1111,0.1111,0.1111,0.5556,0.1111\}$
		\item $W_2^*=\{0.4286,0.1429,0.1429,0.2381,0.0476\}$.
	\end{enumerate}
	A visual comparison of these optimal weight sets is provided in Fig. \ref{fig:bar_example4}.\\\\
	In this example, we get two optimal weight sets. It is important to note that for $((\tilde{A}_{B}^*)_1,(\tilde{A}_{W}^*)_1)$, we have $\tilde{a}_{45}^*>\tilde{a}_{15}^*=\tilde{a}_{BW}$. This results in a lower weight for the best criterion $c_1$ compared to $c_4$ in $W_1^*$, which violates the intended ordering and makes $W_1^*$ less preferable than $W_2^*$. Moreover, for $W_2^*$, the condition $w_1^* > w_5^*$ holds, which confirms that the intended best-worst ordering is preserved. Furthermore, the constraint deviations for $((\tilde{A}_B^*)_2,(\tilde{A}_W^*)_2)$ are as follows:
	\begin{align*}
		&|(\tilde{a}_{12}^*)_2 - a_{12}| + |(\tilde{a}_{25}^*)_2 - a_{25}| = 4, \quad\quad\quad\ 
		|(\tilde{a}_{13}^*)_2 - a_{13}| + |(\tilde{a}_{35}^*)_2 - a_{35}| = 4, \\
		&|(\tilde{a}_{14}^*)_2 - a_{14}| + |(\tilde{a}_{45}^*)_2 - a_{45}| = 0.8, \quad\quad\  
		|(\tilde{a}_{15}^*)_2 - a_{15}| = 0.
	\end{align*}
	Thus, the largest constraint deviation occurs for the pairs $(a_{12}, a_{25})$ and $(a_{13}, a_{35})$.\\\\
	\begin{figure}[th]
		\centering
		\begin{tikzpicture}[
			declare function={
				f1(\x) = abs(9-\x);
				f2(\x) = 2*sqrt(\x)-2;
				f3(\x) = 2*sqrt(\x)-2;
				f4(\x) = (\x>=1 && \x<=25) ? abs(1-\x/5) : 2*sqrt(\x)-6;
				f(\x) = f1(\x) + f2(\x) + f3(\x) + f4(\x);
			}
			]
			\begin{axis}[
				width=12cm,
				height=7cm,
				xmin=1, xmax=25,
				ymin=0, ymax=42,
				ytick={0, 6, 12, 18, 24, 30, 36, 42}, 
				xlabel={$x$}, 
				ylabel={$f(x)$}, 
				legend style={at={(0.05,0.95)}, anchor=north west, cells={anchor=west}},
				ymajorgrids=true,
				xmajorgrids=true,
				grid style=dashed,
				thick,
				xtick={5, 13, 17, 21, 25},
				extra x ticks={1, 9},
				extra x tick labels={{\color{violet}$x_0=1$}, {\color{violet}$x_2=9$}},
				extra x tick style={tick label style={font=\bfseries}},
				extra y ticks={8.8},
				extra y tick labels={{\color{violet}$8.8$}},
				extra y tick style={tick label style={font=\bfseries}}
				]
				
				\addplot[
				color=blue,
				mark=none,
				domain=1:25,
				samples=300, 
				very thick 
				] {f(x)}; 
				\addlegendentry{$f(x)$}
				
				\addplot[
				color=red,
				dashed, 
				mark=none,
				domain=1:25,
				thick
				] {8.8};
				\addlegendentry{Global minimum value}
				
				\node[circle, fill=black, inner sep=1.5pt] at (axis cs:1, 8.8) {};
				\node[below right, font=\scriptsize, fill=white, inner sep=1pt, yshift=-2pt] at (axis cs:1, 8.8) {{\color{violet}$(1, 8.8)$}};
				
				\node[circle, fill=black, inner sep=1.5pt] at (axis cs:9, 8.8) {};
				\node[below, font=\scriptsize, fill=white, inner sep=1pt, yshift=-2pt] at (axis cs:9, 8.8) {{\color{violet}$(9, 8.8)$}};
				
			\end{axis}
		\end{tikzpicture}
		\caption{Graph of $f$ in $[1, 25]$ for Example 4}
		\label{example_4_figure}
	\end{figure}
	\begin{figure}[th]
		\centering
		\begin{tikzpicture}
			\begin{axis}[
				ybar=3pt,
				bar width=12pt,
				area legend, 
				width=12cm,
				height=7cm,
				ymin=0, 
				ymax=0.70, 
				ylabel={Optimal Weights},
				xlabel={Criteria},
				symbolic x coords={$c_1$, $c_2$, $c_3$, $c_4$, $c_5$},
				xtick=data,
				nodes near coords,
				every node near coord/.append style={
					font=\scriptsize, 
					rotate=90, 
					anchor=west, 
					/pgf/number format/fixed, 
					/pgf/number format/precision=4
				},
				legend style={at={(0.95,0.95)}, anchor=north east}, 
				enlarge x limits=0.15,
				ymajorgrids=true,
				grid style=dashed
				]
				
				\addplot[fill=blue!60, draw=black] coordinates {
					($c_1$, 0.1111) ($c_2$, 0.1111) ($c_3$, 0.1111) ($c_4$, 0.5556) ($c_5$, 0.1111)
				};
				
				\addplot[fill=red!60, draw=black] coordinates {
					($c_1$, 0.4286) ($c_2$, 0.1429) ($c_3$, 0.1429) ($c_4$, 0.2381) ($c_5$, 0.0476)
				};
				
				\legend{$W_1^*$, $W_2^*$}
			\end{axis}
		\end{tikzpicture}
		\caption{Comparison of optimal weight sets for Example 4}
		\label{fig:bar_example4}
	\end{figure}
	\textbf{Example 5:} Let $C=\{c_1,c_2,\ldots,c_5\}$ be the set of decision criteria with $c_1$ as the best and $c_5$ as the worst criterion, and let $A_{B}=(1,2,2,2,9)$ and $A_{W}=(9,3,3,3,1)^T$ be the best-to-other and the other-to-worst vectors respectively.\\\\
	Step 1: By \eqref{function}, we have 
	\begin{align*}
		f_1(x)&=|9-x|,\\
		f_2(x)&=\begin{cases}
			\left|2-\frac{x}{3}\right|\quad\ \  \text{if } 1\leq x\leq 9,\\
			2\sqrt{x}-5\quad \text{otherwise},
		\end{cases}\\
		f_3(x)&=\begin{cases}
			\left|2-\frac{x}{3}\right|\quad\ \  \text{if } 1\leq x\leq 9,\\
			2\sqrt{x}-5\quad \text{otherwise},
		\end{cases}\\
		f_4(x)&=
		\begin{cases}
			\left|2-\frac{x}{3}\right|\quad\ \  \text{if } 1\leq x\leq 9,\\
			2\sqrt{x}-5\quad \text{otherwise},
		\end{cases}\\
		f(x)&=f_1(x)+f_2(x)+f_3(x)+f_4(x)\text{ for }x\in [1,\infty).
	\end{align*}\\
	Step 2: From \eqref{set}, we get $X=\{x_0,x_1\}=\{6,9\}$.\\\\
	Step 3: Theorem \ref{minima} implies that 
	\begin{align*}
		\displaystyle\min_{x\in [1,\infty)}f(x) &=\min\{f(6),f(9)\}\\
		&=\min\{3,3\}\\
		&=3\\
		&=f(6)\\
		&=f(9).
	\end{align*}
	So, the global minimum value of $f$ is $3$, attained at $x_0=6$ and $x_1=9$. Thus, $\epsilon^*=3$. Fig. \ref{example_5_figure} shows the graph of $f$ in the interval $[1,21]$, which supports this conclusion and validates Theorem \ref{minima}.\\\\
	Step 4: Using \eqref{cr}, we get CR $=\frac{3}{24}=0.1250$.\\
	Since CR $\leq$ Threshold$_{9}(5)=0.2526$ (from Table \ref{threshold}), the given PCS is admissible for real-world scenarios.\\\\
	Step 5: $f$ attains its global minimum value at $x_0=6$ and $x_1=9$. Since $f(x)=3$ for $6\leq x\leq 9$, all possible values of $\tilde{a}_{BW}^*$ are $[6,9]$.\\\\
	Step 6: From \eqref{optimal_pcs}, we get infinitely many optimally modified PCSs given by $\tilde{A}_{B}^*=(1,\frac{a}{3},\frac{a}{3},\frac{a}{3},a)$, $\tilde{A}_{W}^*=(a,3,3,3,1)^T$, $a\in [6,9]$.\\\\
	Step 7: Using \eqref{pcs_to_weights}, we get infinitely many optimal weight sets $$W^*=\biggl\{\frac{a}{a+10},\frac{3}{a+10},\frac{3}{a+10},\frac{3}{a+10},\frac{1}{a+10}\biggr\},\ a\in [6,9].$$
	A visual representation of variation in these optimal weights with respect to parameter `$a$' is provided in Fig. \ref{fig:line_example5}.\\\\
	In this example, we get infinitely many optimal weight sets. It is important to note that for all optimal weight sets, the condition $w_1^* > w_5^*$ holds, which confirms that the intended best-worst ordering is preserved. Furthermore, the constraint deviations for the optimally modified PCSs are as follows:
	\begin{align*}
		&|\tilde{a}_{12}^* - a_{12}| + |\tilde{a}_{25}^* - a_{25}| = \frac{a-6}{3}, \quad
		|\tilde{a}_{13}^* - a_{13}| + |\tilde{a}_{35}^* - a_{35}| = \frac{a-6}{3}, \\
		&|\tilde{a}_{14}^* - a_{14}| + |\tilde{a}_{45}^* - a_{45}| = \frac{a-6}{3}, \quad 
		|\tilde{a}_{15}^* - a_{15}| = 9-a.
	\end{align*}
	Thus, for the optimal weight sets, the location of the largest constraint deviation depends on $a$:
	\begin{itemize}
		\item For $a\in[6,8.25)$, it occurs for $a_{25}$.
		\item At $a=8.25$, it occurs for $a_{25}$ and the pairs $(a_{12}, a_{25})$, $(a_{13}, a_{35})$ and $(a_{14}, a_{45})$.
		\item For $a\in (8.25,9]$, it occurs for $a_{25}$ and the pairs $(a_{12}, a_{25})$, $(a_{13}, a_{35})$ and $(a_{14}, a_{45})$.
	\end{itemize}
	\begin{figure}[th]
		\centering
		\begin{tikzpicture}[
			declare function={
				f1(\x) = abs(9-\x);
				f2(\x) = (\x>=1 && \x<=9) ? abs(2-\x/3) : 2*sqrt(\x)-5;
				f3(\x) = (\x>=1 && \x<=9) ? abs(2-\x/3) : 2*sqrt(\x)-5;
				f4(\x) = (\x>=1 && \x<=9) ? abs(2-\x/3) : 2*sqrt(\x)-5;
				f(\x) = f1(\x) + f2(\x) + f3(\x) + f4(\x);
			}
			]
			\begin{axis}[
				width=12cm,
				height=7cm,
				xmin=1, xmax=21,
				ymin=0, ymax=35, 
				xlabel={$x$}, 
				ylabel={$f(x)$}, 
				legend style={at={(0.05,0.95)}, anchor=north west, cells={anchor=west}},
				ymajorgrids=true,
				xmajorgrids=true,
				grid style=dashed,
				thick,
				xtick={1, 11, 16, 21},
				extra x ticks={6, 9},
				extra x tick labels={{\color{violet}$x_0=6$}, {\color{violet}$x_1=9$}},
				extra x tick style={tick label style={font=\bfseries}},
				extra y ticks={3},
				extra y tick labels={{\color{violet}$3$}},
				extra y tick style={tick label style={font=\bfseries}}
				]
				
				\addplot[
				color=blue,
				mark=none,
				domain=1:21,
				samples=300, 
				very thick 
				] {f(x)}; 
				\addlegendentry{$f(x)$}
				
				\addplot[
				color=red,
				dashed, 
				mark=none,
				domain=1:21,
				thick
				] {3};
				\addlegendentry{Global minimum value}
				
				\node[circle, fill=black, inner sep=1.5pt] at (axis cs:6, 3) {};
				\node[below, font=\scriptsize, fill=white, inner sep=1pt, yshift=-2pt] at (axis cs:6, 3) {{\color{violet}$(6, 3)$}};
				
				\node[circle, fill=black, inner sep=1.5pt] at (axis cs:9, 3) {};
				\node[below, font=\scriptsize, fill=white, inner sep=1pt, yshift=-2pt] at (axis cs:9, 3) {{\color{violet}$(9, 3)$}};
				
				\node[above, font=\scriptsize, fill=white, inner sep=1pt, yshift=4pt] at (axis cs:7.5, 3.1) {{\color{violet}$(a,3), 6<a<9$}};

			\end{axis}
		\end{tikzpicture}
		\caption{Graph of $f$ in $[1, 21]$ for Example 5}
		\label{example_5_figure}
	\end{figure}
	\begin{figure}[th]
		\centering
		\begin{tikzpicture}
			\begin{axis}[
				width=12cm,
				height=7cm,
				domain=6:9,
				samples=1000,
				ymin=0, 
				ymax=0.65, 
				xlabel={Parameter $a$},
				ylabel={Optimal Weights},
				legend style={at={(0.05,0.95)}, anchor=north west}, 
				ymajorgrids=true,
				xmajorgrids=true,
				grid style=dashed,
				thick,
				clip=false 
				]
				
				\addplot[color=blue, mark=none] {x / (x + 10)};
				\addlegendentry{$w_1^*$}
				
				\addplot[color=red, mark=none] {3 / (x + 10)};
				\addlegendentry{$w_2^*, w_3^*, w_4^*$}
				
				\addplot[color=teal, mark=none] {1 / (x + 10)};
				\addlegendentry{$w_5^*$}
				
				\node[circle, fill=blue, inner sep=1.2pt] at (axis cs:6, 0.3750) {};
				\node[below, xshift=3pt, yshift=-2pt, fill=white, inner sep=1pt, font=\scriptsize, text=blue] at (axis cs:6, 0.3750) {0.3750};
				
				\node[circle, fill=red, inner sep=1.2pt]  at (axis cs:6, 0.1875) {};
				\node[below, xshift=3pt, yshift=-2pt, fill=white, inner sep=1pt, font=\scriptsize, text=red]  at (axis cs:6, 0.1875) {0.1875};
				
				\node[circle, fill=teal, inner sep=1.2pt] at (axis cs:6, 0.0625) {};
				\node[below, xshift=3pt, yshift=-2pt, fill=white, inner sep=1pt, font=\scriptsize, text=teal] at (axis cs:6, 0.0625) {0.0625};
				
				\node[circle, fill=blue, inner sep=1.2pt] at (axis cs:9, 0.4737) {};
				\node[below, xshift=-3pt, yshift=-2pt, fill=white, inner sep=1pt, font=\scriptsize, text=blue] at (axis cs:9, 0.4737) {0.4737};
				
				\node[circle, fill=red, inner sep=1.2pt]  at (axis cs:9, 0.1579) {};
				\node[below, xshift=-3pt, yshift=-2pt, fill=white, inner sep=1pt, font=\scriptsize, text=red]  at (axis cs:9, 0.1579) {0.1579};
				
				\node[circle, fill=teal, inner sep=1.2pt] at (axis cs:9, 0.0526) {};
				\node[below, xshift=-3pt, yshift=-2pt, fill=white, inner sep=1pt, font=\scriptsize, text=teal] at (axis cs:9, 0.0526) {0.0526};
				
			\end{axis}
		\end{tikzpicture}
		\caption{Variation in the optimal weights for Example 5}
		\label{fig:line_example5}
	\end{figure}
	\textbf{Example 6:} Let $C=\{c_1,c_2,\ldots,c_{10}\}$ be the set of decision criteria with $c_1$ as the best and $c_{10}$ as the worst criterion, and let $A_{B}=(1,2,3,5,3,6,4,4,2,9)$ and $A_{W}=(9,4,3,2,4,2,4,2,6,1)^T$ be the best-to-other and the other-to-worst vectors respectively.\\\\
	Step 1: By \eqref{function}, we have 
	\begin{align*}
		f_1(x)&=|9-x|,\\
		f_2(x)&=
		\begin{cases}
			\left|2-\frac{x}{4}\right|\quad\ \  \text{if } 1\leq x\leq 16,\\
			2\sqrt{x}-6\quad \text{otherwise},
		\end{cases}\\
		f_3(x)&=
		\begin{cases}
			\left|3-\frac{x}{3}\right|\quad\ \  \text{if } 1\leq x\leq 9,\\
			2\sqrt{x}-6\quad \text{otherwise},
		\end{cases}\\
		f_4(x)&=
		\begin{cases}
			\left|2-\frac{x}{5}\right|\quad\ \  \text{if } 1\leq x\leq 25,\\
			2\sqrt{x}-7\quad \text{otherwise},
		\end{cases}\\
		f_5(x)&=
		\begin{cases}
			\left|3-\frac{x}{4}\right|\quad\ \  \text{if } 1\leq x\leq 16,\\
			2\sqrt{x}-7\quad \text{otherwise},
		\end{cases}\\
		f_6(x)&=
		\begin{cases}
			\left|2-\frac{x}{6}\right|\quad\ \  \text{if } 1\leq x\leq 36,\\
			2\sqrt{x}-8\quad \text{otherwise},
		\end{cases}\\
		f_7(x)&=
		\begin{cases}
			\left|4-\frac{x}{4}\right|\quad\ \  \text{if } 1\leq x\leq 16,\\
			2\sqrt{x}-8\quad \text{otherwise},
		\end{cases}\\
		f_8(x)&=
		\begin{cases}
			\left|2-\frac{x}{4}\right|\quad\ \  \text{if } 1\leq x\leq 16,\\
			2\sqrt{x}-6\quad \text{otherwise},
		\end{cases}\\
		f_9(x)&=
		\begin{cases}
			\left|2-\frac{x}{6}\right|\quad\ \  \text{if } 1\leq x\leq 36,\\
			2\sqrt{x}-8\quad \text{otherwise},
		\end{cases}\\
		f(x)&=f_1(x)+f_2(x)+\ldots+f_9(x)\text{ for }x\in [1,\infty).
	\end{align*}\\
	Step 2: From \eqref{set}, we get $X=\{x_0,x_1,x_2,x_3,x_4\}=\{8,9,10,12,16\}$.\\\\
	Step 3: Theorem \ref{minima} implies that 
	\begin{align*}
		\displaystyle\min_{x\in [1,\infty)}f(x) &=\min\{f(8),f(9),f(10),f(12),f(16)\}\\
		&=\min\{6.0667,4.2,4.9912,7.3282,16.5333\}\\
		&=4.2\\
		&=f(9).
	\end{align*}\\
	So, the global minimum value of $f$ is $4.2$, attained at $x_1=9$. Thus, $\epsilon^*=4.2$. Fig. \ref{example_6_figure} shows the graph of $f$ in the interval $[1,21]$, which supports this conclusion and validates Theorem \ref{minima}.\\\\
	Step 4: Using \eqref{cr}, we get CR $=\frac{4.2}{64}=0.0656$.\\
	Since CR $\leq$ Threshold$_{9}(10)=0.2394$ (from Table \ref{threshold}), the given PCS is admissible for real-world scenarios.\\\\
	Step 5: There are no consecutive $x_j$ at which $f$ attains its global minimum value. Therefore, the only possible value of $\tilde{a}_{BW}^*$ is $9$.\\\\
	Step 6: From \eqref{optimal_pcs}, we get two optimally modified PCSs as follows:
	\begin{enumerate}
		\item $(\tilde{A}_{B}^*)_1=(1,2.25,3,5,2.25,6,4,4,1.5,9)$,\\ $(\tilde{A}_{W}^*)_1=(9,4,3,1.8,4,1.5,2.25,2.25,6,1)^T$
		\item $(\tilde{A}_{B}^*)_2=(1,2.25,3,5,2.25,6,2.25,4,1.5,9)$,\\ $(\tilde{A}_{W}^*)_2=(9,4,3,1.8,4,1.5,4,2.25,6,1)^T$.
	\end{enumerate}
	Step 7: Using \eqref{pcs_to_weights}, we get the corresponding optimal weight sets as follows:
	\begin{enumerate}
		\item $W_1^*=\left\{0.2586,0.1149,0.0862,0.0517,0.1149,0.0431,0.0647,0.0647,0.1724,\right.$ $\left.\quad\ \ \quad \ \ 0.0287\right\}$
		\item $W_2^*=\left\{0.2462,0.1094,0.0821,0.0492,0.1094,0.0410,0.1094,0.0616,0.1642,\right.$ $\left. \quad\ \ \quad \ \ 0.0274\right\}$.
	\end{enumerate}
	A visual comparison of these optimal weight sets is provided in Fig. \ref{fig:bar_example6}.\\\\
	In this example, we get two optimal weight sets. It is important to note that for both weight sets, the condition $w_1^* > w_{10}^*$ holds, which confirms that the intended best-worst ordering is preserved. Furthermore, the constraint deviations for both optimally modified PCSs are as follows:
	\begin{align*}
		&|\tilde{a}_{12}^* - a_{12}| + |\tilde{a}_{2,10}^* - a_{2,10}| = 0.25, \quad
		|\tilde{a}_{13}^* - a_{13}| + |\tilde{a}_{3,10}^* - a_{3,10}| = 0, \\
		&|\tilde{a}_{14}^* - a_{14}| + |\tilde{a}_{4,10}^* - a_{4,10}| = 0.2, \quad\ 
		|\tilde{a}_{15}^* - a_{15}| + |\tilde{a}_{5,10}^* - a_{5,10}| = 0.75, \\
		&|\tilde{a}_{16}^* - a_{16}| + |\tilde{a}_{6,10}^* - a_{6,10}| = 0.5, \quad\ 
		|\tilde{a}_{17}^* - a_{17}| + |\tilde{a}_{7,10}^* - a_{7,10}| = 1.75, \\
		&|\tilde{a}_{18}^* - a_{18}| + |\tilde{a}_{8,10}^* - a_{8,10}| = 0.25, \quad
		|\tilde{a}_{19}^* - a_{19}| + |\tilde{a}_{9,10}^* - a_{9,10}| = 0.5, \\
		&|\tilde{a}_{1,10}^* - a_{1,10}| = 0.
	\end{align*}
	Thus, the largest constraint deviation occurs for the pair $(a_{17}, a_{7,10})$.
	\begin{figure}[th]
		\centering
		\begin{tikzpicture}[
			declare function={
				f1(\x) = abs(9-\x);
				f2(\x) = (\x>=1 && \x<=16) ? abs(2-\x/4) : 2*sqrt(\x)-6;
				f3(\x) = (\x>=1 && \x<=9)  ? abs(3-\x/3) : 2*sqrt(\x)-6;
				f4(\x) = (\x>=1 && \x<=25) ? abs(2-\x/5) : 2*sqrt(\x)-7;
				f5(\x) = (\x>=1 && \x<=16) ? abs(3-\x/4) : 2*sqrt(\x)-7;
				f6(\x) = (\x>=1 && \x<=36) ? abs(2-\x/6) : 2*sqrt(\x)-8;
				f7(\x) = (\x>=1 && \x<=16) ? abs(4-\x/4) : 2*sqrt(\x)-8;
				f8(\x) = (\x>=1 && \x<=16) ? abs(2-\x/4) : 2*sqrt(\x)-6;
				f9(\x) = (\x>=1 && \x<=36) ? abs(2-\x/6) : 2*sqrt(\x)-8;
				f(\x) = f1(\x) + f2(\x) + f3(\x) + f4(\x) + f5(\x) + f6(\x) + f7(\x) + f8(\x) + f9(\x);
			}
			]
			\begin{axis}[
				width=12cm,
				height=7cm,
				xmin=1, xmax=21,
				ymin=0, ymax=40, 
				xlabel={$x$}, 
				ylabel={$f(x)$}, 
				legend style={at={(0.05,0.95)}, anchor=north west, cells={anchor=west}},
				ymajorgrids=true,
				xmajorgrids=true,
				grid style=dashed,
				thick,
				xtick={1, 6, 11, 16, 21},
				extra x ticks={9},
				extra x tick labels={{\color{violet}$x_1=9$}},
				extra x tick style={tick label style={font=\bfseries}},
				extra y ticks={4.2},
				extra y tick labels={{\color{violet}$4.2$}},
				extra y tick style={tick label style={font=\bfseries}}
				]
				
				\addplot[
				color=blue,
				mark=none,
				domain=1:21,
				samples=400, 
				very thick 
				] {f(x)}; 
				\addlegendentry{$f(x)$}
				
				\addplot[
				color=red,
				dashed, 
				mark=none,
				domain=1:21,
				thick
				] {4.2};
				\addlegendentry{Global minimum value}
				
				\node[circle, fill=black, inner sep=1.5pt] at (axis cs:9, 4.2) {};
				\node[below, font=\scriptsize, fill=white, inner sep=1pt, yshift=-2pt] at (axis cs:9, 4.2) {{\color{violet}$(9, 4.2)$}};
			\end{axis}
		\end{tikzpicture}
		\caption{Graph of $f$ in $[1, 21]$ for Example 6}
		\label{example_6_figure}
	\end{figure}
	\begin{figure}[th]
		\centering
		\begin{tikzpicture}
			\begin{axis}[
				ybar=2pt,
				bar width=10pt, 
				area legend, 
				width=12cm, 
				height=7cm,
				ymin=0, 
				ymax=0.35, 
				ylabel={Optimal Weights},
				xlabel={Criteria},
				symbolic x coords={$c_1$, $c_2$, $c_3$, $c_4$, $c_5$, $c_6$, $c_7$, $c_8$, $c_9$, $c_{10}$},
				xtick=data,
				nodes near coords,
				every node near coord/.append style={
					font=\scriptsize, 
					rotate=90, 
					anchor=west, 
					/pgf/number format/fixed, 
					/pgf/number format/precision=4
				},
				legend style={at={(0.95,0.95)}, anchor=north east}, 
				enlarge x limits=0.08, 
				ymajorgrids=true,
				grid style=dashed
				]
				
				\addplot[fill=blue!60, draw=black] coordinates {
					($c_1$, 0.2586) ($c_2$, 0.1149) ($c_3$, 0.0862) ($c_4$, 0.0517) ($c_5$, 0.1149)
					($c_6$, 0.0431) ($c_7$, 0.0647) ($c_8$, 0.0647) ($c_9$, 0.1724) ($c_{10}$, 0.0287)
				};
				
				\addplot[fill=red!60, draw=black] coordinates {
					($c_1$, 0.2462) ($c_2$, 0.1094) ($c_3$, 0.0821) ($c_4$, 0.0492) ($c_5$, 0.1094)
					($c_6$, 0.0410) ($c_7$, 0.1094) ($c_8$, 0.0616) ($c_9$, 0.1642) ($c_{10}$, 0.0274)
				};
				
				\legend{$W_1^*$, $W_2^*$}
			\end{axis}
		\end{tikzpicture}
		\caption{Comparison of optimal weight sets for Example 6}
		\label{fig:bar_example6}
	\end{figure}
	\subsection{Selection strategy for multiple optimal weight sets}
	From Eq. \eqref{optimal_pcs}, it is clear that the taxicab BWM leads to multiple optimal weight sets in two scenarios: (i) the existence of multiple values for $\tilde{a}_{BW}^*$ (Examples 3, 4, and 5), and (ii) the presence of an upside criterion $i$ corresponding to $\tilde{a}_{BW}^*$ with $a_{Bi}=a_{iW}$ (Examples 2 and 6). A combination of both scenarios may also arise.\\\\
	For the first scenario with finite optimal weight sets, there is often a natural preference because the intended ordering is violated in all the other optimal weight sets (Examples 3 and 4). For the first scenario with interval weights, the optimal weight set corresponding to the midpoint of the parameter interval is a natural choice (e.g., $a=7.5$ in Example 5). However, no such straightforward selection criterion exists for the second scenario.\\\\
	Eq. \eqref{optimal_pcs} suggests that the second scenario always leads to two optimal modification strategies: one retains $a_{Bi}$ and lowers $a_{iW}$, resulting in $a_{Bi}^*>a_{iW}^*$, and the other lowers $a_{Bi}$ and retains $a_{iW}$, resulting in $a_{Bi}^*<a_{iW}^*$. Therefore, the issue can be resolved by simply asking the decision-maker an auxiliary question: ``Which is higher for you: the preference of the best criterion over criterion $i$, or the preference of criterion $i$ over the worst criterion?" If the decision-maker chooses the first option, then the optimally modified PCS with $a_{Bi}^*>a_{iW}^*$---and thus corresponding optimal weight set---is preferred over the other, and vice versa. Consequently, if the decision-maker chooses the first option, $W_1^*$ is preferred over $W_2^*$ for both Examples 2 and 6, and vice versa.
	\subsection{Comparative analysis of computation time}
	In this subsection, we compare the computational efficiency of the proposed analytical approach against a conventional optimization-based approach for calculating the optimal weights. For a fair comparison, we randomly generated $1,000$ PCSs for each value of $n=3,4,\ldots,10$ and solved them using both approaches. Both approaches were implemented in MATLAB 2018a—the proposed method as a dedicated program and the conventional method by solving optimization model \eqref{optimization_2} using the `fmincon' solver. All computations were performed on a laptop with an Intel Core i3-4005U processor and 4 GB of RAM. The results, measured in seconds, are summarized in Table \ref{comparison_1}.\\\\
	The data unequivocally demonstrates the superior time-efficiency of the proposed approach.  For every value of $n$, the total, average, maximum, and minimum computation times for our method are orders of magnitude lower than those required by the optimization software. It is essential to note that the time recorded for the conventional optimization software is the cost of computing just a single optimal solution.\\\\
	This profound disparity in computational load highlights a fundamental advantage of our method. The proposed closed-form approach, by leveraging a direct analytical formula, eliminates the need for iterative optimization routines. This not only avoids convergence issues but also reduces the computational overhead to a negligible level. Consequently, our method is exceptionally suited for large-scale, real-time applications.
	\begin{sidewaystable}[ph]
		\caption{Computation time comparison between the proposed approach and an optimization software (in seconds)}\label{comparison_1}
		\centering		
		\begin{tabular}{@{}ccccccccccc@{}}
			\toprule
			\multirow{2}{*}{$n$}&&\multicolumn{4}{c}{Proposed approach}&&\multicolumn{4}{c}{Optimization software}\\
			\cmidrule{3-6}\cmidrule{8-11}
			&&Total&Average&Maximum&Minimum&&Total&Average&Maximum&Minimum\\
			\midrule
			$3$&&$1.0400$&$0.0010$&$0.0236$&$0.0004$&&$52.8747$&$0.0529$&$2.7555$&$0.0232$\\
			$4$&&$1.3264$&$0.0013$&$0.0241$&$0.0005$&&$139.7537$&$0.1398$&$2.7349$&$0.0257$\\
			$5$&&$1.4101$&$0.0017$&$0.0311$&$0.0005$&&$104.4773$&$0.1045$&$3.2969$&$0.0278$\\
			$6$&&$1.6038$&$0.0016$&$0.0447$&$0.0006$&&$154.8411$&$0.1548$&$3.7520$&$0.0334$\\
			$7$&&$1.6781$&$0.0017$&$0.0359$&$0.0005$&&$175.1214$&$0.1751$&$3.7821$&$0.0369$\\
			$8$&&$1.7970$&$0.0018$&$0.0450$&$0.0006$&&$191.8566$&$0.1919$&$3.5184$&$0.0423$\\
			$9$&&$2.0698$&$0.0021$&$0.0678$&$0.0007$&&$246.0111$&$0.2460$&$2.9385$&$0.0488$\\
			$10$&&$1.8879$&$0.0019$&$0.0421$&$0.0006$&&$280.2735$&$0.2803$&$3.4931$&$0.0526$\\		
			\bottomrule				
		\end{tabular}
		\par\smallskip \small \textit{The proposed approach reduces computation time by orders of magnitude compared to optimization software.}
	\end{sidewaystable}
	\subsection{Comparison between the nonlinear BWM and the taxicab BWM}
	In this subsection, we compare the outcomes of the nonlinear BWM and the taxicab BWM in particular situations.\\\\
	Consider a decision problem with the set of criteria $C=\{c_1,c_2,c_3,c_4\}$, where $c_1$ and $c_4$ are the best and the worst criterion respectively. The best-to-other vector is $A_{B}=(1,2,4,8)$ and the other-to-worst vector is $A_{W}=(8,4,2,1)^T$. Note that $(A_{B},A_{W})$ is consistent. So, both the nonlinear BWM and the taxicab BWM give the same optimal weight set $W=\{0.5333,0.2667,0.1333,0.0667\}$ as their unique solution.\\\\
	Now, suppose an additional criterion $c_5$ that is neither best nor worst, with pairwise comparisons $(a_{15},a_{54})=(2,2)$ is included in the decision process. Then the revised PCS $(A_{B}',A_{W}')=((1,2,4,8,2),(8,4,2,1,2)^T)$ becomes inconsistent. For this revised PCS, the nonlinear BWM produces multiple optimal solutions, yielding the optimal interval-weights $w_1=[0.4074,0.4605]$, $w_2=[0.1998,0.2725]$, $w_3=[0.0900,0.1340]$, $w_4=[0.0558,0.0631]$, and $w_5=[0.1508,0.1704]$.\cite{wu2023analytical} We then obtain a unique weight set $W_1=\{0.4354,0.2315,0.1122,0.0597,0.1612\}$ by introducing a secondary objective function.\cite{wu2023analytical} On the other hand, the taxicab BWM gives a unique optimal weight set $W_2=\{0.4487,0.2244,0.1122,0.0561,0.1586\}$. The modified PCSs associated with $W_1$ and $W_2$ are given in Table \ref{comparison}.\\\\
	The results suggest that introducing a single inconsistent criterion into an otherwise consistent PCS can lead to multiple weight sets in the nonlinear BWM model. To derive a unique solution, the nonlinear approach employs a secondary objective function, which modifies all seven original comparison values—indicating a global redistribution of inconsistency that distorts even initially consistent judgments. In contrast, the taxicab BWM selectively adjusts only two inconsistent comparisons, preserving the remaining structure. This implies that the nonlinear BWM modifies all comparison values more extensively, which may not always be desirable, whereas the taxicab BWM adopts a more targeted approach by resolving inconsistencies without compromising the integrity of the consistent judgments.\\\\ The targeted approach applies extensive modifications to comparisons where inconsistencies are highly concentrated, while making minimal adjustments where inconsistencies are less pronounced. This adaptive modification strategy enables effective local consistency restoration.\\\\
	Based on this analysis, we propose the following model selection strategy for the decision-makers:
		\begin{itemize}
			\item When inconsistency is concentrated in a few comparisons, prioritize the taxicab BWM as its targeted nature maintains the integrity of the consistent preferences.
			\item When inconsistency is evenly distributed across all comparisons, prioritize the nonlinear BWM as its global modification approach adjusts all comparisons uniformly.
	\end{itemize}
	 
	\begin{table}[t!]
		\caption{Comparison between the nonlinear BWM and the taxicab BWM}\label{comparison}
		\centering		
		\begin{tabular}{@{}cccc@{}}
			\toprule
			\multirow{2}{*}{Comparison}&Original&\multicolumn{2}{c}{Optimally modified comparison values}\\
			\cmidrule{3-4}
			&comparison value&Nonlinear BWM\cite{wu2023analytical}&Taxicab BWM\\		
			\midrule
			$a_{12}$&$2$&$1.8807$&$2$\\
			$a_{13}$&$4$&$3.8807$&$4$\\
			$a_{14}$&$8$&$7.2984$&$8$\\
			$a_{15}$&$2$&$2.7016$&$2.8284$\\
			$a_{24}$&$4$&$3.8807$&$4$\\
			$a_{34}$&$2$&$1.8807$&$2$\\
			$a_{54}$&$2$&$2.7016$&$2.8284$\\			
			\bottomrule				
		\end{tabular}
		\par\smallskip \small \textit{Taxicab BWM selectively adjusts only inconsistent inputs, preserving consistent ones, while nonlinear BWM modifies all values.}
	\end{table}
	 
	\subsection{Real-world scenario: smartphone selection}
	In this subsection, we rank the criteria for purchasing a smartphone in order of their criticality.\\\\
	We consider $30$ criteria divided across $6$ categories as shown in Table \ref{criteria}. The decision-maker provided the pairwise comparisons between categories and between the criteria within the same category, which are given in Table \ref{comparison_values}. We first performed consistency analysis, which confirmed that all PCSs have CR values below their respective thresholds and are therefore admissible. We then calculated the category weights and local weights of the criteria using the proposed framework. For each instance, except for the local weights of criteria from the category `core technical specifications' and the category `build, design \& durability,' we obtained a unique optimal weight set.\\\\
	For the criteria from the category `core technical specifications' $(c_{21}-c_{27})$, we obtained two optimal weight sets as follows:
		\begin{enumerate}
			\item $W_1^*=\{0.2378,0.0991,0.1329,0.2973,0.0743,0.0595,0.0991\}$
			\item $W_2^*=\{0.2204,0.1653,0.1232,0.2754,0.0689,0.0551,0.0918\}$.
		\end{enumerate}
		For the auxiliary question \textbf{Which is higher: the preference of `battery capacity' over `RAM capacity', or the preference of `RAM capacity' over `operating system'}, the decision-maker responded \textbf{the preference of `RAM capacity' over `operating system'}. Thus, in this scenario, $W_2^*$ is preferred over $W_1^*$.\\\\
For the criteria from the category `build, design \& durability' $(c_{51}-c_{54})$, we obtained two optimal weight sets as follows:
	\begin{enumerate}
		\item $W_1^*=\{0.2308,0.3846,0.3077,0.0769\}$
		\item $W_2^*=\{0.1429,0.4286,0.3429,0.0857\}$.
	\end{enumerate}
	For the auxiliary question \textbf{Which is higher: the preference of `weight \& ergonomics' over `build materials', or the preference of `build materials' over `biometric security'}, the decision-maker responded \textbf{the preference of `weight \& ergonomics' over `build materials'}. Thus, in this scenario, $W_2^*$ is preferred over $W_1^*$.\\\\
The category weights, the global and local weights of the criteria (the preferred one in cases of multiple solutions), and the final criteria ranking are provided in Table \ref{weights}.\\\\
Our analysis revealed that `5G \& network bands' is the most critical criterion, while `optical zoom' is the least critical criterion among those ranked.
\begin{table}[t!]
	\caption {List of criteria for smartphone selection}\label{criteria}
	\centering
	\begin{tabular}{lll}
			\toprule
			Category&&Criteria\\
			\midrule
			$c_1$: Financial \& commercial&&$c_{11}$: Base price\\
			&&$c_{12}$: Cost of accessories\\
			&&$c_{13}$: Warranty \& repair\\
			&&$c_{14}$: Resale value\\
			$c_2$: Core technical specifications&&$c_{21}$: Processor (CPU/GPU)\\
			&&$c_{22}$: RAM capacity\\
			&&$c_{23}$: Internal storage\\
			&&$c_{24}$: Battery capacity\\
			&&$c_{25}$: Charging speed\\
			&&$c_{26}$: Operating system\\
			&&$c_{27}$: Thermal management\\
			$c_3$: Camera \& optics&&$c_{31}$: Primary/rear camera\\
			&&$c_{32}$: Front/selfie camera\\
			&&$c_{33}$: Video recording\\
			&&$c_{34}$: Low light performance\\
			&&$c_{35}$: Optical zoom\\
			&&$c_{36}$: Ultrawide capability\\
			$c_4$: Display \& multimedia&&$c_{41}$: Screen size \& resolution\\
			&&$c_{42}$: Display technology\\
			&&$c_{43}$: Refresh rate\\
			&&$c_{44}$: Peak brightness\\
			&&$c_{45}$: Audio quality\\
			$c_5$: Build, design \& durability&&$c_{51}$: Build materials\\
			&&$c_{52}$: Weight \& ergonomics\\
			&&$c_{53}$: IP rating\\
			&&$c_{54}$: Biometric security\\
			$c_6$: Connectivity \& network&&$c_{61}$: 5G \& network bands\\
			&&$c_{62}$: Wi-fi \& bluetooth\\
			&&$c_{63}$: NFC support\\
			&&$c_{64}$: SIM capabilities\\
			\bottomrule
	\end{tabular}
\end{table}
\begin{table}[t!]
	\caption {Pairwise comparison values provided by the decision-maker}\label{comparison_values}
	\centering
	\begin{subtable}{1\textwidth}
		\caption {Comparison values for categories}
		\centering
		\begin{tabular}{ccc}
				\toprule
				Best&Worst&\multirow{2}{*}{PCS}\\
				category&category&\\
				\midrule
				$c_2$&$c_3$&$A_B=(3,1,7,5,6,2)$, $A_W=(2,7,1,2,4,6)^T$\\
				\bottomrule
		\end{tabular}
	\end{subtable}\\
	\vspace{0.5cm}
	\begin{subtable}{1\textwidth}
		\caption {Comparison values for the criteria within the same category}
		\centering
		\begin{tabular}{cccc}
				\toprule
				\multirow{2}{*}{Criteria}&Best&Worst& \multirow{2}{*}{PCS}\\
				&criteria&criteria&\\
				\midrule
				$c_{11}$-$c_{14}$&$c_{11}$&$c_{14}$&$A_B=(1,7,6,9)$, $A_W=(9,3,3,1)^T$\\
				$c_{21}$-$c_{27}$&$c_{24}$&$c_{26}$&$A_B=(1,3,2,1,4,5,3)$, $A_W=(4,3,2,5,2,1,2)^T$\\
				$c_{31}$-$c_{36}$&$c_{31}$&$c_{35}$&$A_B=(1,2,3,4,7,5)$, $A_W=(7,6,4,3,1,2)^T$\\
				$c_{41}$-$c_{45}$&$c_{43}$&$c_{44}$&$A_B=(3,2,1,8,4)$, $A_W=(6,5,8,1,6)^T$\\
				$c_{51}$-$c_{54}$&$c_{52}$&$c_{54}$&$A_B=(3,1,1,5)$, $A_W=(3,5,4,1)^T$\\
				$c_{61}$-$c_{64}$&$c_{61}$&$c_{63}$&$A_B=(1,2,7,3)$, $A_W=(7,6,1,5)^T$\\
				\bottomrule	
		\end{tabular}
	\end{subtable}
\end{table}
 
\begin{table}[t!]
	\caption {Category weights, criteria weights, and criteria ranking \label{weights}}
	\centering
	{\small\begin{tabular}{@{}cccccc@{}}
			\toprule
			Category &Weight&Criteria&Local weight&Global weight&Rank\\
			\midrule
			$c_1$&$0.1235$&$c_{11}$&$0.7039$&$0.0869$&$4$\\
			&&$c_{12}$&$0.1006$&$0.0124$&$22$\\
			&&$c_{13}$&$0.1173$&$0.0145$&$20$\\
			&&$c_{14}$&$0.0782$&$0.0097$&$24$\\
			&&$\epsilon^*$&$3.2143$&-&-\\
			&&CR&$0.2009$&-&-\\
			\midrule
			$c_2$&$0.3704$&$c_{21}$&$0.2204$&$0.0816$&$6$\\
			&&$c_{22}$&$0.1653$&$0.0612$&$7$\\
			&&$c_{23}$&$0.1232$&$0.0456$&$8$\\
			&&$c_{24}$&$0.2754$&$0.1020$&$2$\\
			&&$c_{25}$&$0.0689$&$0.0255$&$11$\\
			&&$c_{26}$&$0.0551$&$0.0204$&$14$\\
			&&$c_{27}$&$0.0918$&$0.0340$&$9$\\
			&&$\epsilon^*$&$7.6388$&-&-\\
			&&CR&$0.3819$&-&-\\
			\midrule
			$c_3$&$0.0529$&$c_{31}$&$0.3310$&$0.0175$&$15$\\
			&&$c_{32}$&$0.2837$&$0.0150$&$19$\\
			&&$c_{33}$&$0.1891$&$0.0100$&$23$\\
			&&$c_{34}$&$0.0827$&$0.0044$&$27$\\
			&&$c_{35}$&$0.0473$&$0.0025$&$30$\\
			&&$c_{36}$&$0.0662$&$0.0035$&$28$\\
			&&$\epsilon^*$&$3.9333$&-&-\\
			&&CR&$0.1639$&-&-\\
			\midrule
			$c_4$&$0.0741$&$c_{41}$&$0.2308$&$0.0171$&$16$\\
			&&$c_{42}$&$0.1923$&$0.0143$&$21$\\
			&&$c_{43}$&$0.3077$&$0.0228$&$12$\\
			&&$c_{44}$&$0.0385$&$0.0029$&$29$\\
			&&$c_{45}$&$0.2308$&$0.0171$&$17$\\
			&&$\epsilon^*$&$4.7333$&-&-\\
			&&CR&$0.2254$&-&-\\
			\midrule
			$c_5$&$0.0617$&$c_{51}$&$0.1429$&$0.0088$&$25$\\
			&&$c_{52}$&$0.4286$&$0.0264$&$10$\\
			&&$c_{53}$&$0.3429$&$0.0212$&$13$\\
			&&$c_{54}$&$0.0857$&$0.0053$&$26$\\
			&&$\epsilon^*$&$1.5833$&-&-\\
			&&CR&$0.1979$&-&-\\
			\midrule
			$c_6$&$0.3175$&$c_{61}$&$0.3684$&$0.1170$&$1$\\
			&&$c_{62}$&$0.3158$&$0.1003$&$3$\\
			&&$c_{63}$&$0.0526$&$0.0167$&$18$\\
			&&$c_{64}$&$0.2632$&$0.0836$&$5$\\
			&&$\epsilon^*$&$2.4333$&-&-\\
			&&CR&$0.2433$&-&-\\
			\midrule
			$\epsilon^*$&$4.6$&-&-&-&-\\
			CR&$0.1917$&-&-&-&-\\
			\bottomrule
	\end{tabular}}
	\par\smallskip \small \textit{`5G \& network bands $(c_{61})$' and `optical zoom $(c_{35})$' are the most and least critical criteria, respectively.}	
\end{table}

\section{Conclusions and Future Directions}
The BWM is a recent MCDM method that has been effectively applied to numerous real-world applications, drawing significant attention from researchers. In this paper, we propose an analytical framework for a model of BWM called taxicab BWM by formulating an equivalent optimal modification based model. We develop an algorithm to obtain the optimal weights, illustrate its effectiveness through numerical examples, and demonstrate its applicability using a real-world smartphone selection scenario.\\\\
This research provides fundamental clarifications to the theory of the taxicab BWM. A central contribution is the constructive characterization of its solution space. While previous studies reported a unique solution, our analytical approach reveals and formally explains the conditions under which the model yields multiple optimal weight sets—whether finitely many (Examples 2, 3, 4, and 6) or infinitely many (Example 5). In such instances, determining the exact number of optimal weight sets and obtaining them all numerically through optimization software can be challenging, particularly when there are finitely many due to the discrete nature of the solution space. In this research, we analytically derive all possible optimal weight sets, thereby eliminating the need for optimization software while significantly enhancing computational efficiency over the traditional approach, making the framework more accessible to non-experts and suitable for large-scale, dynamic decision-making environments. Additionally, we develop a decision-maker-aided strategy to select the most appropriate solution when multiple optimal weight sets arise. Based on this framework, we formulate a mixed-integer optimization model to compute the values of CI. We also establish threshold values for the consistency ratio, providing a practical criterion for assessing the admissibility of PCSs in real-world applications.\\\\
This analytical framework also reveals a inherent limitation of the taxicab BWM. By its nature of minimizing the TD, the model does not explicitly control the worst local violation. While the targeted approach of the taxicab BWM is advantageous in certain contexts, it is not always preferable.\\\\
Building on this analytical foundation, several key research directions emerge. First, it would be valuable to develop a hybrid model that integrates the inconsistency localization capability of the taxicab BWM with the uniform modification approach of the nonlinear BWM. Second, formulating an analytical framework for the Euclidean BWM to study its solution uniqueness and properties would provide important theoretical insights.
\section*{Acknowledgment}
The first author would like to express his gratitude to the Council of Scientific and Industrial Research (CSIR), India for providing the fellowship that supported this research work.
\section*{Declaration of Conflict of Interest}
The authors declare that they have no known conflict of financial interests or personal relationships that could have appeared to influence the work reported in this paper.

\section*{ORCID}
\noindent Harshit M. Ratandhara - \url{https://orcid.org/0009-0007-1791-1475}

\noindent Mohit Kumar - \url{https://orcid.org/0000-0001-7609-3360}

\section{Appendix}

\textbf{Proof of Proposition \ref{geq1}:}\\\\
Since $(\tilde{A}_{B},\tilde{A}_{W})$ is consistent, we have $\tilde{a}_{Bi}\times\tilde{a}_{iW}=\tilde{a}_{BW}<1$ for all $i\in D$. Also, $a_{Bi},a_{iW}\geq 1$ gives $a_{Bi}\times a_{iW}\geq1$. This implies $\tilde{a}_{Bi}\times\tilde{a}_{iW}<a_{Bi}\times a_{iW}$. Let $|\tilde{a}_{Bi}-a_{Bi}|=\zeta_{Bi}$ and $|\tilde{a}_{iW}-a_{iW}|=\zeta_{iW}$. Then there are four cases:
\begin{enumerate}
	\item $\tilde{a}_{Bi}=a_{Bi}+\zeta_{Bi}$, $\tilde{a}_{iW}=a_{iW}+\zeta_{iW}$\\
	Since $\zeta_{Bi},\zeta_{iW}\geq 0$, we get $\tilde{a}_{Bi}\geq a_{Bi}$ and $\tilde{a}_{iW}\geq a_{iW}$. This gives $\tilde{a}_{Bi}\times\tilde{a}_{iW}\geq a_{Bi}\times a_{iW}$, which is not possible.
	\item $\tilde{a}_{Bi}=a_{Bi}+\zeta_{Bi}$, $\tilde{a}_{iW}=a_{iW}-\zeta_{iW}$\\
	In this case, we have $\tilde{a}_{Bi}\geq a_{Bi}$, which implies $a_{Bi}\times a_{iW}\leq \tilde{a}_{Bi}\times a_{iW}$. Take $\tilde{a}_{Bi}'=\tilde{a}_{Bi}$ and $\tilde{a}_{iW}'=\frac{1}{\tilde{a}_{Bi}}$. So, $|\tilde{a}_{Bi}'-a_{Bi}|=|\tilde{a}_{Bi}-a_{Bi}|$. Note that $\tilde{a}_{Bi}\times\tilde{a}_{iW}<1=\tilde{a}_{Bi}\times \tilde{a}_{iW}'\leq a_{Bi}\times a_{iW}\leq \tilde{a}_{Bi}\times a_{iW}$. This gives $\tilde{a}_{iW}<\tilde{a}_{iW}'\leq a_{iW}$. So, we get $|\tilde{a}_{iW}'-a_{iW}|=a_{iW}-\tilde{a}_{iW}'<a_{iW}-\tilde{a}_{iW}=|\tilde{a}_{iW}-a_{iW}|$.
	\item $\tilde{a}_{Bi}=a_{Bi}-\zeta_{Bi}$, $\tilde{a}_{iW}=a_{iW}+\zeta_{iW}$\\
	Take $\tilde{a}_{Bi}'=\frac{1}{\tilde{a}_{iW}}$ and $\tilde{a}_{iW}'=\tilde{a}_{iW}$. By reasoning similarly to 2, we obtain $|\tilde{a}_{Bi}'-a_{Bi}|<|\tilde{a}_{Bi}-a_{Bi}|$ and $|\tilde{a}_{iW}'-a_{iW}|=|\tilde{a}_{iW}-a_{iW}|$.
	\item $\tilde{a}_{Bi}=a_{Bi}-\zeta_{Bi}$, $\tilde{a}_{iW}=a_{iW}-\zeta_{iW}$\\
	If $\tilde{a}_{Bi}\times a_{iW}>1$, then take $\tilde{a}_{Bi}'=\tilde{a}_{Bi}$ and $\tilde{a}_{iW}'=\frac{1}{\tilde{a}_{Bi}}$. By arguing similarly to 2, we get $|\tilde{a}_{Bi}'-a_{Bi}|=|\tilde{a}_{Bi}-a_{Bi}|$ and $|\tilde{a}_{iW}'-a_{iW}|<|\tilde{a}_{iW}-a_{iW}|$. If $\tilde{a}_{Bi}\times a_{iW}\leq 1$, then take $\tilde{a}_{Bi}'=\frac{1}{a_{iW}}$ and $\tilde{a}_{iW}'=a_{iW}$. So, $|\tilde{a}_{iW}'-a_{iW}|=0\leq |\tilde{a}_{iW}-a_{iW}|$. Now, $\tilde{a}_{Bi}\times a_{iW}\leq 1 =\tilde{a}_{Bi}'\times a_{iW}\leq a_{Bi}\times a_{iW}$ implies $\tilde{a}_{Bi}\leq \tilde{a}_{Bi}'\leq a_{Bi}$, which gives $|\tilde{a}_{Bi}'-a_{Bi}|=a_{Bi}-\tilde{a}_{Bi}'\leq a_{Bi}-\tilde{a}_{Bi}=|\tilde{a}_{Bi}-a_{Bi}|$. 
\end{enumerate}
Now, take $\tilde{a}_{BW}'=1$. Since $\tilde{a}_{Bi}'\times \tilde{a}_{iW}'=1$, $(\tilde{A}_{B}',\tilde{A}_{W}')$ is consistent. Also, $\tilde{a}_{BW}<1=\tilde{a}_{BW}'\leq a_{BW}$ gives $|\tilde{a}_{BW}'-a_{BW}|=a_{BW}-\tilde{a}_{BW}'< a_{BW}-\tilde{a}_{BW}=|\tilde{a}_{BW}-a_{BW}|$, which completes the proof.\\\\
\textbf{Proof of Proposition \ref{min_1}:}\\\\
There are 8 possibilities for $(x,y,z)\in \mathbb{R}^3$ such that $(a+x)\times(b+y)=c+z$:\\
\begin{enumerate*}[series=possibilities]
	\item $x\geq 0$, $y\leq 0$, $z\geq 0$\hspace{1cm}
	\item $x\leq 0$, $y\leq 0$, $z\geq 0$\hspace{1cm}
	\item $x\leq 0$, $y\geq 0$, $z\geq 0$\\
\end{enumerate*}
\begin{enumerate*}[resume=possibilities]
	\item $x\geq 0$, $y\geq 0$, $z\geq 0$\hspace{1cm}
	\item $x\leq 0$, $y\leq 0$, $z\leq 0$\hspace{1cm}
	\item $x\geq 0$, $y\leq 0$, $z\leq 0$\\
\end{enumerate*}
\begin{enumerate*}[resume=possibilities]
	\item $x\leq 0$, $y\geq 0$, $z\leq 0$\hspace{1cm}
	\item $x\geq 0$, $y\geq 0$, $z\leq 0$.\\
\end{enumerate*}\\
Here, we shall prove that for possibilities 1 to 7, statement (2) holds. Out of these seven possibilities, we shall discuss only possibility 1, 2, 4 and 5 as for the other possibilities, proof is similar to one of these four possibilities.\\\\
Possibility 1: Here, we have $x\geq 0$, $y\leq 0$, $z\geq 0$. If $y=z=0$, then statement (1) holds. Now, consider the case that at least one of $y$ and $z$ is non-zero. So, we get $(a+x)\times b - c>0$. Let $x'$ be such $(a+x')\times b-c=0$. This gives $0<x'<x$, and so, $|x'|<|x|$. Taking $y'=z'=0$, we get $(a+x')\times(b+y')=c+z'$ and $|x'|+|y'|+|z'|<|x|+|y|+|z|$.\\\\
Possibility 2: Here, we have $x\leq 0$, $y\leq 0$, $z\geq 0$. If $a+x\geq0$, then $b+y\geq0$. Since $(a+x)\times (b+y)-c-z=0$, we get $a\times b\geq c$, which is contradiction. So, $a+x<0$ and $b+y<0$. Take $x''=-2a-x$ and $y''=-2b-y$. Now, it is sufficient to check $|x''|< |x|$, $|y''|<|y|$, $a+x'',b+y''>0$ and $(a+x'')\times (b+y'')=c+z$, i.e., this possibility can be transformed into one of the possibility 1, 3 or 4. Note that $a+x''=-(a+x)$ and $b+y''=-(b+y)$. So, we get $a+x'',b+y''>0$ and $(a+x'')\times (b+y'')=c+z$. Now, observe that $|x|=-x$ and $|x''|=\begin{cases}
	-2a-x\ \text{ if } x<-2a,\\
	2a+x\quad  \text { if } x\geq -2a.	
\end{cases}$ For $x<-2a$, we get $0\leq -2a-x<-x$, and for $x\geq -2a$, $x+a<0$ implies $2a+x<-x$. This gives $|x''|<|x|$. Similarly, it follows that $|y''|<|y|$. \\\\
Possibility 4: Here, we have $x\geq 0$, $y\geq 0$, $z\geq 0$. If $z=0$, then statement (1) holds. Now, consider the case $z\neq 0$. Then we get $(a+x)\times (b+y)-c>0$. If $a\times (b+y)-c\geq 0$, then take $x'=z'=0$, and let $y'$ be such that $a\times (b+y')-c=0$. Then $0<y'\leq y$. If $a\times (b+y)-c< 0$, then take $y'=y$, $z'=0$, and let $x'$ be such that $(a+x')\times (b+y)-c=0$. Then $0<x'<x$. Observe that, in either case, we get $x',y'\geq 0$, $z'\leq 0$, $(a+x')\times(b+y')=c+z'$ and $|x'|+|y'|+|z'|<|x|+|y|+|z|$. \\\\
Possibility 5: Here, we have $x\leq 0$, $y\leq 0$, $z\leq 0$. If $x=y=0$, then statement (1) holds. Now, consider the case that at least one of $x$ and $y$ is non-zero. If $c+z=0$, then $|z|=c$. Take $x'=y'=0$ and $z'=a\times b-c$. So, we get $(a+x')\times(b+y')=c+z'$ and $|x'|+|y'|+|z'|=c-a\times b<c=|z|\leq |x|+|y|+|z|$. Thus, we are done. Now, assume that $c+z\neq 0$. This implies $a+x\neq 0$ and $b+y\neq 0$. If $a+x, b+y >0$, then $a\times b-c-z>0$. Let $z'$ be such that $a\times b-c-z'=0$. Then $z<z'<0$. This given $|z'|<|z|$. Take $x'=y'=0$. So, we get $(a+x')\times(b+y')=c+z'$ and $|x'|+|y'|+|z'|<|x|+|y|+|z|$. If $a+x,b+y<0$, then it suffices to prove that there exist $|x''|+|y''|<|x|+|y|$, $a+x'',b+y''>0$ and $(a+x'')\times (b+y'')=c+z$. Take $x''=-2a-x$ and $y''=-2b-y$. By possibility 2, $x''$ and $y''$ satisfy all the requirement. Hence the proof.\\\\
\textbf{Proof of Theorem \ref{a=b}:}\\\\
First, assume that $z=0$. This gives $x\neq y\neq x'$. Without loss of generality, we may assume that $x<y$. Consider $f(w)= (a+w)\times (a+w)-c$, $w\in[0,\infty)$. Note that $f$ strictly increasing and $f(x)<0$. Also, $f(\frac{x+y}{2})=(a+\frac{x+y}{2})\times (a+\frac{x+y}{2})=a^2+a(x+y)+\frac{1}{4}(x+y)^2-c=\frac{1}{4}(x+y)^2-xy=\frac{1}{4}(x-y)^2$. Now, $x\neq y$ gives $(x-y)^2>0$, and so, $f(\frac{x+y}{2})>0$. Since $f$ is strictly increasing and $f(x')=0$, we get $x'<\frac{x+y}{2}$, i.e., $2x'<x+y=x+y+z$.\\\\
Now, assume that $z\neq 0$. Observe that $\max\{a+x,a+y\}\geq 1$.\\ 
Case 1: Let $\max\{a+x,a+y\}>1$. Then, without loss of generality, we may assume that $a+x>1$. To prove Theorem, it is sufficient to prove that there exist $x'',y''\geq 0$ such that $(a+x'')\times (a+y'')=c$ and $x''+y''<x+y+z$. We have $(a+x)\times (a+y)-c+z=0$. So, we get $(a+x)\times (a+y)-c+(a+x)z>0$. This implies $(a+x)\times (a+y+z)-c>0$. Let $y''$ be such that $(a+x)\times (a+y'')-c=0$. Since $(a+x)\times a<c$, we get $y''>0$. Now, $(a+x)\times (a+y)-c+(a+x)z>0$ gives $y''<y+z$. Take $x''=x$. So, we get $(a+x'')\times (a+y'')=c$. Also, $y''<y+z$ implies $x''+y''<x+y+z$.\\
Case 2: Let $\max\{a+x,a+y\}=1$. So, we get $a=1$, $x=y=0$, $z=c-1$ and $x'=\sqrt{c}-1$. We also get $c>1$, which gives $(\sqrt{c}-1)^2>0$. Thus, $2\sqrt{c}-2<c-1$, i.e., $2x'<z=x+y+z$. This completes the proof.\\\\
\textbf{Proof of Theorem \ref{a<b}:}\\\\
First, assume $b\geq \sqrt{c}$. Let $(x,y,z)\neq (x',0,0)$ be such that $x,y,z\geq 0$ and $(a+x)\times (b+y)=c-z$. So, at least one of $y$ and $z$ is non-zero. Suppose, if possible, $a+x>b$. Then we get $(a+x)\times b>c$, which is not possible. Also, if $a+x=b$, then $y=z=0$, which is not possible. So, we have $a+x<b$. Now, $(a+x)\times (b+y)-c+z=0$, along with $b>1$, implies $b(a+x)+by-c+bz>0$. This gives $(a+x+y+z)\times b-c>0$. Thus, we get $x'<x+y+z$.\\\\
Now, assume that $b< \sqrt{c}$. Let $(x,y,z)\neq (b-a+y',y',0)$ be such that $x,y,z\geq 0$ and $(a+x)\times (b+y)=c-z$.\\
Case 1: Let $a+x\geq b$. Then $a+x=b+d$ for some $d\geq 0$. This gives $(b+d)\times(b+y)=c-z$. From hypothesis, if $y=y'$ and $z=0$, then $x\neq b-a+y'$, i.e., $d\neq y'$. This implies $(d,y,z)\neq (y',y',0)$. So, by Theorem \ref{a=b}, we get $2y'<d+y+z$. This gives $b-a+2y'<x+y+z$.\\
Case 2: Let $a+x<b$. It is sufficient to prove $x'<x+y+z$ as $b<\sqrt{c}$ implies $(a+x')> b$ and so, from Case 1, we get $b-a+2y'<x'<x+y+z$. Here, we have $(a+x)\times(b+y)-c+z=0$. Now, $a+x<b$ implies $(a+x+y+z)\times b-c>0$. Thus, we get $x'<x+y+z$. This completes the proof.\\\\
\textbf{Proof of Proposition \ref{min_2}:}\\\\
The proof is similar to the proof of Proposition \ref{min_1}.\\\\
\textbf{Proof of Theorem \ref{a<=b}:}\\\\
First, consider the case $a<b$. Let $(x,y,z)\neq (x',0,0)$ be such that $x,y,z\geq 0$, $a-x,b-y>0$ and $(a-x)\times(b-y)=c+z$. So, at least one of $y$ and $z$ is non-zero. Now, $(a-x)\times (b-y)-(c+z)=0$, along with $b>1$, gives $(a-(x+y+z))\times b-c<0$. So, we get $x'<x+y+z$.\\\\
Now, consider the case $a=b$. Let $(x,y,z)\neq (x',0,0)\neq (0,x',0)$ be such that $x,y,z\geq 0$, $a-x,a-y>0$ and $(a-x)\times(a-y)=c+z$.\\\\
First, assume that $z=0$. This implies $x,y\neq 0$. Here, we have $(a-x)\times(a-y)=c$. We also have $(a-x')\times a=c$. This gives $(a-x)\times (a-y)=(a-x')\times a$. So, $-a(x+y)+xy=-ax'$. Since $x,y\neq 0$, we get $-a(x+y)<-ax'$, i.e., $x'<x+y=x+y+z$.\\\\
Now, assume that $z\neq 0$. To prove result, it is sufficient to prove that there exist $x'',y''\geq 0$ such that $a-x'',a-y''>0$, $(a-x'')\times (a-y'')=c$ and $x''+y''<x+y+z$. Since $a\leq c$, we have $a-x,a-y\leq c$. If $a-x,a-y\leq 1$, then $(a-x)\times (a-y)\leq 1\leq c<c+z$, which is contradiction. So, at least one of $a-x$ and $a-y$ is greater than $1$. Without loss of generality, we may assume that $a-x>1$. Now, $(a-x)\times (a-y)-c-z=0$ implies $(a-x)\times (a-y)-c-(a-x)z<0$, i.e., $(a-x)\times (a-y-z)-c<0$. Let $y''$ be such that $(a-x)\times (a-y'')-c=0$. Since $(a-x)\times (a-y)-(c+z)=0$, we get $0\leq y<y''$. Also, $(a-x)\times (a-y-z)-c<0$ gives $y''<y+z$. Take $x''=x$. So, $a-x=a-x''>0$, and consequently, $(a-y'')>0$. Also, $(a-x'')\times (a-y'')=c$ and $x''+y''<x+y+z$. Hence the proof.\\\\
\textbf{Proof of Theorem \ref{minima}:}\\\\
If $x_0=1$, then $[1,x_0]=\{x_0\}$. So, $\displaystyle\min_{x\in [1,x_0]}f(x)=f(x_0)$. Now, assume that $x_0\neq 1$. Observe that $[1,x_0]\subset [1,a_{Bi}\times a_{iW}]$ and $[1,x_0]\subset[1,a_{BW}]$ for all $i\in D$. So, $f_i(x)=\begin{cases}
	a_{iW}-\frac{x}{a_{Bi}}\ \text{ if } a_{iW}\leq a_{Bi}\\
	a_{Bi}-\frac{x}{a_{iW}}\ \text{ if } a_{Bi}\leq a_{iW}
\end{cases}$ and $f_{B}(x)=a_{BW}-x$, $1\leq x\leq x_0$, for all $i\in D$. Thus, $f(x)=\displaystyle\sum_{\substack{i\in D\\a_{iW}\leq a_{Bi}}}\bigg(a_{iW}-\frac{x}{a_{Bi}}\bigg)+\displaystyle\sum_{\substack{i\in D\\a_{Bi}\leq a_{iW}}}\bigg(a_{Bi}-\frac{x}{a_{iW}}\bigg)+a_{BW}-x$, i.e., $f(x)$ is of the form $bx+c$, where $b\in \mathbb{R}_{<0}$ and $c\in \mathbb{R}_{>0}$, for $1\leq x\leq x_0$. Thus, $f'(x)=b<0$ for $x\in (1,x_0)$. So, $f$ is strictly decreasing in $(1,x_0)$. Since $f$ is continuous, we get $\displaystyle\min_{x\in [1,x_0]}f(x)=f(x_0)$. Thus, in either case, we get $\displaystyle\min_{x\in [1,x_0]}f(x)=f(x_0)$. Furthermore, $x_0$ is the only point in $[1,x_0]$ at which $f$ attains this minimum value.\\\\
By similar argument, it can be proven that $f$ is strictly increasing in $[x_m,\infty)$. So, $\displaystyle\min_{x\in [x_m,\infty)}f(x)=f(x_m)$, and $x_m$ is the only point in $[x_m,\infty)$ at which $f$ attains this minimum value.\\\\
Fix $j\in\{1,2,\ldots,m\}$. Now, the fact that $[x_{j-1},x_j]$ is either subset of $[1,a_{Bi}\times a_{iW}]$, $[a_{Bi}\times a_{iW},\max\{a_{Bi}^2,a_{iW}^2\}]$, or $[\max\{a_{Bi}^2,a_{iW}^2\},\infty)$ implies that $f_i(x)$ is of the form $a\sqrt{x}+bx+c$ for $x_{j-1}\leq x\leq x_j$, where $a\in \mathbb{R}_{\geq 0}$, $b$, $c\in \mathbb{R}$. Similarly, the fact that $[x_{j-1},x_j]$ is either subset of $[1,a_{BW}]$ or $[a_{BW},\infty)$ implies that $f_{B}(x)$ is of the form $bx+c$ for $x_{j-1}\leq x\leq x_j$, where $b,c\in \mathbb{R}$. Thus, $f(x)$ is of the form $a\sqrt{x}+bx+c$ for $x_{j-1}\leq x\leq x_j$, where $a\in \mathbb{R}_{\geq 0}$, $b,c\in \mathbb{R}$. So, $f'(x)=\frac{a}{2\sqrt{x}}+b$, $x_{j-1}<x<x_j$. If $a=b=0$, then $f$ is constant on $[x_{j-1},x_j]$. So, $\displaystyle\min_{x\in [x_{j-1},x_j]}f(x)=\min\{f(x_{j-1}),f(x_j)\}$. Now, assume that $f$ is nonconstant on $[x_{j-1},x_j]$. If $a=0$, then $b\neq 0$. So, $f$ is strictly increasing if $b>0$ and strictly decreasing if $b<0$. This gives $\displaystyle\min_{x\in [x_{j-1},x_j]}f(x)=\min\{f(x_{j-1}),f(x_j)\}$. If $a\neq 0$, then $f'$ is strictly decreasing. Suppose, if possible, $f$ has a local minimum at some $x_{j-1}<x'<x_j$. This implies that $f'(x')=0$, $f'(x)<0$ for $x'-\delta<x<x'$, and $f'(x)>0$ for $x'<x<x'+\delta$ for some $\delta>0$, which is not possible as $f'$ is strictly decreasing. This gives $\displaystyle\min_{x\in [x_{j-1},x_j]}f(x)=\min\{f(x_{j-1}),f(x_j)\}$, and there is no other point in $[x_{j-1},x_j]$ at which $f$ attains this minimum value.\\\\
From the above discussion, we get $\displaystyle\min_{x\in [1,\infty)}f(x)=\min\{f(x_j):j=0,1,\ldots,m\}$. Thus, $f$ attains its global minimum at some $x_j\in X$. Also, if $f$ is nonconstant on each interval $[x_{j-1},x_j]$ for $j=1,2,\ldots,m$, then this global minimum is achieved only at some $x_j\in X$.

\end{document}